\documentclass[onefignum,onetabnum]{siamart171218}
\usepackage{amssymb}
\definecolor{Darkgreen}{RGB}{00,64,00}
\definecolor{Brightpink}{RGB}{255,00,255}
\usepackage{rotating}

\usepackage{lipsum}
\usepackage{amsfonts}
\usepackage{graphicx}
\usepackage{epstopdf}
\usepackage{algorithmic}
\usepackage{comment}
\ifpdf
  \DeclareGraphicsExtensions{.eps,.pdf,.png,.jpg}
\else
  \DeclareGraphicsExtensions{.eps}
\fi

\newsiamremark{remark}{Remark}
\newsiamremark{hypothesis}{Hypothesis}
\crefname{hypothesis}{Hypothesis}{Hypotheses}
\newsiamthm{claim}{Claim}

\headers{Generalized reduction for SWPM}{MJ Goeckner, D Harcey, RQ Pederson, A Niyonzima and J Zweck}

\title{A generalized reduction scheme for the Stochastic Weighted Particle Method\thanks{Submitted to the editors \today.}}

\author{Matthew Goeckner\thanks{Department of Physics, University of Texas at Dallas, Richardson TX 
  (\email{goeckner@utdallas.edu}})
\and Donovan Harcey\thanks{Departments of Physics and Mathematical Sciences, University of Texas at Dallas, Richardson, TX 
  (\email{donovan.harcey@utdallas.edu}})
\and Rainier Q Pederson\thanks{Departments of Computer Science and Mathematical Sciences, University of Texas at Dallas, Richardson TX 
  (\email{rainier.pederson@utdallas.edu}})
\and Axel Niyonzima\thanks{Departments of Physics and Computer Science, University of Texas at Dallas, Richardson TX 
  (\email{axel.niyonzima@utdallas.edu}})
\and John Zweck\thanks{Department of Mathematics, New York Institute of Technology, New York NY
  (\email{jzweck@nyit.edu}})}

\usepackage{amsopn}

\makeatletter
\newcommand*{\addFileDependency}[1]{
  \typeout{(#1)}
  \@addtofilelist{#1}
  \IfFileExists{#1}{}{\typeout{No file #1.}}
}
\makeatother

\newcommand*{\myexternaldocument}[1]{%
    \externaldocument{#1}%
    \addFileDependency{#1.tex}%
    \addFileDependency{#1.aux}%
}

\ifpdf
\hypersetup{
  pdftitle={A generalized reduction scheme for the Stochastic Weighted Particle Method},
  pdfauthor={Matthew Goeckner, Donovan Harcey, Rainier Q Pederson, Axel Niyonzima, John Zweck}
}
\fi

\myexternaldocument{ex_supplement}

\usepackage{verbatim}
\usepackage{hyperref}
\usepackage{nameref}

\begin{document}

\maketitle
\begin{abstract}
The Stochastic Weighted Particle Method (SWPM) of Rjasanow and Wagner is a generalization of the Direct Simulation Monte Carlo method for computing the probability density function of the velocities of a system of interacting particles for applications that include rarefied gas dynamics and plasma processing systems. Key components of a SWPM simulation are a particle grouping technique and particle reduction scheme. These are periodically applied to reduce the computational cost of simulations due to the gradual increase in the number of stochastic particles. A general framework for designing particle reduction schemes is introduced that enforces the preservation of a prescribed set of moments of the distribution through the construction and explicit solution of a system of linear equations for particle weights in terms of particle velocities and the  moments to be preserved. This framework is applied to preserve all moments of the distribution up to order three. Numerical simulations are performed to verify the scheme and quantify the degree to which even higher-order moments and tail functionals are preserved. These results reveal an unexpected trade off between the preservation of these higher-order moments and tail functionals. 
\end{abstract}

\section{Introduction}

The Direct Simulation Monte Carlo (DSMC) method, first proposed by Bird in 1963~\cite{Bird1963}, is a widely used  stochastic method for solving kinetic models of fluid and fluid-like systems. With DSMC, a system of interacting physical particles is modeled using stochastic particles, each of which represents a group of physical particles that are in close proximity in phase space. Collisions between stochastic particles are designed to approximate the collision processes between the physical particles they represent.  
The fraction of physical particles that are represented by a given stochastic particle is called the weight of the stochastic particle. In a DSMC simulation, every stochastic particle is assumed to have the same weight.  The simplicity of this approach has allowed researchers to use DSMC methods to model a wide variety of physical systems. For example, recent applications of DSMC include studies of gaseous diffusion as a correlated random walk~\cite{Yuste2024}, shear stress in Couette flow~\cite{Johnson2024}, physisorption of gas on dust~\cite{Yu2024}, and spacecraft reentry dynamics~\cite{Maigler2024}.  

In fluid and fluid-like systems in which very low probability events drive major physical or chemical processes, a standard DSMC simulation is not an efficient or sufficient method.  For example, the atmosphere of Venus is 96.5 mol\% CO$_2$, 3.5\% N$_2$, and approximately 200 parts per million (ppm) SO$_2$, 70 ppm Ar, and $3 \times 10^{-4}$ ppm Xe~\cite{Basilevsky_2003,Donahue1981}. Of particular interest is the Ar/Xe ratio and the ratios of the common isotopes of both species as those values can give insight into the Venus' origins, past geodynamics and atmospheric loss~\cite{Rabinovitch2023, Rabinovitch2019}. To collect samples and return them to earth, one would need to pass a spacecraft through Venus' upper atmosphere at hypersonic speeds. Such a collection technique is known to induce changes to gas and isotope ratios via differential diffusion of species through the bow shock~\cite{Rabinovitch2023,Rabinovitch2019}.  Exactly how much this differential diffusion affects the collected concentration of an isotope is difficult to ascertain with a standard DSMC simulation.  For example, making use of a series of modified atmospheric models, Rabinovitch et al.~\cite{Rabinovitch2019} show that the ratio of the average predicted concentration from 10 independent ensembles, relative to the standard deviation of those values, ranged from 25 when Xe isotope concentrations were assumed to be 1000~ppm to as low as 1.6 when Xe isotope concentrations were assumed to be 10~ppm,  both of which are well over the expected value on the order of $3 \times 10^{-4}$~ppm. These statistical errors were induced by low probability events. Despite using 1.2 billion stochastic particles in their simulations, these large uncertainties occurred because at the lowest concentrations, many ensembles produced only 0 or 1 stochastic Xe atoms. 

As with some traditional fluids, low probability events play a critical role in plasmas. Plasmas, particularly low-temperature plasmas (LTPs), are often not in thermal equilibrium.  For example in many LTPs, the ions and neutrals will often have effective temperatures close to room temperature, $\frac{1}{40}$~eV, while the electrons may have an effective temperature near 3~eV~\cite{Lieberman2ndEd}.  In addition, such systems often have velocity probability density functions (pdfs) that only approximate a true Maxwellian,
\begin{equation}
    f\left(\mathbf{v}\right) 
    = {\left( {\frac{m}{{2\pi kT}}} \right)^{3/2}}e^{-m(\mathbf{v}-\mathbf{v}_0)^2/2kT} \\
    = {\frac{1}{{(2\pi v^2_{\rm th}})^{3/2}}}e^{-(\mathbf{v}-\mathbf{v}_0)^2/2v^2_{\rm th}}.
\label{eqn:Maxwellian}\end{equation}
Here, $\mathbf{v}$ is the velocity, $\mathbf{v}_0$ is the drift velocity, $m$ is the mass, $k$ is Boltzmann's constant and $T$ is the temperature. The variance of the distribution is  the square of the thermal velocity, $v_{\text{th}}=(kT/m)^{1/2}$. In fact, experimental measurements of both electron and ion velocity distributions in LTPs, hot fusion plasmas, and even in solar wind have found distribution functions that have been characterized as being bi-Maxwellian~\cite{Hellinger2010, Kuznetsov2022} or Druyvesteyn distributions~\cite{Khrapak2010, Liao2017}, and that exhibit drifts~\cite{SHERIDAN1998} and have many other structures. Finally, a given plasma system can exhibit multiple different distributions, depending on various external control parameters~\cite{Chung2002,Donnelly2004}. The subtle differences between Maxwellian-like and truly Maxwellian distributions can play a major role in how such a plasma operates. For instance, it is often only the particles in the tail of the distribution that have enough energy to break chemical bonds and ionize neutral atoms.  Critically, such tail particles have velocities far from the drift velocity and may represent 1:10,000 or less of the total population. In a non-drifting laboratory plasma, it is often high-energy electrons ($E \sim 5\ kT$) which cause electron-ion pair production, determining the intensity and efficacy of the plasma. Here, the electrons may have an effective temperature of $kT\sim3\ eV$ (33000 K) while the ionization energy of Ar is $15.7\ eV$. This phenomenon is similar to hypersonic flows, where it is often the high-energy neutrals ($E \sim 6\ kT$) that drive chemical reactions.  For such flows, the gas may reach an effective temperature of $kT\sim0.5\ eV$ (5500 K), while the dissociation energy of $O_2$ is $5.2\ eV$.

Citing large density variations often observed in hypersonic systems, in 1996, Rjasanow and Wagner proposed the stochastic weighted particle method (SWPM) which is an adapation of the DSMC method~\cite{Rjasanow1996}.  SWPM overcomes the limitations found in DSMC by associating a variable weight, $w$, to each stochastic particle, thereby allowing different stochastic particles to represent different numbers of physical particles. By using variable weights it is possible to initialize the velocity distribution so as to obtain an approximately equal number of stochastic particles in each bin of the velocity histogram, which enables highly accurate representations of the velocity distribution function, especially in the low-probability tails. The SWPM algorithm is designed to preserve this feature of the histogram as physical time increases during the simulation. 
Because each stochastic particle represents a different number of physical particles, for collisions to be accurately modeled, the weights of the two colliding particles must be the same.  If the collision involves the $i$-th and $j$-th stochastic particles, then the $i$-th stochastic particle is split into two particles, one with weight $\gamma$ and the other with weight $w_i-\gamma$, and similarly for the $j$-th particle. The two new particles with weight $\gamma$ collide while the other two particles do not. Consequently, in each collision the number of stochastic particles increases by two. This increase in particle count causes the computation to slow down. To maintain computationally efficiency, several times over the course of a SWPM simulation, the stochastic particles are grouped and each group of particles is replaced by a reduced number of stochastic particles. Equally importantly, the grouping and reduction processes can be designed to keep the particles uniformly spread over phase space, thus maintaining the fraction of stochastic particles in the low probability tails, which decreases their statistical uncertainty.

The algorithms used for the particle grouping and reduction processes both impact the accuracy of the simulation~\cite{RWBook2005}. However in this article we focus on particle reduction schemes.  Rjasanow and Wagner proposed several such schemes~\cite{RWBook2005}, including methods that preserve (i) the zeroth order moment (the total weight), (ii) the zeroth and first order moments (total weight and drift velocity), and (iii) the total energy and the zeroth and first order moments. In ~\cite{Lama2020}, we proposed a reduction scheme that preserves the full pressure tensor (full second order moment) of the group, together with all of the lower order moments. The numerical results in ~\cite{Lama2020, RWBook2005} demonstrate that the accuracy of a SWPM simulation increases as more moments of the group are preserved. 

In this article, we propose a generalized particle reduction method for SWPM that explicitly preserves higher order moments. In principle, this new method can be used to preserve any desired collection of moments of a group, including elements of the third order tensor moment (and higher) and as such includes all of the previously studied reduction schemes as special cases \cite{Lama2020, RWBook2005}. We cast the problem of preserving a given collection of moments as a linear algebra problem of the following form. Given a vector, $\boldsymbol\mu$, whose entries correspond to the desired set of moments to be preserved, construct a matrix, $\mathbf P$, in terms of a set of velocities of particles in the reduced group so that the equation
\begin{equation}
      \mathbf P \mathbf w = \boldsymbol{\mu}
    \label{eqn:muPw}
\end{equation}
has a solution vector, $\mathbf w$, whose elements are the weights of the reduced particles. We call $\mathbf P$ the progenitor matrix, $\boldsymbol\mu$ the moment vector, and $\mathbf w$ the weight vector. Each scalar equation in this system enforces a condition that a particular moment of the reduced group is equal to the corresponding moment of the original group. The challenge is to choose the particle velocities used to define $\mathbf P$ in such a way that the components, $w_i$, of the weight vector $\mathbf w$ of \eqref{eqn:muPw} are all positive. In this paper, we use this approach to develop a  reduction scheme that preserves all moments of order three or less. In particular, we derive explicit analytical formulae for the particle velocities and weights in terms of a few parameters. This approach could also be applied to preserve desired higher-order moments or even tail functionals. 
 
We note that our methodology is not the only way in which the higher order moments can be examined/preserved.  While preparing this manuscript we learned that Oblapenko~\cite{Oblapenko2024} also formulated the moment-preservation problem as a linear system of the form~\eqref{eqn:muPw}, which he solved numerically using non-negative least squares algorithm. The advantages of our analytical approach is that we identify the minimal number of particles required in each reduced group and that we can tailor the choice of particle velocities to increase computational efficiency. Moreover, while Oblapenko shows results for the preservation of the 4-th total moment, we investigate the degree to which our reduction scheme preserves both 4-th and 5-th order moments and tail functionals. Additionally, although deterministic methods for kinetic models suffer from the curse of dimensionality, in some situations it is now possible to use full-order models that represent the solution in a well-chosen basis to achieve similar performance to standard stochastic methods~\cite{gamba2018galerkin}. More recently, attention has been focused on using adaptive tensor networks to compute low rank approximations of the velocity distribution function for a range of kinetic equations~\cite{einkemmer2024review,sands2024high}. While these methods have the potential for greatly increasing computational efficiency, they can be complicated to design and it is not yet clear whether the error inherent in a low rank representation of the distribution can be made small enough to access the high energy tails.

In section~\ref{sec:EnsureSolutions}, we show how to choose the particle velocities so that the progenitor matrix is guaranteed to be invertible.  In section~\ref{sec:StandardizedP} we show how this process is simplified in a standardized reference frame, and in section~\ref{sec:PosWeights}, we adapt this construction to ensure that the weights are nonnegative. In section~\ref{sec:Results} we determine how accurately the statistics of the velocity distribution can be computed prior to reduction using DSMC-like and SWPM-like simulations. Then we quantify any additional degree of uncertainty that is added by the reduction scheme. Finally, we show that for a particular grouping technique, we can preserve the pdf out to the $6$-{th} tail functional, provided we have a sufficient number of stochastic particles. 

\section{The general construction of the progenitor matrix}\label{sec:GeneralProgenitor}

\subsection[Ensuring existence of solutions]{Ensuring existence of solutions}\label{sec:EnsureSolutions}
Given an original group consisting of a large number of particles, here we  show how to create a reduced group of particles that has the same moments up to the $K$-{th} order as the original group. In \eqref{eqn:muPw}, the moment vector, $\boldsymbol{\mu}$, is determined from the original particles and hence is known.  In comparison, the weight vector, $\mathbf w$, describes the weights of the reduced system of particles and hence is unknown.  Careful construction of $\mathbf P$ allows us to map between these two vectors via the velocities of the reduced particles.

One of the simplest ways to solve \eqref{eqn:muPw} is if $\mathbf P$ is invertible. To construct an invertible progenitor matrix we first need to consider the general structure of $\mathbf P$.  The $j$-th column of $\mathbf P$ represents the contributions of the $j$-{th} reduced particle to the components of the moment vector, $\boldsymbol{\mu}$. The ordering of the moments in $\boldsymbol\mu$ determines the ordering of the rows of $\mathbf P$. To ensure that $\mathbf P$ is invertible we choose the number of particles in the reduced group to be equal to the total number, $N_K$, of moments up to the $K$-{th} order.    This yields a system of $N_K$ equations in $N_K$ unknowns. To further simplify this process, we choose the ordering and the velocities of the reduced particles so that $\mathbf P$ is an upper-triangular block matrix whose diagonal blocks are all invertible. Consequently, the entire matrix, $\mathbf P$, is invertible. 

Knowing how to choose the velocities of the reduced particles requires understanding how $\mathbf P$ maps $\mathbf w$ to $\boldsymbol{\mu}$.  The moments of a group of $N$ particles are defined in terms of the particle velocities, $\mathbf v_j = (v_{x,j},v_{y,j}, v_{z,j})$, and particle weights, $w_j$, for $j=1,\dots, N$. For each triple of non-negative integers, $k_{x},k_{y},k_{z}$, we associate  a component,
\begin{equation}
    M_{k_{x} k_{y} k_{z}} = \sum\limits_{j=1}^{N} w_j\, v_{x,j}^{k_x} \,v_{y,j}^{k_y}\, v_{z,j}^{k_z},
    \label{eqn:MomentOfTriple}
\end{equation}
of the $K$-th order  moment tensor, where  $K=k_x+k_y+k_z$. 
To resolve possible occurrences  of the indeterminate form $0^0$ in \eqref{eqn:MomentOfTriple}, we set $0^0=1$. If the $i$-{th} rows of $\boldsymbol\mu$ and $\mathbf P$ correspond to the moment $M_{k_{x} k_{y} k_{z}}$, and the $j$-{th} column of $\mathbf P$ corresponds to the $j$-{th} reduced particle, then the $(i,j)$-entry of $\mathbf P$ is given by
\begin{equation}
    P_{ij} = v_{x,j}^{k_x} \,v_{y,j}^{k_y}\,v_{z,j}^{k_z}.
    \label{eqn:Pijv}
\end{equation}
Therefore, the problem of determining the weights given the velocities and moments involves solving the linear system~\eqref{eqn:muPw}.

To illustrate how to choose the particle velocities to ensure that $\mathbf P$ is invertible, we first consider the simple case of moments of order at most $K=1$. (We say such a system is of order $K = 1$.) Since a system of order one in three velocity dimensions has $N_1^{(3)}$ = 4 moments, \eqref{eqn:muPw} is a system of four equations for the four weights, which in this case we denote by $w_0$, $w_x$, $w_y$, and $w_z$. Since the $0$-{th} moment,
\begin{equation}
   M_{000} =w_0 + w_x + w_y +w_z,
\end{equation}
is independent of the particle velocities, we choose the first reduced  particle to have zero velocity,  $\mathbf{v}_0=\left(0,0,0\right)^T$. Consequently, the first column of the progenitor matrix will have zeros below the (1,1)-entry. The first moment of the distribution is given by
\begin{equation}
        \begin{bmatrix} 
            M_{100} \\
            M_{010} \\
            M_{001}
        \end{bmatrix}
        \,\,=\,\, 
        w_0 \mathbf v_0 + w_x \mathbf v_x + w_y \mathbf v_y + w_z \mathbf v_z,
\end{equation}
where $\mathbf v_x$, $\mathbf v_y$, and $\mathbf v_z$ denote the velocities of the remaining three reduced particles. The simplest way to ensure that the progenitor matrix is invertible is to choose these three velocities to be of the form 
\begin{equation}
   \mathbf{v}_x = \begin{bmatrix} v_{x}&0&0 \end{bmatrix}^T, 
   \qquad
   \mathbf{v}_y = \begin{bmatrix} 0& v_{y} & 0 \end{bmatrix}^T,
   \qquad
   \mathbf{v}_z = \begin{bmatrix} 0& 0 & v_{z} \end{bmatrix}^T,
   \label{eqn:vK1}
\end{equation}
where $\mathbf v^T$ denotes the transpose of $\mathbf v$. In this case, the progenitor matrix is the upper triangular matrix,
\begin{equation}
\mathbf P \,\,=\,\,
\begin{bmatrix}
1&1&1&1 \\
0&v_{x}&0&0 \\
0&0&v_{y}&0 \\
0&0&0&v_{z} \\
\end{bmatrix},  
\label{eq:K1DiagP}
\end{equation}
which is invertible provided that the three velocities in \eqref{eqn:vK1} are chosen to be nonzero.  

To ensure that the weights are positive, we simply require that the
signs of the velocities match those of the moments in that
\begin{equation}
    \operatorname{sign}\left(v_x\right)=\operatorname{sign}\left(M_{100}\right), \quad
    \operatorname{sign}\left(v_y\right)=\operatorname{sign}\left(M_{010}\right), \quad
    \operatorname{sign}\left(v_y\right)=\operatorname{sign}\left(M_{001}\right).
\label{eqn:1D_sign_state}
\end{equation}
Finally, the velocities can be scaled so that 
\begin{equation}
    w_0 \,\,=\,\, M_{000} - \left( 
        \frac{v_{x}}{M_{100}} +
        \frac{v_{y}}{M_{010}} +
        \frac{v_{z}}{M_{001}}
    \right) \,\, > \,\, 0.
\end{equation}
In the event that a first moment is zero, that moment is omitted from the moment vector, $\boldsymbol\mu$, and the corresponding particle is removed from the group of reduced particles.  This solution is equivalent to Rjasanow and Wagner's~\cite{RWBook2005} method for preserving the drift velocity during reduction. 

Next, we consider the case that $K=2$. The full second moment is given by
\begin{equation}
 \begin{bmatrix} 
            M_{200} & M_{110} & M_{101}\\
            M_{110} & M_{020} & M_{011}\\
            M_{101} & M_{011} & M_{002}
        \end{bmatrix}
        \,\,=\,\,
        \sum\limits_{j} \,\, w_j
        \begin{bmatrix} 
            v_{x,j}^2 & v_{x,j} v_{y,j} & v_{x,j} v_{z,j}\\
            v_{x,j} v_{y,j} & v_{y,j}^2 & v_{y,j} v_{z,j}\\
            v_{x,j} v_{z,j} & v_{y,j} v_{z,j} & v_{z,j}^2
        \end{bmatrix}. 
\end{equation}
Since we have six independent second-order moments, we require an additional six particles, for a total of 10 particles.  We observe that three of these moments are functions of only one component of the velocity, while the other three are functions of two components. Following the pattern established for the first moment, we choose velocities of the additional six particles to be of the form
$\mathbf v_x$, $\mathbf v_y$, and $\mathbf v_z$, in \eqref{eqn:vK1}, and 
\begin{equation}
   \mathbf{v}_{xy} = \begin{bmatrix} v_{x}&v_{y} &0 \end{bmatrix}^T, 
   \qquad
   \mathbf{v}_{xz} = \begin{bmatrix} v_{x} & 0 & v_{z} \end{bmatrix}^T,
   \qquad
   \mathbf{v}_{yz} = \begin{bmatrix} 0&  v_{y}& v_{z} \end{bmatrix}^T.
   \label{eqn:vK2}
\end{equation}
With these choices  the progenitor matrix contains many zero entries. In fact, by suitably permuting the rows and columns of $\mathbf P$, we can ensure that $\mathbf P$ is block upper triangular. First, we arrange the components of the moment vector so that
\begin{equation}\label{eq:muK2}
    \boldsymbol{\mu}^T \,\,=\,\,
    \begin{bmatrix}
        \mu_0 &
        \boldsymbol{\mu}_x^T &
        \boldsymbol{\mu}_y^T &
        \boldsymbol{\mu}_z^T &
        {\mu}_{xy} &
        {\mu}_{xz} &
        {\mu}_{yz} 
    \end{bmatrix},
\end{equation}
where $\mu_0 = M_{000}$,
\begin{align}
   \boldsymbol{\mu}_x &= \begin{bmatrix}
        M_{100}\\M_{200}
    \end{bmatrix}, 
    \qquad 
    \boldsymbol{\mu}_y = \begin{bmatrix}
        M_{010}\\M_{020}
\end{bmatrix},
\qquad
\boldsymbol{\mu}_z = \begin{bmatrix}
        M_{001}\\M_{002}
    \end{bmatrix},
    \\ 
    {\mu}_{xy} &= \,\, M_{110},\,\,\,
  \qquad
{\mu}_{xz}=\,\,M_{101},\,\,\,\,
\qquad
{\mu}_{yz}=\,\,M_{011}.
\end{align}
In \eqref{eq:muK2}, the subscripts denote the components of the velocity vectors that contribute to the associated moments. Ordering the reduced particle velocities in the same manner we obtain a block upper triangular progenitor matrix. Specifically, we choose the reduced particle velocities to be given by the columns of the matrix
\begin{equation}
    \begin{bmatrix}
        0&v_{x,2}&v_{x,3}&0       &0      &0       &0       &v_{x,8}&v_{x,9}&0\\
        0&0      &0      &v_{y,4} &v_{y,5}&0       &0       &v_{y,8}&0      &v_{y,10}\\
        0&0      &0      &0       &0      &v_{z,6} &v_{z,7} &0      &v_{z,9}&v_{z,10}
    \end{bmatrix}.
\label{eqn:v_mat_2}
\end{equation}
With these choices, we can apply \eqref{eqn:Pijv} to conclude that \eqref{eqn:muPw} is of the form
\begin{equation}
    \begin{bmatrix}
        1&\mathbf{1}&\mathbf{1}&\mathbf{1}&{1}&{1}&1\\ 
        \mathbf{0}&\mathbf{P}_{x,x}&\mathbf{0}&\mathbf{0}&\mathbf{P}_{x,xy}&\mathbf{P}_{x,xz}&0\\ 
        \mathbf{0}&\mathbf{0}&\mathbf{P}_{y,y}&\mathbf{0}&\mathbf{P}_{y,xy}&0&\mathbf{P}_{y,yz}\\
        \mathbf{0}&\mathbf{0}&\mathbf{0}&\mathbf{P}_{z,z}&\mathbf{0}&\mathbf{P}_{z,xz}&\mathbf{P}_{z,yz}\\
        0&\mathbf{0}&\mathbf{0}&\mathbf{0}&P_{xy,xy}&0&0\\
        0&\mathbf{0}&\mathbf{0}&\mathbf{0}&0&P_{xz,xz}&0\\
        0&\mathbf{0}&\mathbf{0}&\mathbf{0}&0&0&P_{yz,yz}
    \end{bmatrix}
    \begin{bmatrix}
        w_0\\\mathbf{w}_x\\\mathbf{w}_y\\\mathbf{w}_z\\{w}_{xy}\\ {w}_{xz}\\ {w}_{yz}
    \end{bmatrix}
    =\begin{bmatrix}\mu_0\\\boldsymbol{\mu}_x\\\boldsymbol{\mu}_y\\\boldsymbol{\mu}_z\\{\mu}_{xy}\\{\mu}_{xz}\\{\mu}_{yz}\end{bmatrix},
    \label{eqn:Pwmu2}
\end{equation}
where $\mathbf{w}_x$, $\mathbf{w}_y$, and $\mathbf{w}_z$ are $2\times 1$. Here,
\begin{equation}
    \mathbf P_{x,x}\,\,=\,\,
    \begin{bmatrix}
        v_{x,2}&  v_{x,3} \\
        v_{x,2}^2&  v_{x,3}^2 
    \end{bmatrix},
\end{equation}
with $\mathbf P_{y,y}$ and $\mathbf P_{z,z}$ being defined similarly, and $P_{xy,xy} = v_{x,8}v_{y,8}$, with  $P_{xz,xz}$ and $P_{yz,yz}$ being defined similarly. If the entries, $v_{*,j}$  in \eqref{eqn:v_mat_2} are all chosen to be nonzero and such that 
\begin{equation}
    v_{x,3}\neq v_{x,2}, \qquad v_{y,5}\neq v_{y,4}, \qquad v_{z,7}\neq v_{z,6},
\label{eqn:VanDerMondeInvert2}
\end{equation}
then the diagonal blocks of the matrix $\mathbf P$ in \eqref{eqn:Pwmu2} are all invertible, which implies that $\mathbf P$ itself is invertible. This solution is effectively equivalent to 
the method of Lama \textit{et al.}~\cite{Lama2020} for preserving the pressure tensor during reduction.   

To preserve up to an arbitrary moment order, $K$, we order the components of the moment vector and, emulating the pattern in \eqref{eqn:v_mat_2}, choose the nonzero components of the  velocities of the reduced particles to ensure that the progenitor matrix is block upper triangular.  As before, we set the number of  new particles to
equal to the number of preserved moments, $N_K$.  
Applying a counting argument we find that $N_K$ increases rapidly with $K$.  Specifically, we find in two velocity dimensions, the total number of entries in the moment vector is
\begin{equation}
    N_K^{(2)} \,\,=\,\, 1 \,\,+\,\, 2K \,\,+\,\, 
    \frac{K(K-1)}2,
\end{equation}
and in three velocity dimensions it is given by,
\begin{equation}
    N_K^{(3)} \,\,=\,\, 1 \,\,+\,\, 3K \,\,+\,\, 3\, 
    \frac{K(K-1)}2\,\,+\,\,
    \frac{K(K-1)(K-2)}{6}.
\end{equation}
  Because $N_K$ increases rapidly with $K$, we expect that this method of creating $\mathbf P$ will find practical application for $K\leq4$.

As we will see below, the main difference between the $K\ge3$ and $K=2$ cases is that for $K\ge 3$, $\mathbf P$ has an additional block, $\mathbf P_{xyz,xyz}$, in which all velocity components of all particles  are non-zero. In general, the moment vector is of the form,
\begin{equation}
    \boldsymbol{\mu}^T =
    \begin{bmatrix}
        \mu_0 &
        \boldsymbol{\mu}_x^T &
        \boldsymbol{\mu}_y^T &
        \boldsymbol{\mu}_z^T &
        \boldsymbol{\mu}_{xy}^T &
        \boldsymbol{\mu}_{xz}^T &
        \boldsymbol{\mu}_{yz}^T &
        \boldsymbol{\mu}_{xyz}^T
    \end{bmatrix},
\label{eqn:muK}
\end{equation} 
where $\mu_0 = M_{000}$ and
\begin{equation}
   \boldsymbol{\mu}_x^T = 
   \begin{bmatrix}
        M_{100}& \dots & M_{K00}
    \end{bmatrix}, 
    \boldsymbol{\mu}_y^T = 
    \begin{bmatrix}
        M_{010}&  \dots & M_{0K0}
    \end{bmatrix},
    \boldsymbol{\mu}_z^T = 
    \begin{bmatrix}
        M_{001}& \dots & M_{00K}
    \end{bmatrix}.
\end{equation}
Furthermore, 
\begin{equation}
    \boldsymbol \mu_{xy} = 
    \begin{bmatrix}
        \mathbf a_1 & \mathbf a_2 & \cdots &
        \mathbf a_{K-1}
    \end{bmatrix}^T,
    \label{eqn:muxylexio}
\end{equation}
where 
\begin{equation}
    \mathbf a_k = 
    \begin{bmatrix}
        M_{k10} & M_{k20} & \cdots & M_{k(K-k)0}  
    \end{bmatrix},
\end{equation}
and $\boldsymbol \mu_{xz}$ and $\boldsymbol \mu_{yz}$ are defined similarly. We note that $\boldsymbol \mu_{xy}$ has $K(K-1)/2$ entries. For example, when $K=3$, 
   $\boldsymbol \mu_{xy}^T \,\,=\,\,
    \begin{bmatrix} 
        M_{110}&M_{120}&M_{210}
    \end{bmatrix}$, 
    $\boldsymbol \mu_{xz}^T \,\,=\,\,
    \begin{bmatrix} 
        M_{101}&M_{102}&M_{201}
    \end{bmatrix}$, and 
    $\boldsymbol \mu_{yz}^T \,\,=\,\,
    \begin{bmatrix} 
        M_{011}&M_{012}&M_{021}
    \end{bmatrix}$.
Finally, the vector $\boldsymbol \mu_{xyz}$ is constructed by a lexiographic ordering of all mixed moments $M_{k_x k_y k_z}$, where $k_x\neq 0$, $k_y\neq 0$, $k_z\neq 0$, and $3\leq k_x+k_y+k_z \leq K$. A counting argument shows that $\boldsymbol \mu_{xyz}$ has $K(K-1)(K-2)/6$ entries. When $K=3$, $\boldsymbol \mu_{xyz}=M_{111}$, and when $K=4$, 
$\boldsymbol \mu_{xyz}^T \,\,=\,\,
\begin{bmatrix} 
M_{111} & M_{112} & M_{121} &M_{211}
\end{bmatrix}$.
By appropriately setting selected components of the reduced particle velocities to zero we obtain a block upper triangular system of the form
\begin{equation}
\begin{bmatrix}
    1&\mathbf{1}&\mathbf{1}&\mathbf{1}&\mathbf{1}&\mathbf{1}&\mathbf{1}&\mathbf{1}\\
    \mathbf{0}&\mathbf{P}_{x,x}&\mathbf{0}&\mathbf{0}&\mathbf{P}_{x,xy}&\mathbf{P}_{x,xz}&\mathbf{0}& \mathbf{P}_{x,xyz}\\
    \mathbf{0}&\mathbf{0}&\mathbf{P}_{y,y}&\mathbf{0}&\mathbf{P}_{y,xy}&\mathbf{0}&\mathbf{P}_{y,yz} &\mathbf{P}_{y,xyz}\\
    \mathbf{0}&\mathbf{0}&\mathbf{0}&\mathbf{P}_{z,z}&\mathbf{0}&\mathbf{P}_{z,xz}&\mathbf{P}_{z,yz} &\mathbf{P}_{z,xyz}\\
    \mathbf{0}&\mathbf{0}&\mathbf{0}&\mathbf{0}&\mathbf{P}_{xy,xy}&\mathbf{0}&\mathbf{0}&\mathbf{P}_{xy,xyz}\\
    \mathbf{0}&\mathbf{0}&\mathbf{0}&\mathbf{0}&\mathbf{0}&\mathbf{P}_{xz,xz}&\mathbf{0}&\mathbf{P}_{xz,xyz}\\
    \mathbf{0}&\mathbf{0}&\mathbf{0}&\mathbf{0}&\mathbf{0}&\mathbf{0}&\mathbf{P}_{yz,yz}&\mathbf{P}_{yz,xyz} \\
    \mathbf{0}&\mathbf{0}&\mathbf{0}&\mathbf{0}&\mathbf{0}&\mathbf{0}&\mathbf{0}&\mathbf{P}_{xyz,xyz}
\end{bmatrix}
\begin{bmatrix}
    w_0\\\mathbf{w}_x\\\mathbf{w}_y\\\mathbf{w}_z\\\mathbf{w}_{xy}\\ \mathbf{w}_{xz}\\ \mathbf{w}_{yz} \\ \mathbf{w}_{xyz}
\end{bmatrix}
    =
\begin{bmatrix}
    \mu_0\\\boldsymbol{\mu}_x\\\boldsymbol{\mu}_y\\\boldsymbol{\mu}_z\\
    \boldsymbol{\mu}_{xy}\\
    \boldsymbol{\mu}_{xz}\\
    \boldsymbol{\mu}_{yz}\\
    \boldsymbol{\mu}_{xyz}
    \end{bmatrix}.
    \label{eqn:PwmuK}
\end{equation}
Specifically, to obtain this form, the particle velocities used to generate the columns containing the block $\mathbf P_{x,x}$ must be chosen so that $v_y=v_z=0$, and similarly for $\mathbf P_{y,y}$ and $\mathbf P_{z,z}$. Those used to generate the columns containing $\mathbf P_{xy,xy}$ must be chosen so that $v_z=0$, and similarly for $\mathbf P_{xz,xz}$ and $\mathbf P_{yz,yz}$. Finally the particle velocities used to generate the columns containing the block $\mathbf P_{xyz,xyz}$ must be chosen to have no zero components.  

The progenitor matrix $\mathbf P$ in \eqref{eqn:PwmuK} is invertible if each block on the diagonal is invertible. 
By placing additional constraints on the velocities of the reduced particles, we can ensure that each of the four types of diagonal blocks in \eqref{eqn:PwmuK} is invertible. 
\vskip5pt
\noindent
{\bf Type 1:} The first block is a $1\times 1$ identity matrix which is invertible.
\vskip5pt
\noindent
{\bf Type 2:} The next three blocks $\mathbf P_{x,x}$, $\mathbf P_{y,y}$, and $\mathbf P_{z,z}$ are all $K\times K$ matrices of the form 
\begin{equation}
    \begin{bmatrix}
        v_1 & v_2 & \cdots & v_K \\
        v_1^2 & v_2^2 & \cdots & v_K^2 \\
        \vdots & \vdots & \cdots & \vdots\\
        v_1^K & v_2^K & \cdots & v_K^K \\
    \end{bmatrix}
    \,\,=\,\,
    \begin{bmatrix}
        1 & 1 & \cdots & 1 \\
        v_1 & v_2 & \cdots & v_K \\
        \vdots & \vdots & \cdots & \vdots\\
        v_1^{K-1} & v_2^{K-1} & \cdots & v_K^{K-1} \\
    \end{bmatrix}
    \,
    \begin{bmatrix}
        v_1 & 0 & \cdots & 0\\
        0 & v_2 & \cdots & 0\\
        \vdots & \vdots & \cdots & \vdots\\
        0 & 0 & \cdots & v_K \\
    \end{bmatrix}.
\end{equation}
The first factor on the right hand side is a Vandermonde matrix, $\mathbf V$, whose determinant is given by
\begin{equation}
    \operatorname{det} \mathbf V \,\,=\,\,
    \prod\limits_{0\leq k < l \leq K} (v_k-v_l),
\end{equation}
which is nonzero provided that $v_k\neq v_l$ for all  $k\neq l$. Therefore, these  blocks are invertible provided that $v_k\neq 0$ for all $k$ and $v_k\neq v_l$ for all $k\neq l$. 
\vskip5pt
\noindent
{\bf Type 3:} The three blocks, $\mathbf P_{xy,xy}$, $\mathbf P_{xz,xz}$, and $\mathbf P_{yz,yz}$ all have the same structure. Here we focus on the  block, $\mathbf P_{xy,xy}$, which is a $K(K-1)/2\times K(K-1)/2$ matrix. We order each column of $\mathbf P_{xy,xy}$ using the same lexicographic ordering we used for $\boldsymbol{\mu}_{xy}$ in \eqref{eqn:muxylexio}. For example, if $K=2$, we obtain a scalar of the form $\mathbf P_{xy,xy}=v_xv_y$, which is invertible provided $v_x$ and $v_y$ are both nonzero. 
\vskip5pt
\noindent
{\bf Type 4:} 
In the case $K=3$, the block, $\mathbf P_{xyz,xyz}=v_{x}v_{y}v_{z}$, is a scalar, which is invertible provided $v_x$, $v_y$, and $v_z$ are all nonzero.

\subsection{Standardizing $\bf{P}$}\label{sec:StandardizedP}
While the discussion above describes how to create $\bf{P}$ in a general reference frame, it is straightforward to standardize the velocity pdf, which makes the solution easier to obtain. To do so, we follow the typical method of standardizing a pdf in three dimensions.  First, we normalize the distribution so that $M_{000}=1$.  Second, we translate velocity coordinates so that $M_{100}=M_{010}=M_{001}=0$.  Third, we rotate the coordinate system so that the off-diagonal elements of the covariance matrix are zero, giving $M_{110}=M_{101}=M_{011}=0$. Finally, we standardize the distribution along each of the major principal component axes, ensuring that $M_{200}=M_{020}=M_{002}=1$. Then, using this standardized reference frame, we choose the velocities of the reduced particles, construct $\bf{P}$, and solve \eqref{eqn:muPw} for $\bf{w}$.  Finally, we transform the velocities of the reduced particles back into the original reference frame and unnormalize the weights. Since the resulting equations for the weights are simpler, in the next section we can make use of this process to choose the velocities so that the weights are positive.  

\subsection{Ensuring weight positivity} \label{sec:PosWeights} 
Our generalized reduction scheme, \eqref{eqn:muPw}, does not a priori ensure that the weights of the reduced particles are positive, $w_i > 0$, as is required for the simulation to be consistent with physical reality. In the case $K=1$ we showed that the weights are positive provided that we impose the additional condition \eqref{eqn:1D_sign_state} that the signs of the velocity components agree with the signs of the original moments.  When $K=1$, this sign state matching condition is simple to formulate and satisfy since the equations for the weights $w_x$, $w_y$, and $w_z$ are uncoupled due to the way we chose the particle velocities to form  $\mathbf P$ in \eqref{eq:K1DiagP}. In subsection~\ref{sec:PosWeightsK2}, in the case $K=2$ we will show that by standardizing $\mathbf P$ the system can also mostly be decoupled, which again readily yields a solution with positive weights.

However, when $K\geq 3$ nontrivial couplings between the weights cannot be avoided and it is no longer trivial how to choose the signs of the particle velocities to ensure weight  positivity. For example, for a particle with velocity of the form $\mathbf v = (v_x, v_y,0)$ it is not clear in which quadrant of $\mathbb R^2$ we should place the vector $(v_x, v_y)$, and for a particle with three non-zero velocity components it is not clear in which octant to place the particle velocity. In such a situation a natural way to proceed is to expand the number of reduced particles and to symmetrically place the particle velocities in multiple quadrants/octants to guarantee that the resulting system of moment equations has a solution with positive weights. Increasing the number of particles results in an underdetermined system which has a subspace of solutions of dimension at least one. Geometrically, the condition that the weights are all positive means that the weight vector lies in the first orthant of the solution space. The goal then is to select the sign states of the particle velocities so as to ensure that this solution subspace intersects the first orthant in the space of weight vectors. In section~\ref{sec:PosWeightsK3}, we show that this can be done so as to obtain relatively simple explicit formulae for both the particle velocities and the weights in terms of the moments of the original system. By appropriately choosing the parameters in this construction we can obtain better statistical estimates of higher-order moments ($K\geq 4$) and tail functionals even though our reduction algorithm does not explicitly enforce conservation of these statistics. 

\subsubsection{Case $K=2$ }\label{sec:PosWeightsK2}
Here we show how to choose the particles velocities in the standardized velocity space so that the weights are positive.  The idea is to choose the velocities so that it is easy to solve the system using back substitution while ensuring that the weights we are solving for are positive.  Henceforth, if a particle has weight zero we omit it from the reduced group, which has the effect of reducing the size of the linear system. 

For ease of exposition, we first consider the problem in $\mathbb{R}^2$, in which we have first normalized, centralized, rotated, and standardized the distribution as described in section~\ref{sec:StandardizedP}.  Under these conditions, the system reduces to the form:
\begin{equation}\label{eqn:k=2,2d,prog}
    \begin{bmatrix}
        1 & 1 & 1 & 1 & 1 & 1\\
        0 & v_{x,1} & v_{x,2} & 0 & 0 & v_{x,5} \\
        0 & v_{x,1}^2 & v_{x,2}^2 & 0 & 0 & v_{x,5}^2 \\
        0 & 0 & 0 & v_{y,3} & v_{y,4} & v_{y,5} \\
        0 & 0 & 0 & v_{y,3}^2 & v_{y,4}^2 & v_{y,5}^2 \\
        0 & 0 & 0 & 0 & 0 & v_{x,5}v_{y,5}\\
    \end{bmatrix}
    \begin{bmatrix}
        w_0 \\
        w_1 \\
        w_2 \\
        w_3 \\
        w_4 \\
        w_5 
    \end{bmatrix}
    =\begin{bmatrix}
        M_{000} \\
        M_{100} \\
        M_{200} \\
        M_{010} \\
        M_{020} \\
        M_{110} 
    \end{bmatrix}
    =\begin{bmatrix}
        1 \\
        0 \\
        1 \\
        0 \\
        1 \\
        0 
    \end{bmatrix}.
\end{equation}
Clearly, we can choose $w_5=0$, which means the last column and row can be omitted from $\mathbf P$. To further simplify the system, we introduce a speed parameter, $s>0$, and set $v_{x,1} = v_{y,3}=s$ and $v_{x,2} = v_{y,4}=-s$, thereby ensuring the invertibility of the diagonal blocks. Solving for the weight vector we obtain 
$\mathbf w=\begin{bmatrix} 
    1 - \frac{2}{s^2} &
    \frac{1}{2s^2} & \frac{1}{2s^2} &
    \frac{1}{2s^2} & \frac{1}{2s^2}\end{bmatrix}^T$,
which has positive components, provided $s> \sqrt{2}$.

\vspace{.2cm}
Similarly,  in $\mathbb{R}^3$  the moment vector in \eqref{eqn:Pwmu2} becomes $\mu_0=1$, $\boldsymbol\mu_x = \boldsymbol\mu_y =\boldsymbol\mu_z = \begin{bmatrix}0&1\end{bmatrix}^T$, and $\mu_{xy} = \mu_{xz} = \mu_{yz} = 0$. Just as in $\mathbb R^2$, we can omit the last three particles in \eqref{eqn:v_mat_2}, and set $v_{x,2} = v_{y,4} = v_{z,6} = s$ and $v_{x,3} = v_{y,5} = v_{z,7} = -s$. This gives the solution 
\begin{equation}\label{weightsolnK23D}
    \mathbf w=\begin{bmatrix}
        1 - \frac{3}{s^2} &
        \frac{1}{2s^2} & \frac{1}{2s^2} &
        \frac{1}{2s^2} & \frac{1}{2s^2} &
        \frac{1}{2s^2} & \frac{1}{2s^2}
    \end{bmatrix}^T,
\end{equation}
which has positive components, provided $s\ge \sqrt{3}$. Note that when $s=\sqrt{3}$ this solution is identical to the one found by Lama, et al.~\cite{Lama2020}. 

\subsubsection{Case $K=3$}
\label{sec:PosWeightsK3}
In the case $K=3$, solving the system $\mathbf P\mathbf w=\boldsymbol{\mu}$ for positive weights via back substitution becomes difficult if we use the minimum possible number of reduced particles, i.e. if $\mathbf P$ is a  square matrix. However, by using some additional symmetrically-located reduced particles we can readily ensure that the resulting underdetermined system $\mathbf P\mathbf w=\boldsymbol\mu$ has a solution with positive weights. 

The fundamental idea that drives the symmetrical placement of the reduced particle velocities is that for each preserved moment  there should be  at least one reduced particle whose velocity has a sign state that is compatible with that of the moment. We say the sign state of a  particle velocity $(v_x,v_y,v_z)$ is compatible with the sign state of the moment $M_{n_x,n_y,n_z}$ if
\begin{equation}
    \text{sign}(v_x^{n_x}v_y^{n_y}v_z^{n_z}) = \text{sign}(M_{n_x,n_y,n_z}).
    \label{signstate}
\end{equation}
The total number of compatible sign states of a given moment is  the number of orthants in velocity space containing particles that satisfy \eqref{signstate}. 

To illustrate how to place particles we consider the simplified system $v_xv_yw=M_{110}$. As  in \eqref{eqn:1D_sign_state}, $w>0$ provided that the particle velocity satisfies
\begin{equation}
     \text{sign}(v_xv_y) = \text{sign}(M_{110}). 
     \label{eq:signM110}
\end{equation}
When $\text{sign}(M_{110})>0$ the compatible sign states correspond to velocities in the first and third quadrants, and when $\text{sign}(M_{110})<0$, they are in the second and fourth quadrants.

In the case $K=3$ in  $\mathbb{R}^3$, the $20\times 20$ system  $\mathbf P\mathbf w=\boldsymbol{\mu}$  is  
given in block form by \eqref{eqn:PwmuK}
which we will expand to an underdetermined $20\times 32$ system, $\widetilde{\mathbf P} \widetilde{\mathbf w} = \boldsymbol \mu$, by splitting some of the reduced particles into twins or quadruplets, as needed. The expanded system and its solution can be obtained using block back substitution. First, we split the particle whose velocities define the last  column of $\mathbf P$ into a quadruplet, thereby expanding the subsystem $P_{xyz,xyz} w_{xyz} = \mu_{xyz}$, to a $1\times 4$ system $\widetilde{\mathbf P}_{xyz,xyz} \widetilde{\mathbf w}_{xyz} = \mu_{xyz}$. We choose the velocities in this quadruplet to guarantee that the weights of these four particles are positive. In the next back substitution  step,  we expand and solve the block $\mathbf P_{yz,yz} \mathbf w_{yz} = \boldsymbol \mu_{yz} - \widetilde{\mathbf P}_{xyz,yz} \widetilde{\mathbf w}_{xyz}$, where the modified moment vector on the right hand side is known from the previous step. Continuing this process, we iteratively expand the system and solve for the weights all the way back to $w_0$.  As we will show below for pure moments, such as those in $\boldsymbol \mu_x$ which only involve one velocity component, an analytical solution with positive weights can be obtained without needing to expand the  subsystem involving $\mathbf P_x$.

The method for splitting particles that contribute to the mixed moments of the distribution can be outlined as follows. Each particle with two nonzero velocity components is split into a set of two stochastic particles (twins) of identical weight. Similarly, each particle with three nonzero velocity components is split into a set of four stochastic particles (quadruplets). We associate the reduced particle that determines the $k$-th column of $\mathbf P$ to the $k$-th entry of the moment vector. The twins/quadruplets that are generated from the $k$-th reduced particle are assigned to all the quadrants/octants of velocity space that have sign states which are compatible with the sign of the (modified) $k$-th entry of the moment vector. The resulting underdetermined system, $\mathbf{\widetilde{P}}\mathbf{\widetilde{w}} =\boldsymbol{\mu}$, is of the form
\begin{equation}
  \begin{bmatrix}
        1 & \mathbf{1} & \mathbf{1} & \mathbf{1} & \mathbf{1} & \mathbf{1} & \mathbf{1} & \mathbf{1} \\
        \mathbf{0} & \mathbf{P}_{x} & \mathbf{0}  & \mathbf{0} & \mathbf{\widetilde{P}}_{xy,x} & \mathbf{\widetilde{P}}_{xz,x} &  \mathbf{0} & \mathbf{\widetilde{P}}_{xyz,x}  \\
        \mathbf{0} & \mathbf{0} & \mathbf{P}_y  & \mathbf{0} & \mathbf{\widetilde{P}}_{xy,y} &  \mathbf{0}  & \mathbf{\widetilde{P}}_{yz,y} & \mathbf{\widetilde{P}}_{xyz,y} \\
        \mathbf{0} &  \mathbf{0} & \mathbf{0} & \mathbf{P}_z & \mathbf{0}  & \mathbf{\widetilde{P}}_{xz,z} & \mathbf{\widetilde{P}}_{yz,z} & \mathbf{\widetilde{P}}_{xyz,z}  \\
        \mathbf{0} & \mathbf{0} & \mathbf{0}  & \mathbf{0} & \mathbf{\widetilde{P}}_{xy,xy} & \mathbf{0} & \mathbf{0} & \mathbf{\widetilde{P}}_{xyz,xy}   \\
          \mathbf{0} & \mathbf{0} & \mathbf{0} & \mathbf{0}  & \mathbf{0} & \mathbf{\widetilde{P}}_{xz,xz} & \mathbf{0} & \mathbf{\widetilde{P}}_{xyz,xz}      \\
            \mathbf{0} & \mathbf{0} & \mathbf{0} & \mathbf{0} & \mathbf{0}  & \mathbf{0} & \mathbf{\widetilde{P}}_{yz,yz}   & \mathbf{\widetilde{P}}_{xyz,yz}\\
                 \mathbf{0}  & \mathbf{0}  & \mathbf{0}    & \mathbf{0} & \mathbf{0}  & \mathbf{0} & \mathbf{0} & \mathbf{\widetilde{P}}_{xyz,xyz}  
    \end{bmatrix}\begin{bmatrix}
        w_0\\\mathbf{w}_x\\\mathbf{w}_y\\\mathbf{w}_z\\\mathbf{\widetilde{w}}_{xy}\\\mathbf{\widetilde{w}}_{xz}\\\mathbf{\widetilde{w}}_{yz}\\\mathbf{\widetilde{w}}_{xyz}
    \end{bmatrix}=\begin{bmatrix}\mu_0\\\boldsymbol{\mu}_x\\\boldsymbol{\mu}_y\\\boldsymbol{\mu}_z\\\boldsymbol{\mu}_{xy}\\\boldsymbol{\mu}_{xz}\\\boldsymbol{\mu}_{yz}\\\mu_{xyz}
    \end{bmatrix}.
    \label{k3sysexp}
\end{equation}

Using back substitution, we first consider the bottom right equation in \eqref{eqn:PwmuK} given by $P_{xyz,xyz} w_{xyz}=\mu_{xyz}$. With particle velocity indexing starting at zero, it is the 19-th particle that is used to define $P_{xyz,xyz}$. The compatibility condition between this particle velocity and the moment $\mu_{xyz}$ is that 
\begin{equation}
    \text{sign}(v_{x19}v_{y19}v_{z19}) \,\,=\,\, \text{sign}(\mu_{xyz}) \,\,=\,\, \text{sign}(M_{111}). 
    \label{quadrupsignstate}
\end{equation}
If $M_{111} > 0$ the compatible signs states are 
 $(+,+,+)$,  $(+,-,-)$, $(-,+,-)$, and $(-,-,+)$, while
if $M_{111} < 0$ they are  $(-,-,-)$,  $(-,+,+)$, $(+,-,+)$, and $(+,+,-)$.
To include particle velocities with all four of these sign states, we expand the subsystem to $\mathbf{\widetilde{P}}_{xyz,xyz}\mathbf{\widetilde{w}}_{xyz} = \mu_{xyz}$ where 
\begin{equation}
    \mathbf{\widetilde{P}}_{xyz,xyz} = \begin{bmatrix} v_{x,19}^{(a)}v_{y,19}^{(a)}v_{z,19}^{(a)}&
    v_{x,19}^{(b)}v_{y,19}^{(b)}v_{z,19}^{(b)}&
    v_{x,19}^{(c)}v_{y,19}^{(c)}v_{z,19}^{(c)}&
    v_{x,19}^{(d)}v_{y,19}^{(d)}v_{z,19}^{(d)}
    \end{bmatrix},
\end{equation}
and the weight of each quadruplet particle is a quarter of the weight of the original particle, 
\begin{equation}
    \mathbf{\widetilde{w}}_{xyz}=\begin{bmatrix}
        \frac{1}{4}w_{19} & \frac{1}{4}w_{19} & \frac{1}{4}w_{19} & \frac{1}{4}w_{19}
    \end{bmatrix}^T.
    \label{eqn:w_xyz_Vector}
\end{equation}
To satisfy the compatibility condition \eqref{quadrupsignstate} we choose the particle velocities to be of the form
\begin{equation}
    \begin{bmatrix}
        v_{x,19}^{(a)}& v_{x,19}^{(b)} &  v_{x,19}^{(c)} &  v_{x,19}^{(d)} \\
        v_{y,19}^{(a)}& v_{y,19}^{(b)} &  v_{y,19}^{(c)} &  v_{y,19}^{(d)} \\
        v_{z,19}^{(a)}& v_{z,19}^{(b)} &  v_{z,19}^{(c)} &  v_{z,19}^{(d)} 
 \end{bmatrix}
  \,\, =\,\, 
   s_{\rm{quad}} 
\begin{bmatrix}
    \alpha_{111} & -\alpha_{111} & - \alpha_{111} & \alpha_{111}\\
    1 & -1 & 1 & -1 \\
    1 & 1 &-1 &-1
\end{bmatrix},
    \label{quadvels}
\end{equation}
where $\alpha_{111} = \frac{M_{111}}{|M_{111}|}$ and $s_{\rm{quad}}$ is a speed parameter to be determined. Solving the subsystem using the velocity vectors given by \eqref{quadvels} we find that
\begin{equation}
    w_{19} = \frac{|M_{111}|}{s_{\rm quad}^3},
\end{equation}
which implies that all elements of $\mathbf{\widetilde{w}}_{xyz}$ are nonnegative.

The subsystem involving the 5-th diagonal block of \eqref{k3sysexp}  is of the form
\begin{equation}
    \mathbf{\widetilde{P}}_{xy,xy}\mathbf{\widetilde{w}}_{xy}  = \boldsymbol{\mu}_{xy} - \mathbf{\widetilde{P}}_{xyz,xy}\mathbf{\widetilde{w}}_{xyz}
    = \boldsymbol{\mu}_{xy},
    \label{xyeq}
\end{equation}
where the final equation holds since $\mathbf{\widetilde{P}}_{xyz,xy}\mathbf{\widetilde{w}}_{xyz}=0$. The block $\mathbf{\widetilde{P}}_{xy,xy}$ is obtained from $\mathbf{{P}}_{xy,xy}$ by splitting each of the three original particles into twins, yielding
\begin{equation}
    \mathbf{\widetilde{P}}_{xy,xy} = \begin{bmatrix}
        v_{x,7}v_{y,7} & \tilde{v}_{x,7}\tilde{v}_{y,7} &  v_{x,8}v_{y,8} & \tilde{v}_{x,8}\tilde{v}_{y,8} &  v_{x,9}v_{y,9} & \tilde{v}_{x,9}\tilde{v}_{y,9}\\
         v_{x,7}v_{y,7}^2 & \tilde{v}_{x,7}\tilde{v}_{y,7}^2 &  v_{x,8}v_{y,8}^2 & \tilde{v}_{x,8}\tilde{v}_{y,8}^2 &  v_{x,9}v_{y,9}^2 & \tilde{v}_{x,9}\tilde{v}_{y,9}^2\\
          v_{x,7}^2v_{y,7} & \tilde{v}_{x,7}^2\tilde{v}_{y,7} &  v_{x,8}^2v_{y,8} & \tilde{v}_{x,8}^2\tilde{v}_{y,8} &  v_{x,9}^2v_{y,9} & \tilde{v}_{x,9}^2\tilde{v}_{y,9}\\
    \end{bmatrix}.
\end{equation}
The corresponding block of the weight vector is given by
\begin{equation}
    \mathbf{\widetilde{w}}_{xy} = \begin{bmatrix}
        \frac{1}{2}w_7&\frac{1}{2}w_7&       \frac{1}{2}w_8&\frac{1}{2}w_8&       \frac{1}{2}w_9&\frac{1}{2}w_9
    \end{bmatrix}^T,
\end{equation}
and we recall that
\begin{equation}
    \boldsymbol{\mu}_{xy} = \begin{bmatrix}
        M_{110}&    M_{120}&    M_{210}
    \end{bmatrix}^T.
\end{equation}
Each of the moments in $\boldsymbol{\mu}_{xy}$ gives a compatibility condition. The compatibility condition for $M_{110}$ is given by \eqref{eq:signM110} with $(v_x,v_y) = v_{x,7}, v_{y,7}$, while the one for  $M_{120}$ is given by
\begin{equation}
    \text{sign}(M_{120}) = \text{sign}(v_{x,8}v_{y,8}^2) = \text{sign}(v_{x,8}).
    \label{m120signstate}
\end{equation}
So, for example,  if $M_{120} > 0$, then the particle is placed in the first or fourth quadrants.
Similarly,  $ \text{sign}(M_{210}) = \text{sign}(v_{y,9})$, which means that if $M_{210} > 0$, then the particle is placed in the first or second quadrants.
To satisfy these compatibility conditions we choose the particle velocities to be of the form
\begin{equation}
\begin{split}
    \begin{bmatrix}
        v_{x,7} & \tilde{v}_{x,7} \\
        v_{y,7} & \tilde{v}_{y,7} 
    \end{bmatrix}
    &= 
    s_{\rm twin} \begin{bmatrix}
       1 & -1 \\
       1 & -1 
 \end{bmatrix},
    \\
    \begin{bmatrix}
        v_{x,8} & \tilde{v}_{x,8} \\
        v_{y,8} & \tilde{v}_{y,8} 
    \end{bmatrix}
    &= 
    s_{\rm twin} \begin{bmatrix}
       \alpha_{120}  & \alpha_{120} \\
       1 & -1 
 \end{bmatrix},
\\
    \begin{bmatrix}
        v_{x,9} & \tilde{v}_{x,9} \\
        v_{y,9} & \tilde{v}_{y,9} 
    \end{bmatrix}
    &= 
    s_{\rm twin} \begin{bmatrix}
        1 & -1 \\
       \alpha_{210}  & \alpha_{210} \\
 \end{bmatrix}.
    \label{xyvals}
\end{split}
\end{equation}
where $\alpha_{120} = M_{120}/|M_{120}|$,  $\alpha_{210} = M_{210}/|M_{210}|$, and $s_{\rm{twin}}$ is a speed parameter to be determined. Solving \eqref{xyeq} for $\mathbf{\widetilde{w}}_{xy}$ using \eqref{xyvals} we have
\begin{equation}
    \mathbf{\widetilde{w}}_{xy} = \begin{bmatrix}
        \frac{|M_{110}|}{2s_{\rm twin}^2} & 
     
        \frac{|M_{110}|}{2s^2_{\rm twin}} &

          \frac{|M_{120}|}{2s^3_{\rm twin}} &
     
        \frac{|M_{120}|}{2s^3_{\rm twin}} &
        
         \frac{|M_{210}|}{2s^3_{\rm twin}} &
     
        \frac{|M_{210}|}{2s^3_{\rm twin}} &
        \end{bmatrix}^T,
\end{equation}
whose elements are all clearly nonnegative.

Similarly, the solutions involving the $6$-th and $7$-th blocks of \eqref{k3sysexp} are given by,
\begin{align}
    \mathbf{\widetilde{w}}_{xz} &= 
    \begin{bmatrix}
        \frac{|M_{101}|}{2s_{\rm twin}^2} &
        \frac{|M_{101}|}{2s^2_{\rm twin}} &
          \frac{|M_{102}|}{2s^3_{\rm twin}} &
        \frac{|M_{102}|}{2s^3_{\rm twin}} &     
         \frac{|M_{201}|}{2s^3_{\rm twin}} &  
        \frac{|M_{201}|}{2s^3_{\rm twin}}
        \end{bmatrix}^T,
        \\ 
        \mathbf{\widetilde{w}}_{yz} &= \begin{bmatrix}
        \frac{|M_{011}|}{2s^2_{\rm twin}} &  
        \frac{|M_{011}|}{2s^2_{\rm twin}} &
          \frac{|M_{012}|}{2s^3_{\rm twin}} &
        \frac{|M_{012}|}{2s^3_{\rm twin}} & 
         \frac{|M_{021}|}{2s^3_{\rm twin}} &
        \frac{|M_{021}|}{2s^3_{\rm twin}}
        \end{bmatrix}^T.
\end{align}
and the velocity vectors similarly follow the patterns given in \eqref{xyvals}.

Next we consider the subsystem involving the  second block of \eqref{k3sysexp},
    \begin{equation}
    \mathbf{P}_x\mathbf{w}_x  = \boldsymbol{\mu}_{x}- \mathbf{\widetilde{P}}_{xy,x}\mathbf{\widetilde{w}}_{xy} -\mathbf{\widetilde{P}}_{xz,x}\mathbf{\widetilde{w}}_{xz} - \mathbf{\widetilde{P}}_{xyz,x}\mathbf{\widetilde{w}}_{xyz} =: \boldsymbol{\tilde{\mu}}_x.
    \label{eqxmod}
\end{equation}
where the block matrix $\mathbf{P}_x$ is given by
\begin{equation}
     \mathbf{P}_x = \begin{bmatrix}
         v_{x1}& v_{x2}& v_{x3}\\
          v_{x1}^2& v_{x2}^2& v_{x3}^2\\
           v_{x1}^3& v_{x2}^3& v_{x3}^3\\
     \end{bmatrix},
     \label{pexdot}
\end{equation}
and the corresponding block in weight vector is given by
\begin{equation}
     \mathbf{w}_x = \begin{bmatrix}
         w_1&w_2&w_3
     \end{bmatrix}^T.
     \label{wexdot}
\end{equation}

If we let $s_{\rm{twin}} =\frac{s}{\delta}$ and $s_{\rm{quad}} =\frac{s}{\gamma}$ be given in terms of a speed parameter, $s$, and scaling parameters, $\delta$ and $\gamma$, then the right hand side of \eqref{eqxmod} is given by
\begin{equation}
    \boldsymbol{\widetilde{\mu}}_x=\begin{bmatrix}
        M_{100}-\frac{\delta^2}{s^2}(M_{120}+M_{102})\\
        M_{200}-|M_{110}|-|M_{101}|-\frac{\delta}{s}(|M_{120}|+|M_{210}|+|M_{102}|+|M_{201}|)-\frac{\gamma}{s}|M_{111}|\\
        M_{300}-M_{120}-M_{102}
    \end{bmatrix}.
\label{mexdot}
\end{equation}


To facilitate the solution of \eqref{eqxmod}, we choose  the velocities to be of the form
\begin{equation}
    \begin{bmatrix}v_{x,1} & v_{x,2} & v_{x,3}\end{bmatrix} 
    = 
    \begin{bmatrix}\beta_xs & -\beta_xs & \beta_xl_xs\end{bmatrix}, 
\label{xvalsdot}
\end{equation}
where $l_x$ is a positive constant satisfying $l_x\neq1$ and $\beta_x$ can be either $1$ or $-1$.  The effects of different choices of $\beta_x$ and $l_x$ will be explored in section~\ref{sec:Results}.

We solve $\mathbf{P}_x\mathbf{w}_x = \boldsymbol{\widetilde{\mu}}_x$ by substituting  \eqref{xvalsdot} into \eqref{pexdot} and applying Gaussian elimination. The resulting block of the weight vector is given by
\begin{equation}
\renewcommand{\arraystretch}{1.6}
    \mathbf{w}_x = \begin{bmatrix}
        \frac{1}{2s^3}(a_xs-b_x-c_{x+})\\ 
         \frac{1}{2s^3}(a_xs-b_x-c_{x-})\\ 
         \frac{\beta_x(M_{300}+(M_{102}+M_{120})(\delta^2-1))}{(l_x^2-1)l_xs^3}
    \end{bmatrix},
    \label{wexsoldot}
\end{equation}
where the parameters $a_x$, $b_x$, and $c_{x\pm}$ are combinations of values of the moments, $l_x$, $\beta_x$, $\delta$, and $\gamma$. Explicit values of these parameters are given in appendix ~\ref{subsec::parm}. 

The first two weights in \eqref{wexsoldot} will be positive provided that the speed parameter is chosen so that 
\begin{equation}
    a_x s > b_x+c_{x1} 
    \label{eq:cond111}
\end{equation}
and 
\begin{equation}
    a_x s > b_x+c_{x2}
    \label{eq:cond222}. 
\end{equation}
In our standardized system described in section~\ref{sec:StandardizedP}, $a_x\approx1$. For the third particle, we simply choose 
\begin{equation}
    \beta_x = \text{sign}\left[ \frac{(M_{300}+(M_{102}+M_{120})(\delta^2-1))}{(l_x^2-1)l_xs^3}\right].
\end{equation}
Similarly, the solutions of the subsystems involving the second and third diagonal blocks of \eqref{k3sysexp} are given by,
\begin{equation}
\renewcommand{\arraystretch}{1.6}
    \mathbf{w}_y = \begin{bmatrix}
        \frac{1}{2s^3}(a_ys-b_y-c_{y1})\\ 
         \frac{1}{2s^3}(a_ys-b_y-c_{y2})\\
         \frac{\beta_y(M_{030}+(M_{012}+M_{210})(\delta^2-1))}{(l_y^2-1)l_ys^3}
    \end{bmatrix},\ \ \
    \mathbf{w}_z = \begin{bmatrix}
        \frac{1}{2s^3}(a_zs-b_z-c_{z1})\\
         \frac{1}{2s^3}(a_zs-b_z-c_{z2})\\
         \frac{\beta_z(M_{003}+(M_{021}+M_{201})(\delta^2-1))}{(l_z^2-1)l_zs^3}
    \end{bmatrix}, 
    \label{wexsoldot++}
\end{equation}
and the velocity vectors follow the patterns given in \eqref{xvalsdot} and result in similar weight positivity requirements. 

Finally, the value of $w_0$ is given by
\begin{equation}
    w_0 = 1-\mathbf{1}\cdot\begin{bmatrix}
        \mathbf{w}_x&\mathbf{w}_y&\mathbf{w}_z&\mathbf{\widetilde{w}}_{xy}&\mathbf{\widetilde{w}}_{xz}&\mathbf{\widetilde{w}}_{yz}&\mathbf{\widetilde{w}}_{xyz}
    \end{bmatrix}^T .
    \label{wzero}
\end{equation}
Evaluating \eqref{wzero} we find that
\begin{equation}
    w_0 = \frac{s^3-a_0s+b_0}{2s^3},
    \label{eqn:w0_formula}
\end{equation}
where $a_0$, and $b_0$ are combinations of the values of the moments and of $l_{x,y,z}$, $\beta_{x,y,z}$, $\delta$, and $\gamma$. In appendix \ref{subsec::parm} we provide  explicit formulae for $a_0$, and $b_0$. Finally, we note that $w_0>0$ for a sufficiently large value of $s$
that is easy to determine numerically.

In summary, there are three conditions that must hold in order for all the weights in \eqref{wexsoldot}, \eqref{wexsoldot++} and \eqref{eqn:w0_formula} to be nonnegative. The first condition is that 
\begin{equation}
    s\ge\operatorname{max}\left(\frac{b_x+c_{x1}}{a_x},\frac{b_x+c_{x2}}{a_x},\frac{b_y+c_{y1}}{a_y},\frac{b_y+c_{y2}}{a_y},\frac{b_z+c_{z1}}{a_z},\frac{b_z+c_{z2}}{\alpha_z}\right).
\label{eqn:s_min}\end{equation}
Second, we must choose $\beta_x$ such that 
\begin{equation}
    \beta_x = \text{sign}\left[ \frac{(M_{300}+(M_{102}+M_{120})(\delta^2-1))}{(l_x^2-1)l_xs^3}\right],
\end{equation}
with analogous formulae holding for  $\beta_y$ and $\beta_z$.
Finally, we need to choose $s$ large enough to ensure that $w_0$ is positive.  

\section{Numerical Results}\label{sec:Results}
\subsection{Baseline numerical uncertainty}
The goal of particle simulations like DSMC and SWPM is to accurately and efficiently compute the velocity pdf of a system of particles due to physical processes such as particle collisions, transport, and interactions with surfaces. In practice, rather than attempting to accurately compute the entire velocity pdf, researchers typically aim to accurately compute certain statistics of the pdf, such as moments or tail functionals. 

Since the particle grouping and reduction processes in a SWPM simulation do not have counterparts in physical systems or DSMC simulations, it is important that these processes do not significantly contribute to error in the statistics of the pdf. To assess the efficacy of our new reduction scheme, we must first determine the numerical uncertainty in the empirical measure of a SWPM simulation prior to reduction and compare that to the uncertainty in a DSMC simulation. Once these uncertainties have been quantified we can assess if the numerical uncertainty of the reduced system is on the order of the uncertainty prior to reduction.  Any significant increase in the uncertainty will introduce non-physical numerical errors, while aiming for much less uncertainty is pointless. 

Rather than performing actual DSMC and SWPM simulations that realistically model particle physics, for the results in this paper we manufacture particle systems prior to the grouping and reduction stages of an SWPM simulation by sampling from analytical Maxwellian-like distributions with prescribed skewness and kurtosis. This approach enables us to focus attention on the uncertainty introduced by the particle grouping and reduction processes. 

To model a one-dimensional Maxwellian-like distribution with skewness, we use the distribution~\cite{OHAGAN1976}, 
\begin{equation}
    f(v;\alpha)=2 \phi(v) \Phi(\alpha v),
\end{equation}
where $\alpha$ is a skewness parameter, $\phi(v) = \frac{1}{\sqrt{2\pi}}e^{-v^2/2}$ is a Maxwellian, and $\Phi$ is the associated cumulative distribution function (CDF),
\begin{equation}
    \Phi(v)=\int_{-\infty}^{v} \phi(t) dt = 1 + \operatorname{erf}(v/2).
\end{equation}
To introduce kurtosis, we multiply $f$ by a Druyvesteyn distribution, yielding
\begin{equation}\label{eqn:SkewedMaxwellianAndKurt}
    f(v;\alpha,\beta)=2 \phi(v) \Phi(\alpha v) e^{-\beta v^4/4},
\end{equation}
where the parameter, $\beta$, controls the degree of kurtosis.
To obtain a three dimensional distribution we set
\begin{equation}
    f(\mathbf{v};\boldsymbol\alpha,\boldsymbol\beta) = 
    C f(v_x;\alpha_x,\beta_x) f(v_y;\alpha_y,\beta_y) 
    f(v_z;\alpha_z,\beta_z),
\label{eqn:CompleteDistribution}\end{equation}
where the normalization constant is numerically determined by the condition that
\begin{equation}
    \iiint_{\left|\mathbf{v}\right|\leq v_R}
    f\left(\mathbf{v};\boldsymbol\alpha,\boldsymbol\beta\right) d\mathbf{v} = 1.
\label{eqn:distributionNormalization}
\end{equation}
Unless otherwise noted, we chose $v_R=7$ to allow accurate computation of the tail functional,
\begin{equation}
    \text{Tail}(R) =  \iiint_{\left|\mathbf{v}\right|\geq R} f\left(\mathbf{v}\right) d\left(\mathbf{v}\right),
\label{eqn:TailFunctional}
\end{equation}
out to $R=6$.

\begin{figure}[!ht]
        \includegraphics[width=.495\textwidth]{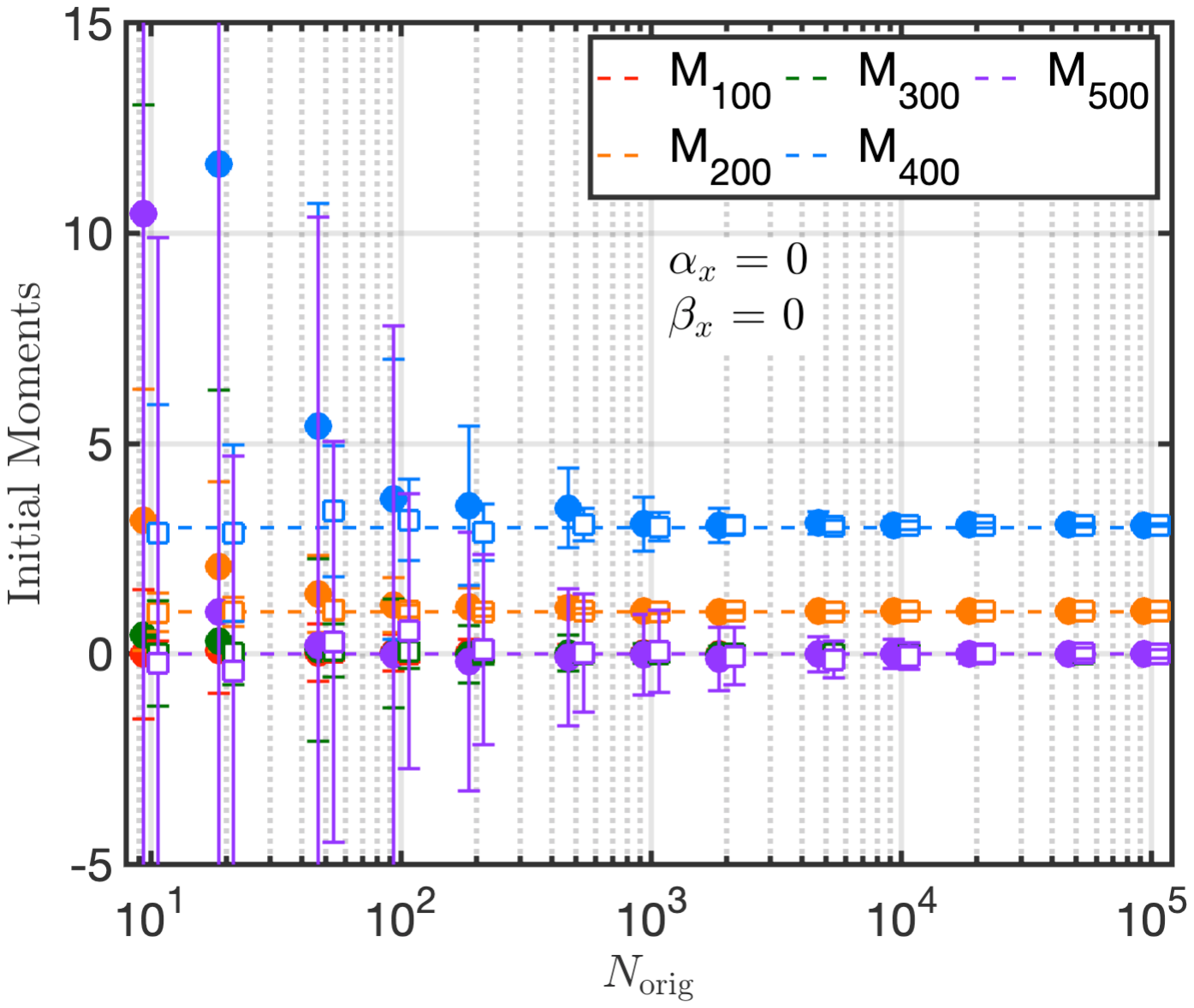}\label{fig:InitialMomentSkew_0,Kurt_0}
        \includegraphics[width=.495\textwidth]{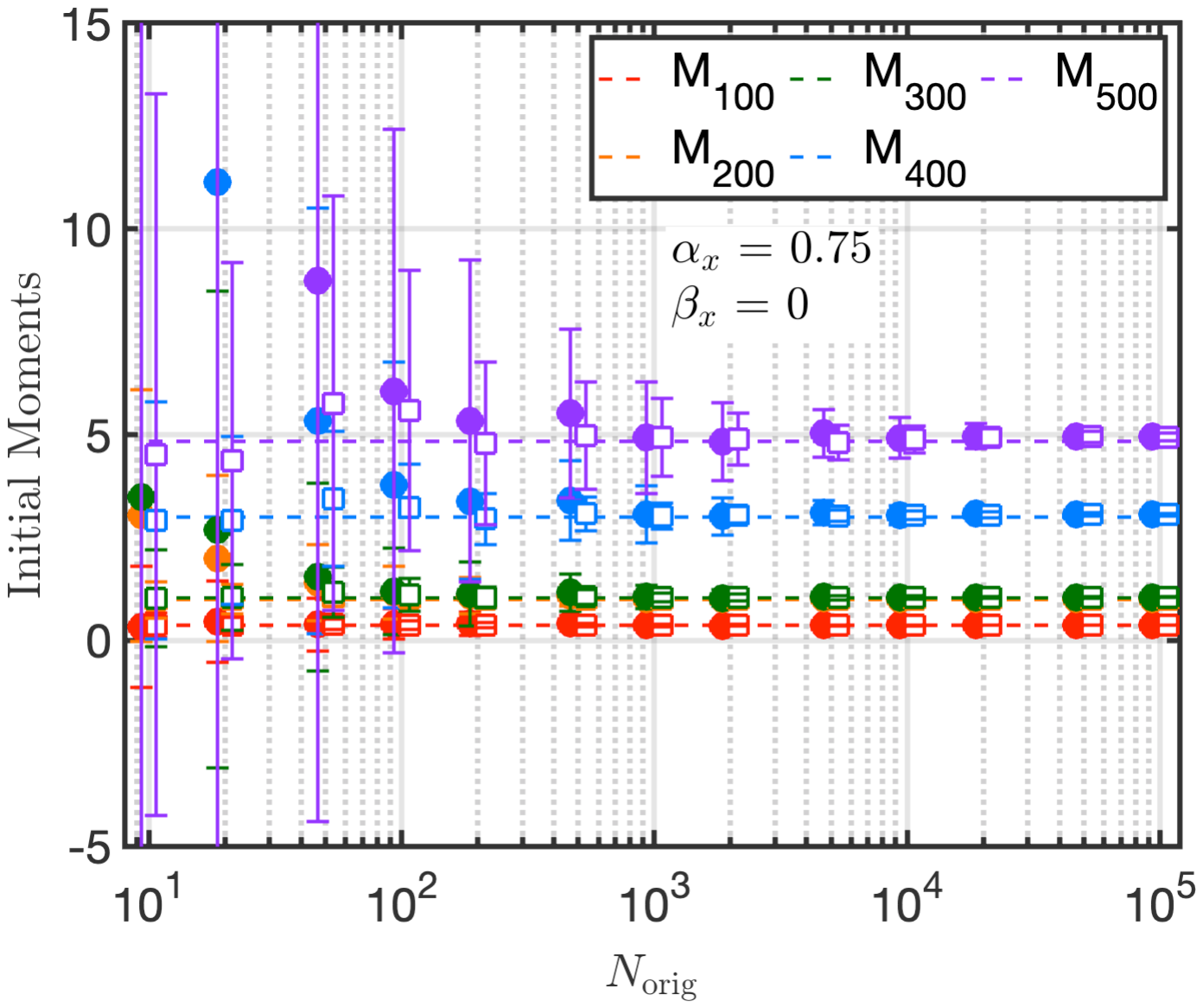}\label{fig:InitialMomentSkew_75,Kurt_0}
        \includegraphics[width=.495\textwidth]{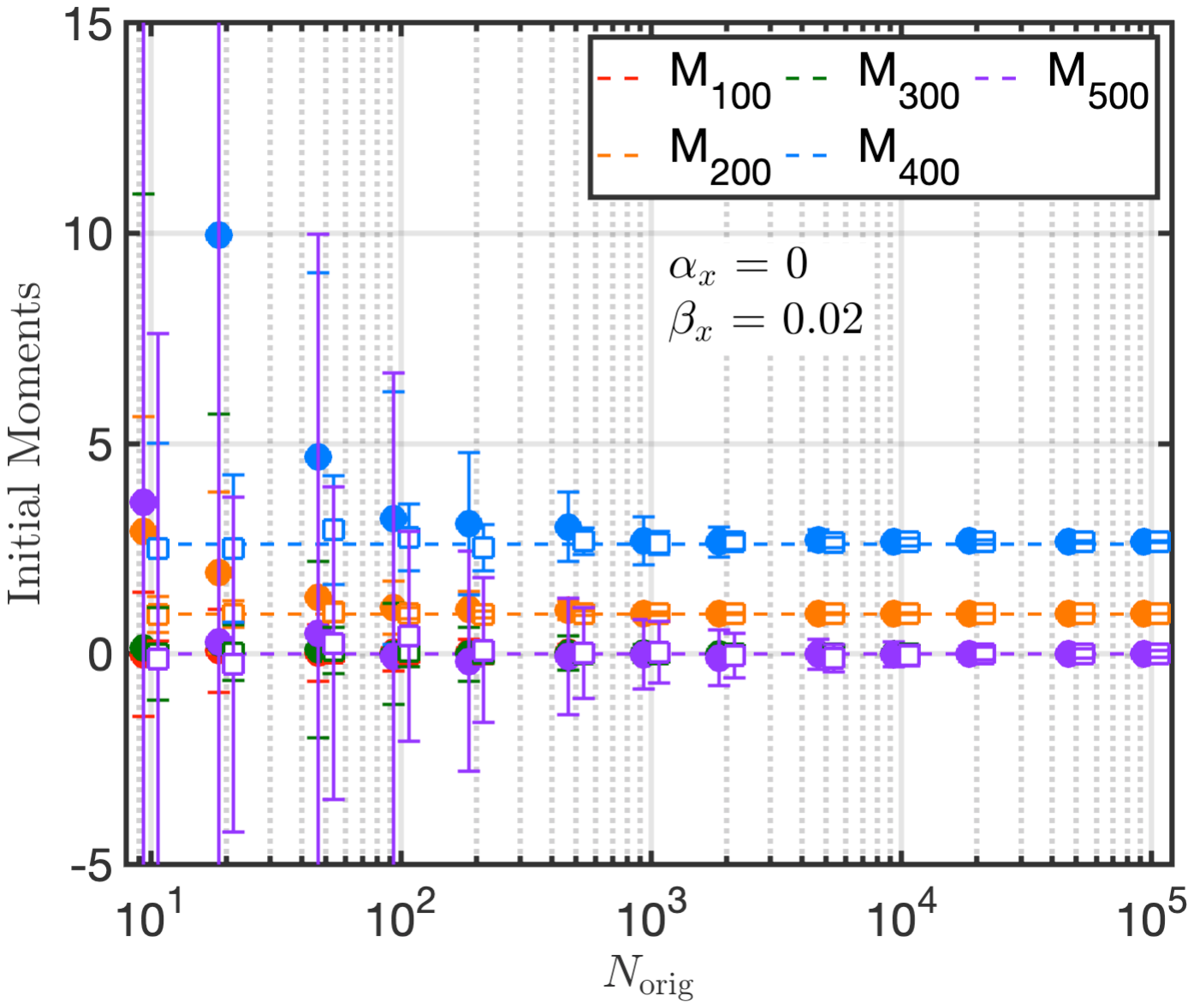}\label{fig:InitialMomentSkew_0,Kurt_02}
        \includegraphics[width=.495\textwidth]{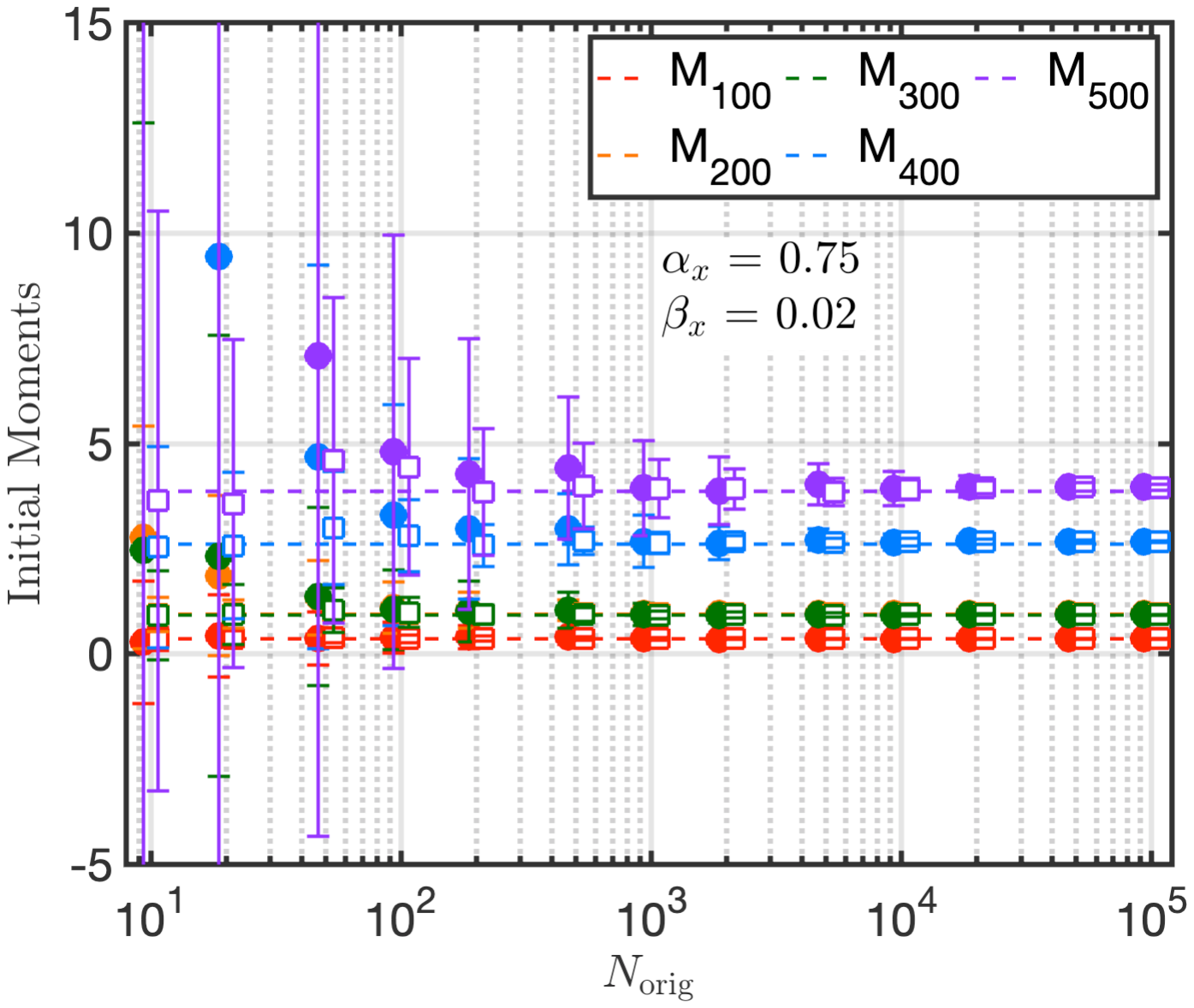}\label{fig:InitialMomentSkew_75,Kurt_02}
\caption{
{Comparison of the  moments obtained using direct integration of \eqref{eqn:CompleteDistribution} (dashed lines) vs. those obtained with 100 ensembles of DSMC-like (open squares) and SWPM-like (closed circles) simulations. For this comparison we consider the moments along the $v_x$ axis, $M_{100}$ ({red}), $M_{200}$ ({orange}), $M_{300}$ ({green}), $M_{400}$ ({blue}) and $M_{500}$ ({violet}) with various levels of skewness and kurtosis along the $v_x$ axis. Top left: $\boldsymbol{\alpha}=(0,0,0)$, $\boldsymbol{\beta}=(0,0,0)$; Top right: $\boldsymbol{\alpha}=(0.75,0,0)$, $\boldsymbol{\beta}=(0,0,0)$; Bottom left: $\boldsymbol{\alpha}=(0,0,0)$, $\boldsymbol{\beta}=(0.02,0,0)$; Bottom right: $\boldsymbol{\alpha}=(0.75,0,0)$, $\boldsymbol{\beta}=(0.02,0,0)$.}
}
\label{fig:InitialMoments}\end{figure}

\begin{figure}[!ht]
\begin{center}
        \includegraphics[width=.495\textwidth]{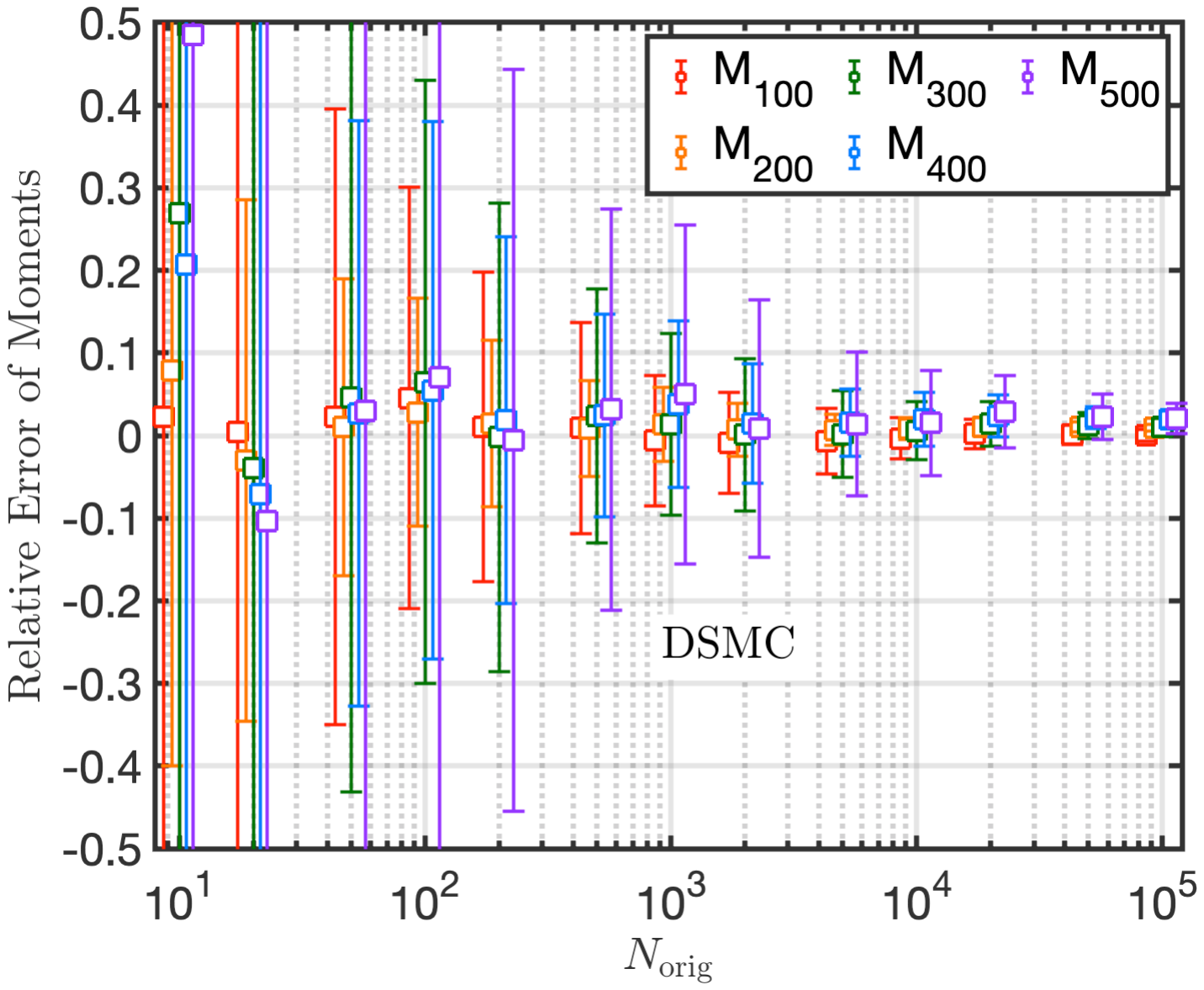}
         \includegraphics[width=.495\textwidth]{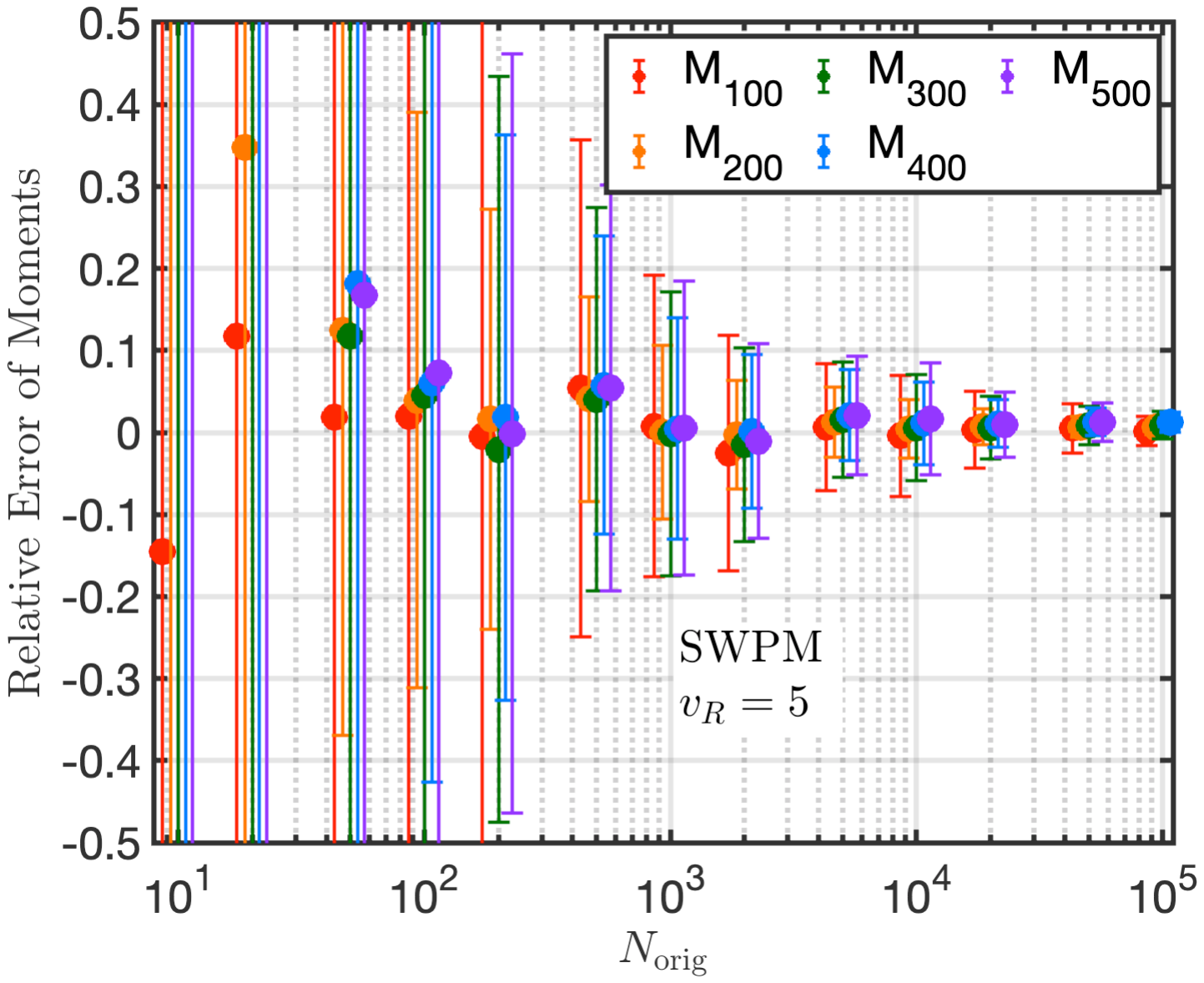}
        \includegraphics[width=.495\textwidth]{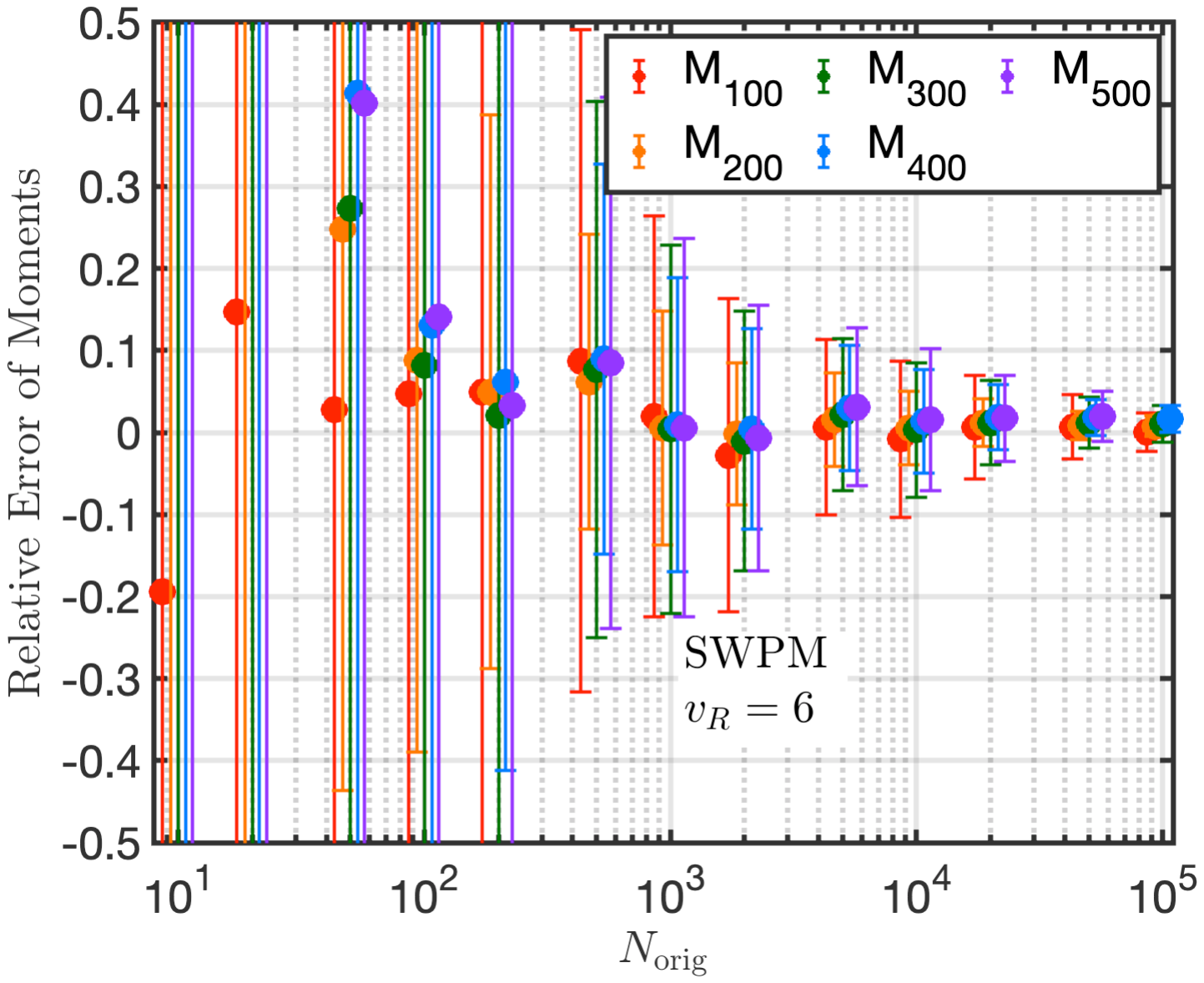}
        \includegraphics[width=.495\textwidth]{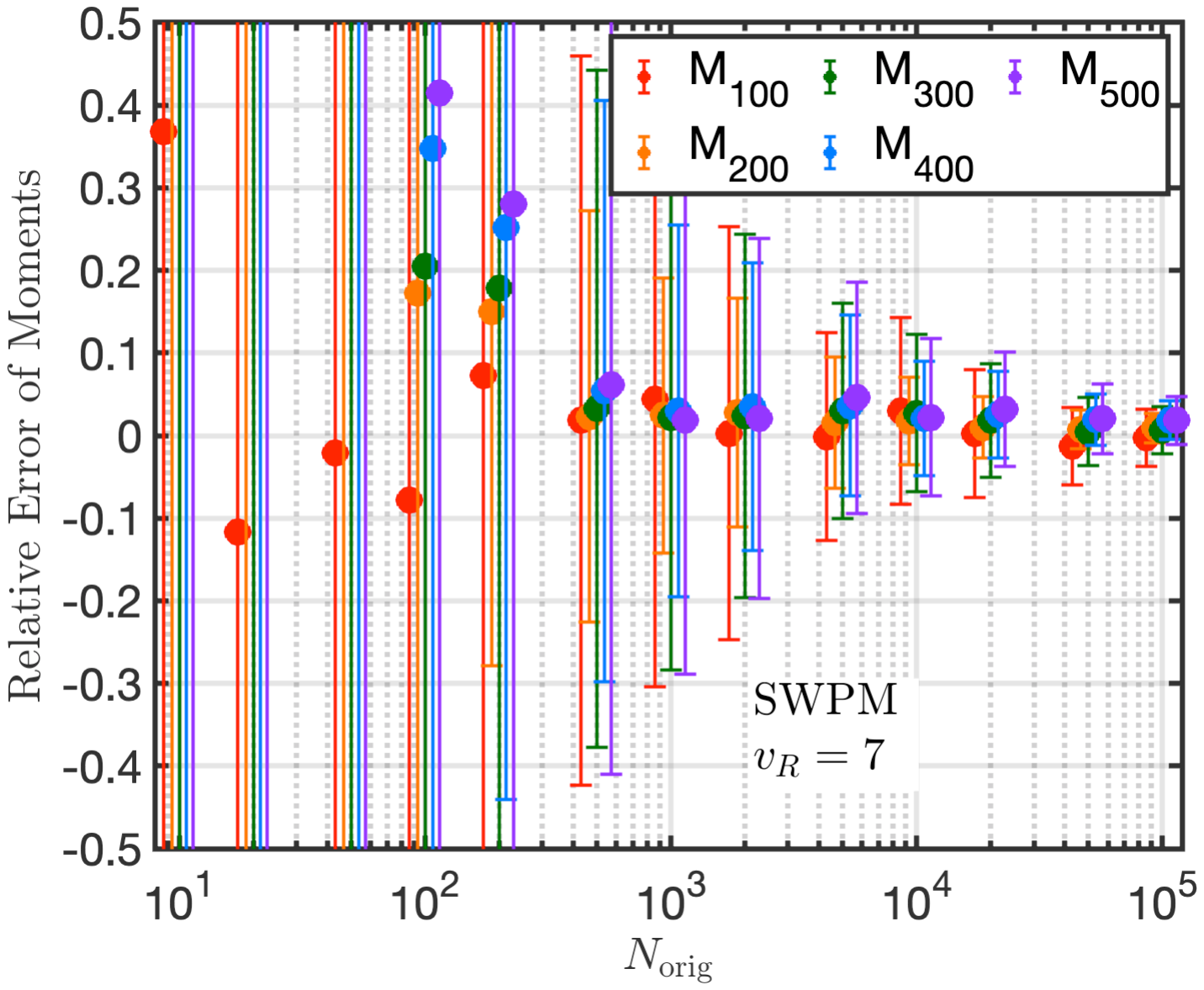}
\end{center}
\caption{
Error in the moments computed using particle simulations relative to those obtained by numerical integration of the analytical distribution for DSMC-like simulations (top left), and SWPM-like simulations with $v_R=5$ (top right), $v_R=6$ (bottom left), and $v_R=7$ (bottom right). Here, $\boldsymbol{\alpha}=(0.75,0,0)$, $\boldsymbol{\beta}=(0.02,0,0)$.
}
\label{fig:RelativeErrorInitialMoments}
\end{figure}

\begin{figure}[!h]
        \includegraphics[width=.495\textwidth]{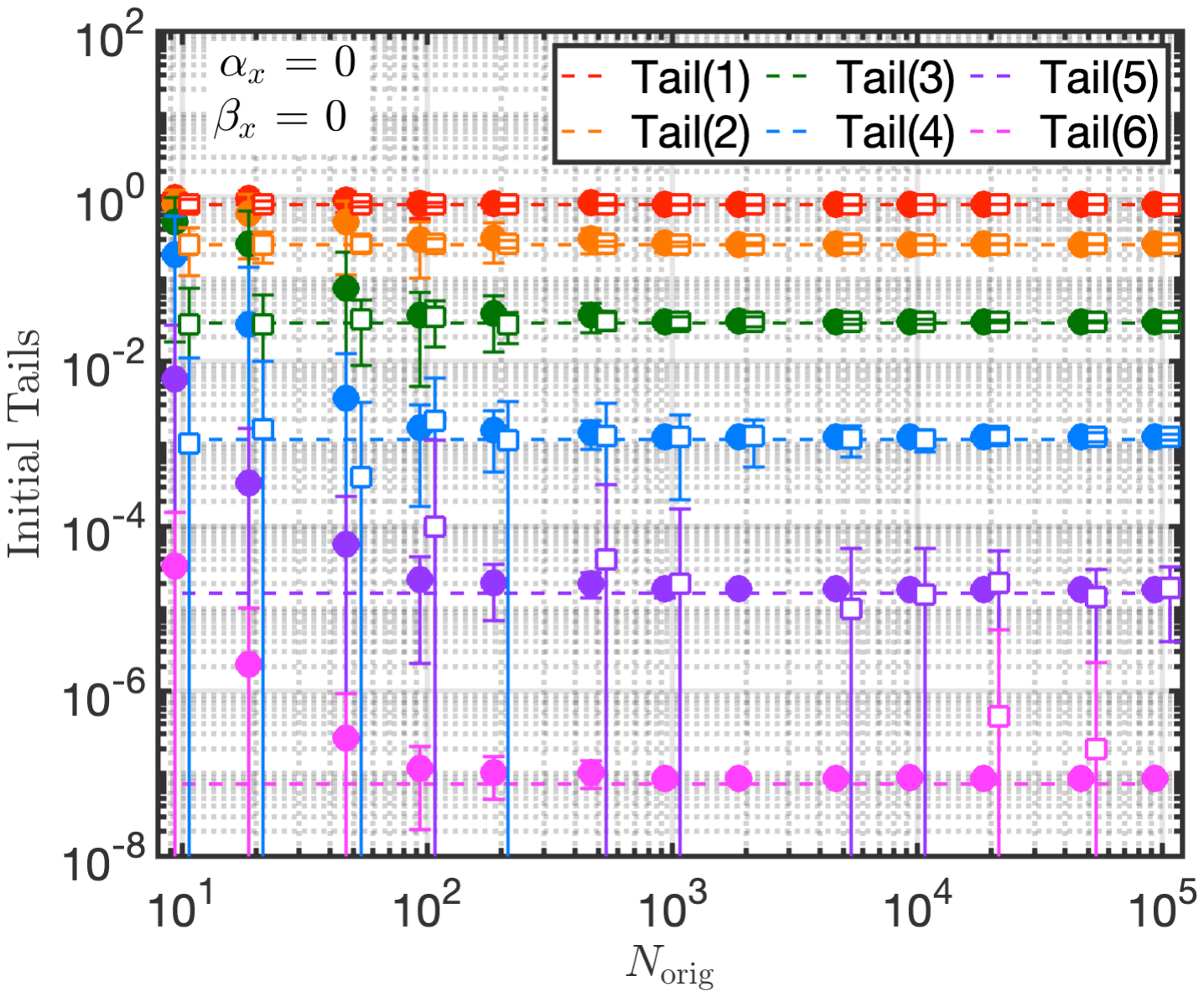}\label{fig:InitialTailSkew_0,Kurt_0}
        \includegraphics[width=.495\textwidth]{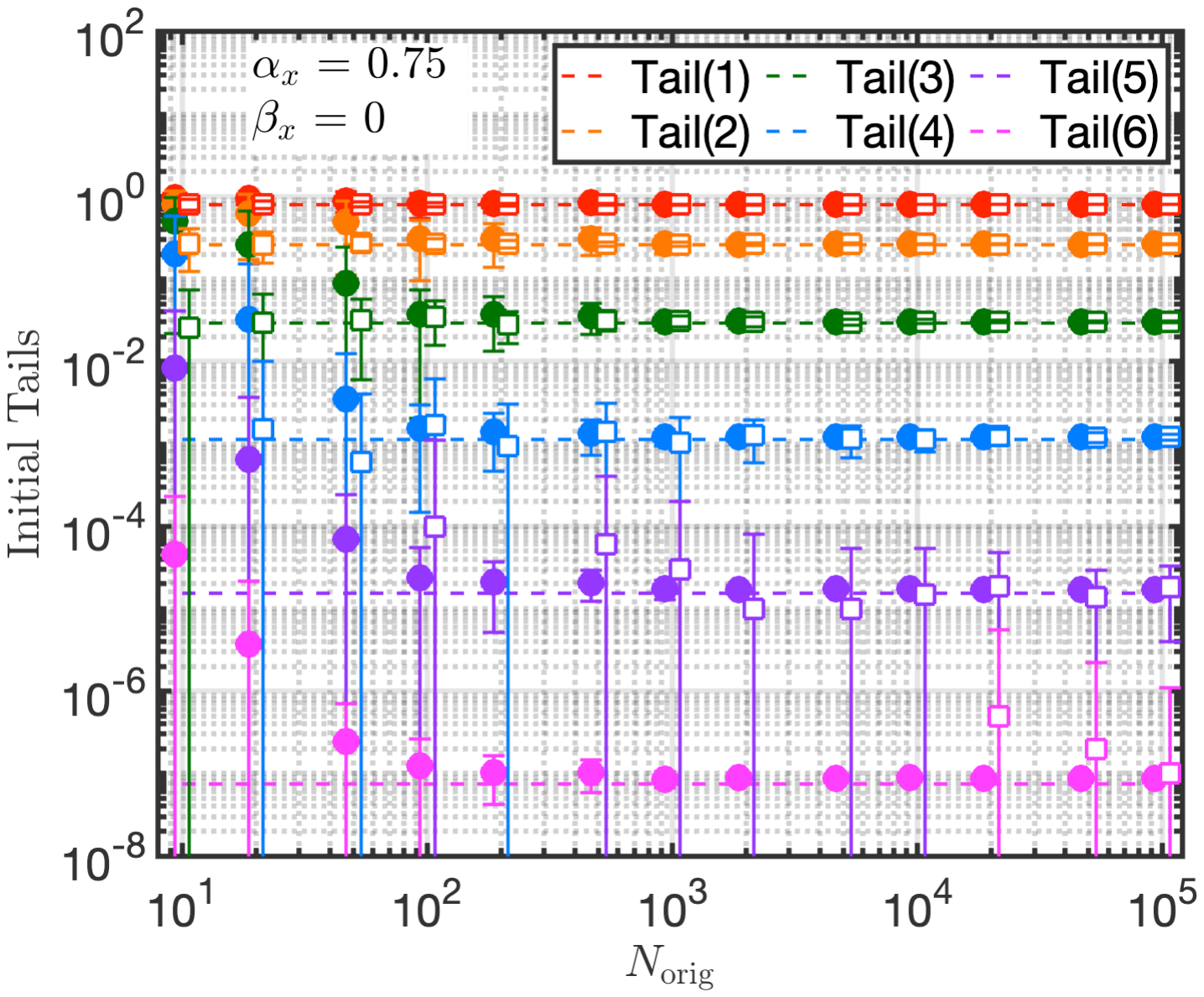}\label{fig:InitialTailSkew_75,Kurt_0}
        \includegraphics[width=.495\textwidth]{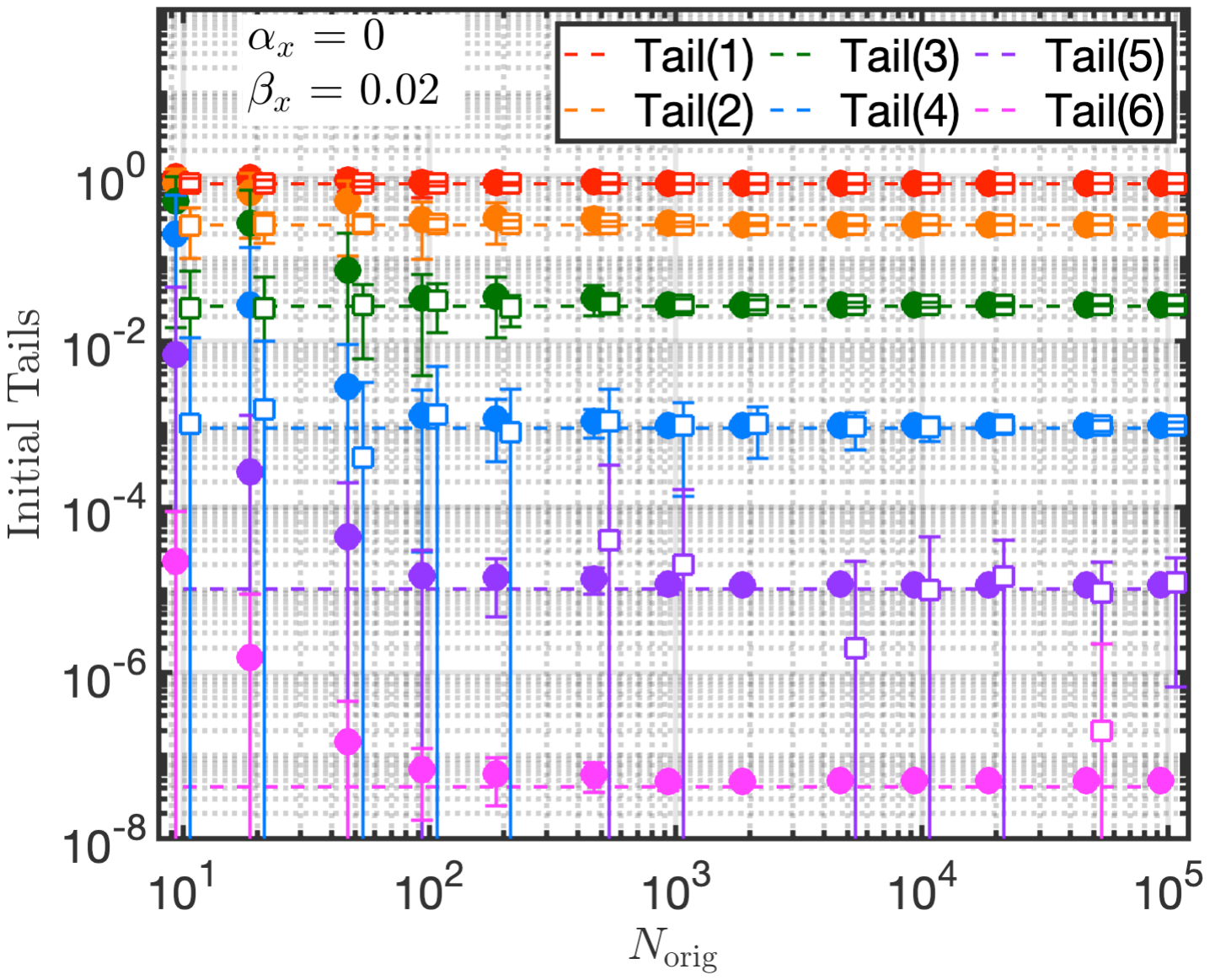}\label{fig:InitialTailSkew_0,Kurt_02}
        \includegraphics[width=.495\textwidth]{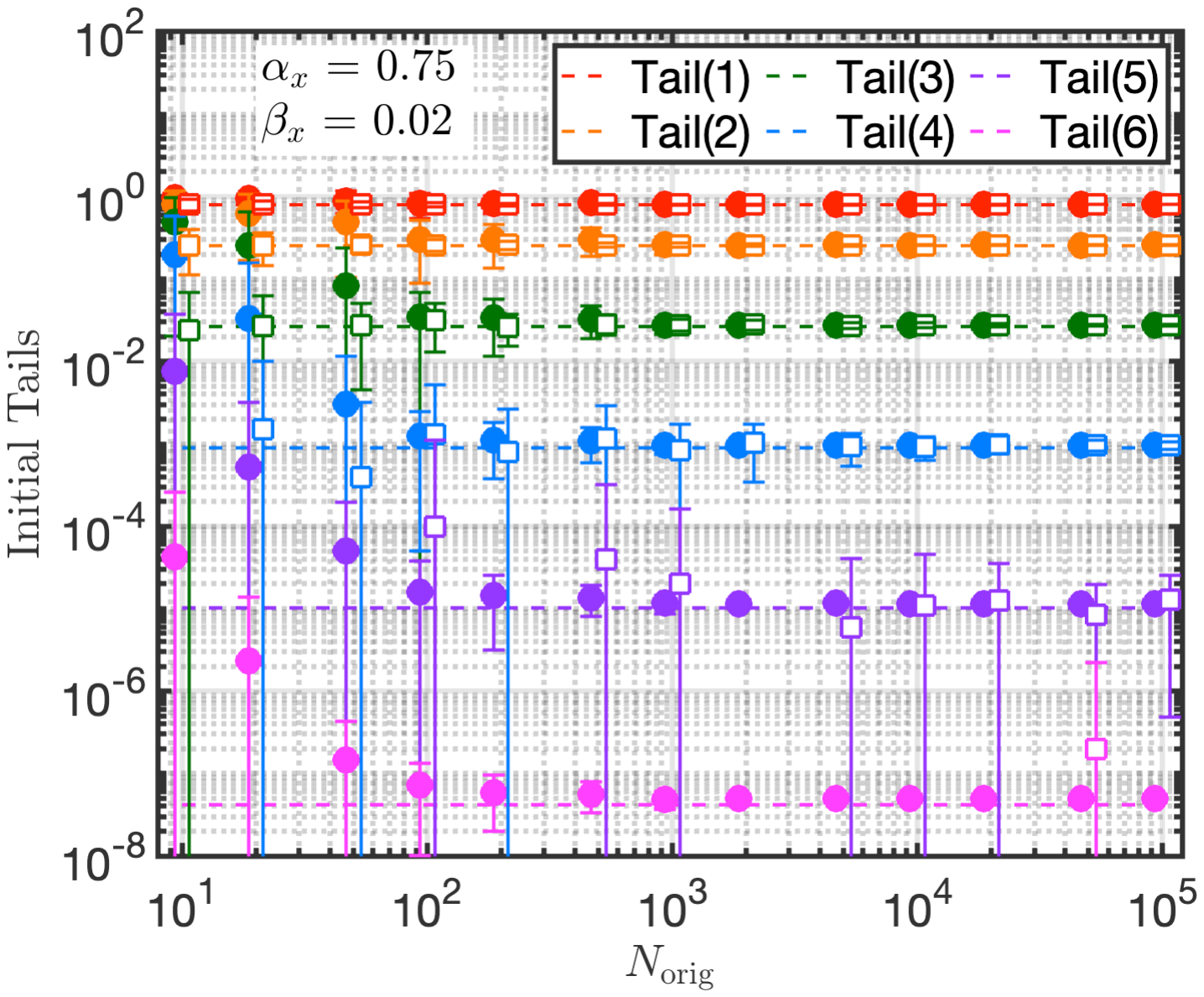}\label{fig:InitialTailSkew_75,Kurt_02}
\caption{Comparison of the tail functionals obtained using direct integration of \eqref{eqn:CompleteDistribution} using \eqref{eqn:TailFunctional} (dashed lines) vs. those one might find with 100 ensembles of DSMC (open squares) or SWPM (closed circles) simulations. 
}
\label{fig:InitialTails}
\end{figure}

Sampling from \eqref{eqn:CompleteDistribution}, we generate systems of stochastic particles that would typically be encountered in a DSMC or SWPM simulation. To generate a DSMC-like system, we randomly sample $N_{\rm{orig}}$ particle velocities using the cumulative distribution functions (CDFs) obtained by numerically integrating \eqref{eqn:SkewedMaxwellianAndKurt} in each direction independently, and setting the weights to $w_i=\frac{1}{N_{\rm orig}}$. This results in stochastic particles that are highly centralized near the peak of \eqref{eqn:CompleteDistribution}. 

To generate a SWPM-like system, we randomly sample $N_{\rm{orig}}$ particle velocities which  are uniformly distributed in a ball of radius $v_R=7$, which enables us to accurately compute $\operatorname{Tail}(6)$. We refer to the ball $|\mathbf v| < v_R$ as the simulation ball, since all the original stochastic particles in the simulation lie in this ball. After the velocities are determined, the weight of each stochastic particle is set to ensure that the particles distribution is given by \eqref{eqn:SkewedMaxwellianAndKurt}. Once all of the particles have been sampled, the weights are normalized so that \eqref{eqn:distributionNormalization} holds.

In Figure~\ref{fig:InitialMoments},  we  examine the accuracy and precision of the moments  computed using the DSMC-like and SWPM-like simulations.  We compute the mean and standard deviation for the moments $M_{100}$, $M_{200}$, $M_{300}$, $M_{400}$, and $M_{500}$, obtained from 100 ensembles, each with $N_{\rm orig}$ stochastic particles for several sets of parameters $(\boldsymbol\alpha,\boldsymbol\beta)$ in \eqref{eqn:SkewedMaxwellianAndKurt}. We chose four parameter sets for which $\boldsymbol{\alpha} = (0, 0, 0)$ or $(0.75, 0, 0)$ and $\boldsymbol{\beta} = (0, 0, 0)$ or $(0.02, 0, 0)$. We computed these statistics with $N_{\rm orig} = 10$, $20$, 50, 100, $\cdots$ 100,000, covering the typical range of the values in a simulation cell. We show the mean and uncertainty of the five moments for the DSMC-like particle systems (open squares) and SWPM-like systems (closed circles) as a function of $N_{\rm orig}$. The true values of the moments, obtained via numerical integration using the analytical formula for the distribution in \eqref{eqn:CompleteDistribution}, are shown with dashed lines.
For each parameter set, the accuracy and precision of the two techniques improve as  $N_{\rm orig}$ increases. 
 
In Figure \ref{fig:RelativeErrorInitialMoments}, we show the relative error in the moments for $\boldsymbol{\alpha}=(0.75,0,0)$ and $\boldsymbol{\beta}=(0.02,0,0)$, which corresponds to the lower right panel in Figure~\ref{fig:InitialMoments}. In the top left panel, we show that for DSMC all the moments are accurate to within a tolerance of $50\%$ for $N_{\rm orig}\ge500$, with a confidence level of $95\%$. For a tighter error tolerance of $20\%$ we require $N_{\rm orig}\ge5000$.  In comparison, for the SWPM method with $v_R=7$ (botttom right), we need $N_{\rm orig}\ge2000$ to achieve an error tolerance of 50\%, which increases to $N_{\rm orig}\ge20,000$ for a $20\%$ error tolerance. In the top right panel, we show that to obtain the same level of accuracy and precision as DSMC for a comparable value of $N_{\rm orig}$, all samples must lie in the ball of radius, $v_R=5$.

In Figure~\ref{fig:InitialTails}, we show the corresponding results for the tail functionals, $\operatorname{Tail}(1)$ through  $\operatorname{Tail}(6)$.  In contrast to the moments, SWPM can generally be used to accurately and precisely determine the first through sixth tail functionals for all $N_{\rm orig}\geq 100$.  A larger value of $v_{R}$ would allow us to probe further out into the tails, with some increase in $N_{\rm orig}$, while a smaller value would decrease how far out we can probe them.  In comparison, while DSMC is accurate and precise in for determining Tail$(1)$ and Tail$(2)$, for almost any $N_{\rm{orig}}$, it is difficult to compute  Tail$(N)$ for $N\geq 5$  with DSMC.  For example, in many cases DSMC cannot be used to accurately and precisely determine Tail$(5)$ unless $N_{\rm{orig}}\gtrsim 100,000$, or Tail$(6)$ unless $N_{\rm{orig}}\gg100,000$. Increasing or decreasing $v_{R}$ does not have much impact on the ability of DSMC to probe the tail functionals, unless $v_{R}\lesssim 4$.

 \begin{figure}[t]  
        \includegraphics[width=.495\textwidth]{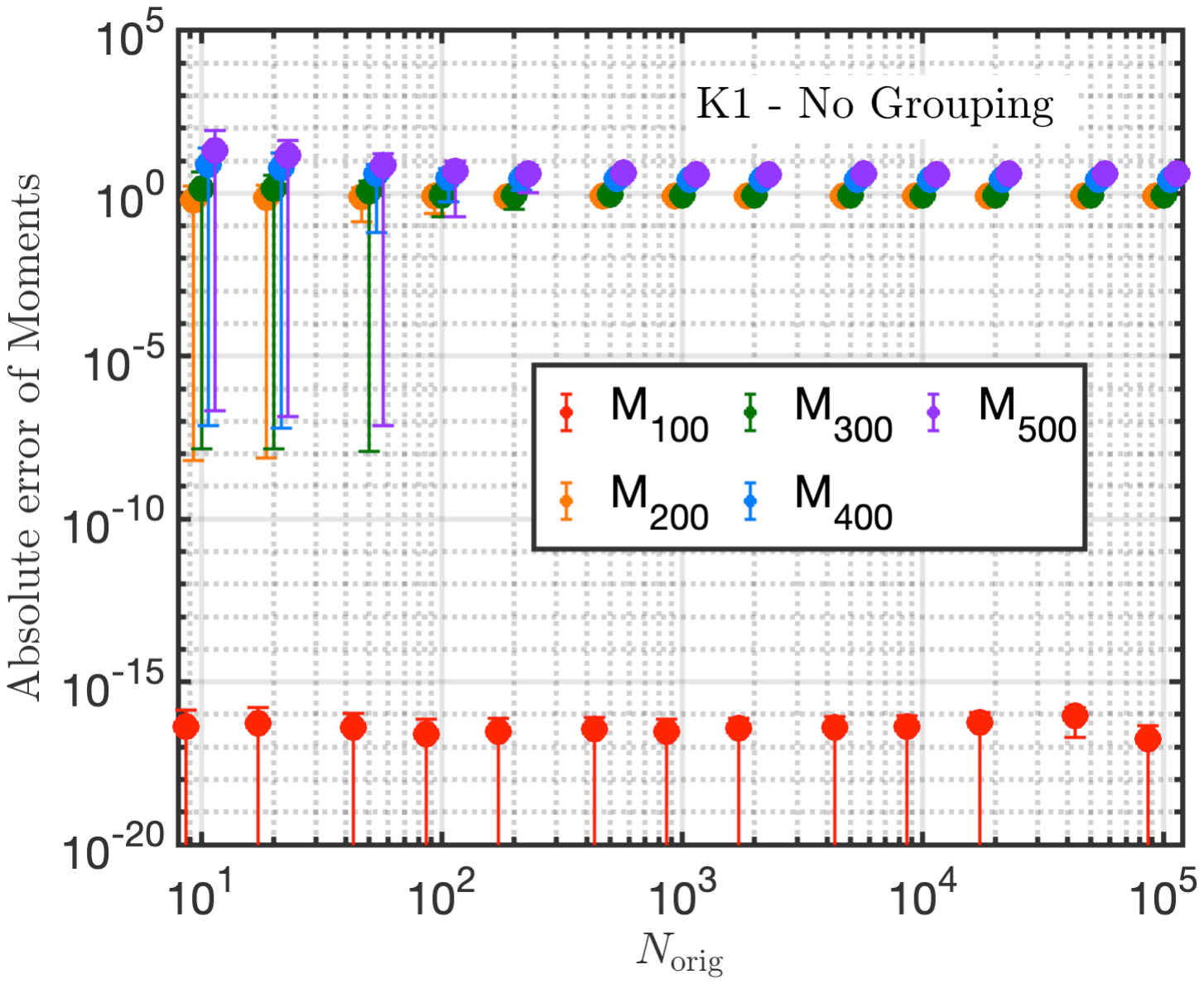}\label{fig:MomentAbsErrorNoGroupK1}
        \includegraphics[width=.495\textwidth]{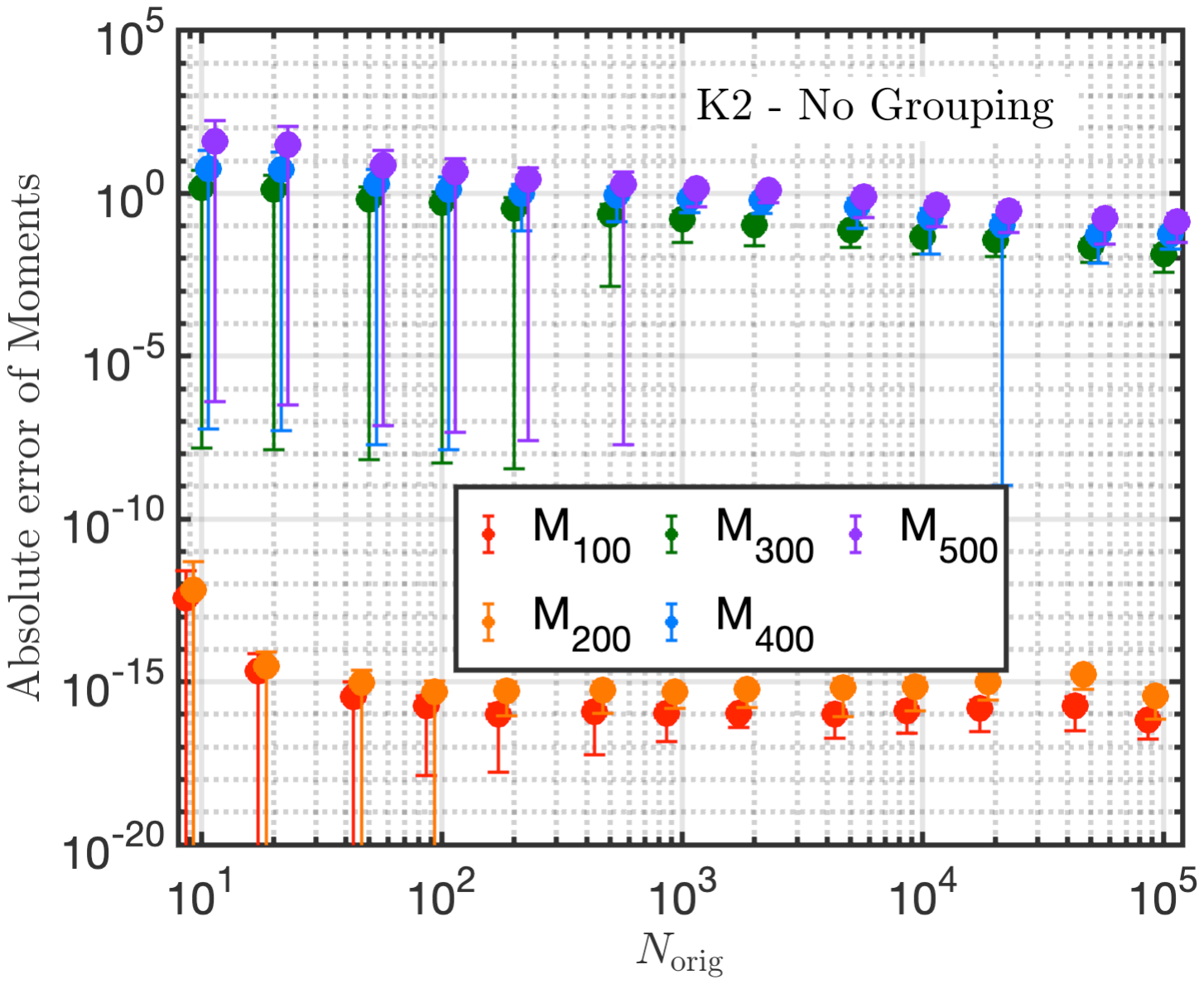}\label{fig:MomentAbsErrorNoGroupK2}
        \includegraphics[width=.495\textwidth]{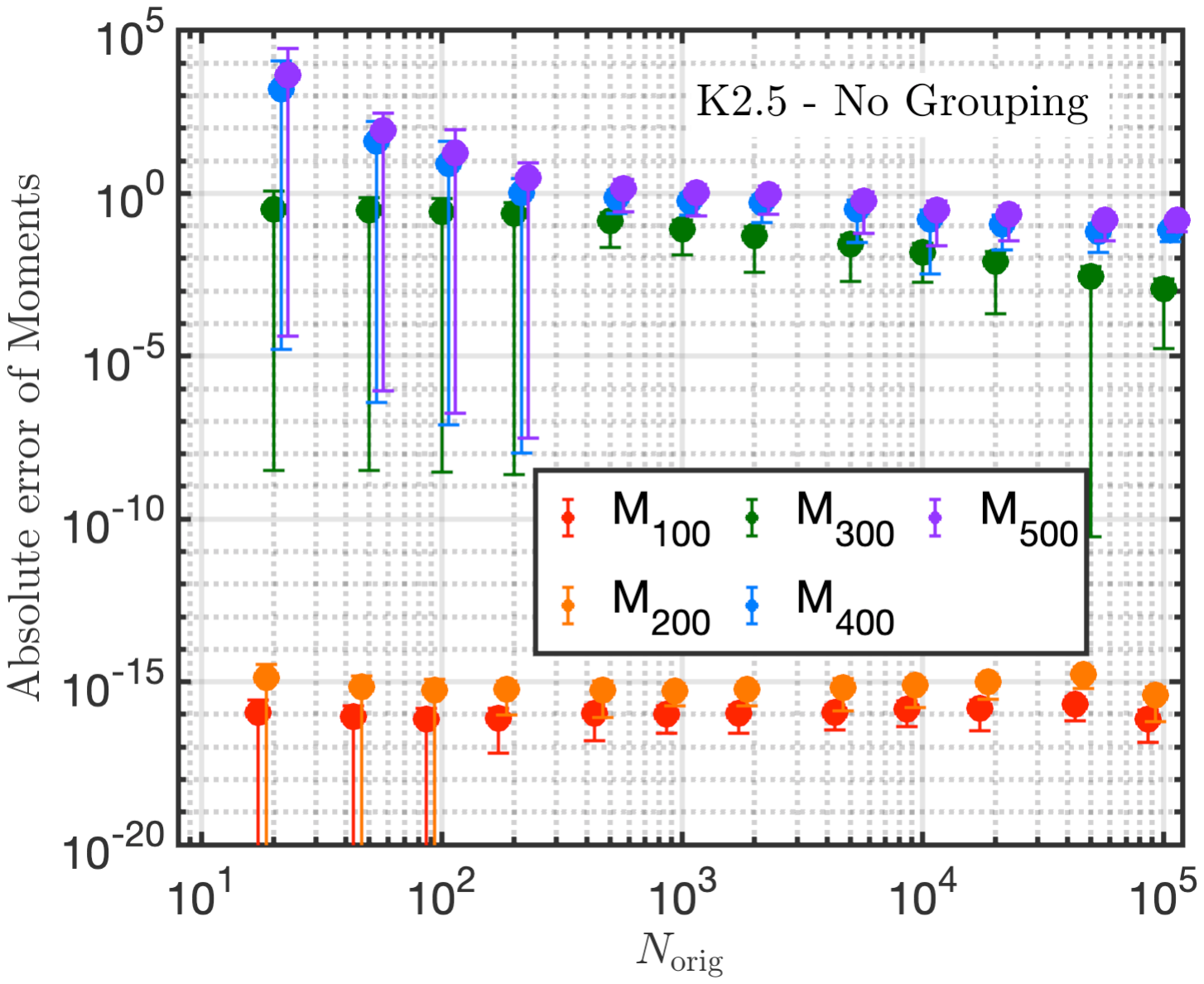}\label{fig:MomentAbsErrorNoGroupK2_5}
        \includegraphics[width=.495\textwidth]{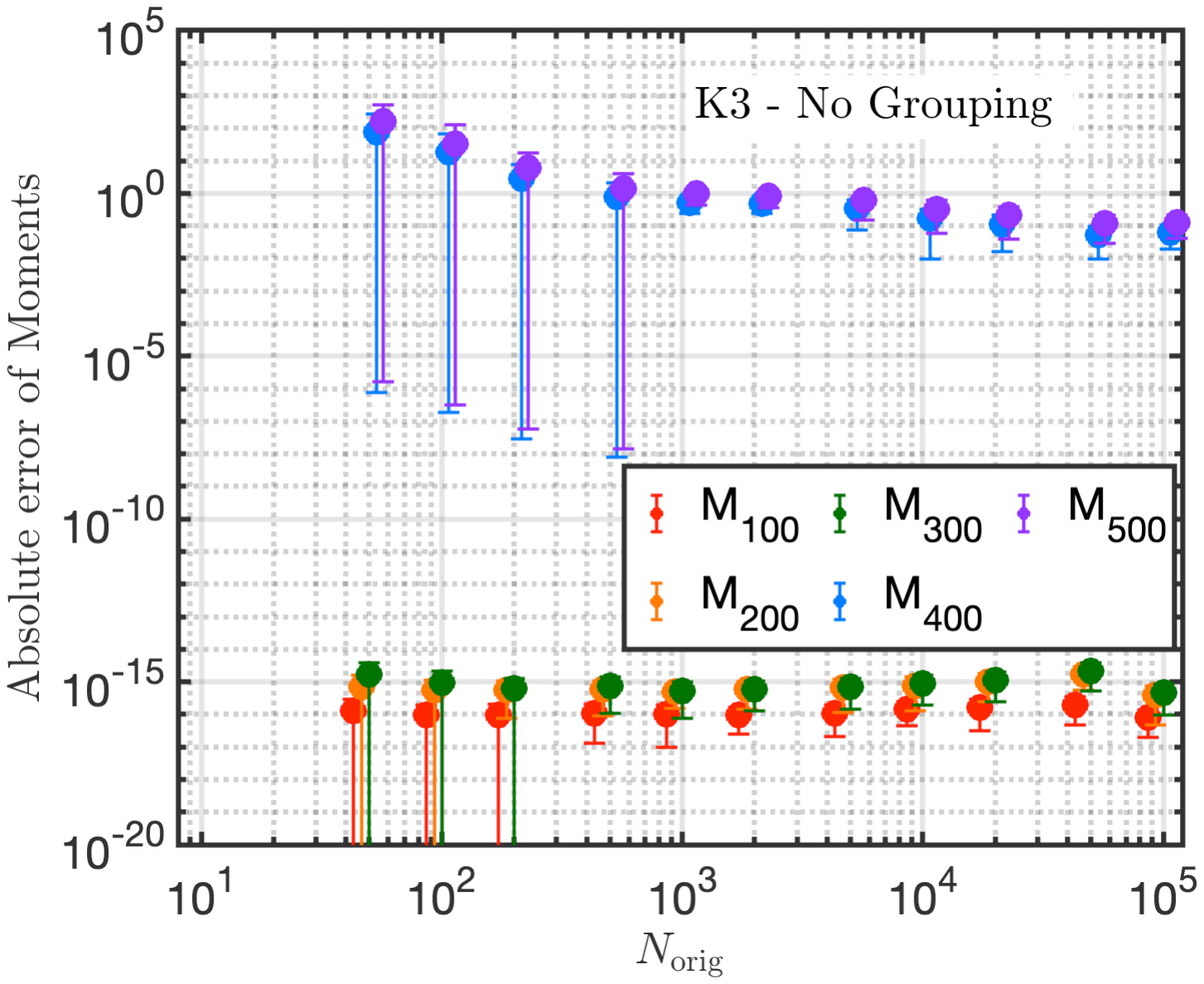}\label{fig:MomentAbsErrorNoGroupK3}
\caption{Absolute error in the moments along the $v_x$ axis due to the four reduction schemes without grouping. Top left: {\bf K1}; Top right: {\bf K2}; Bottom left: {\bf K2.5}; Bottom right: {\bf K3}. Here, $\boldsymbol{\alpha}=(0.75,0,0)$, $\boldsymbol{\beta}=(0.02,0,0)$, and $v_{R}=7$.}
\label{fig:MomentErrorWithoutGrouping}
\end{figure}

\begin{figure}[t]
        \includegraphics[width=.495\textwidth]{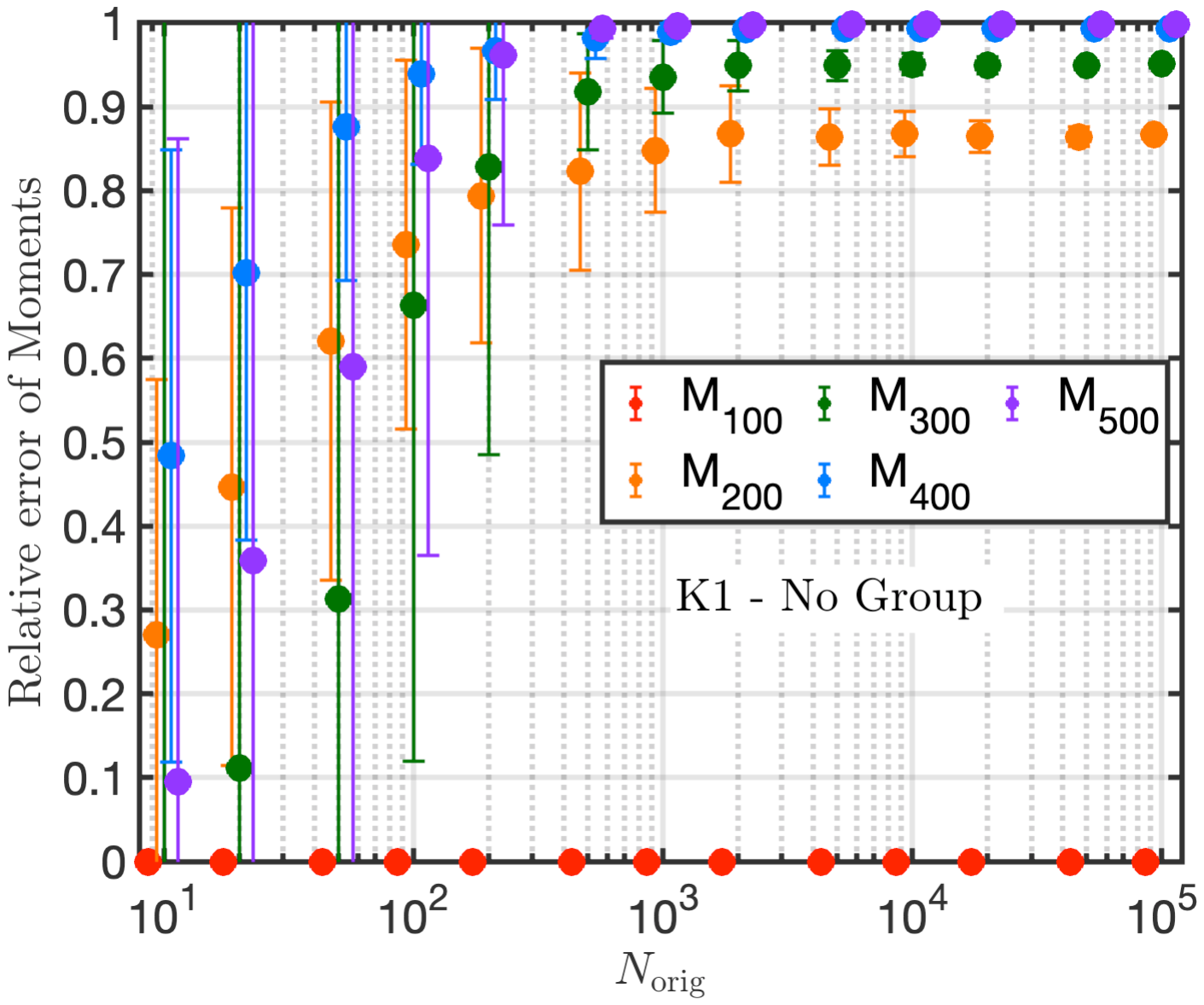}\label{fig:RelMomentErrorNoGroupK1}
        \includegraphics[width=.495\textwidth]{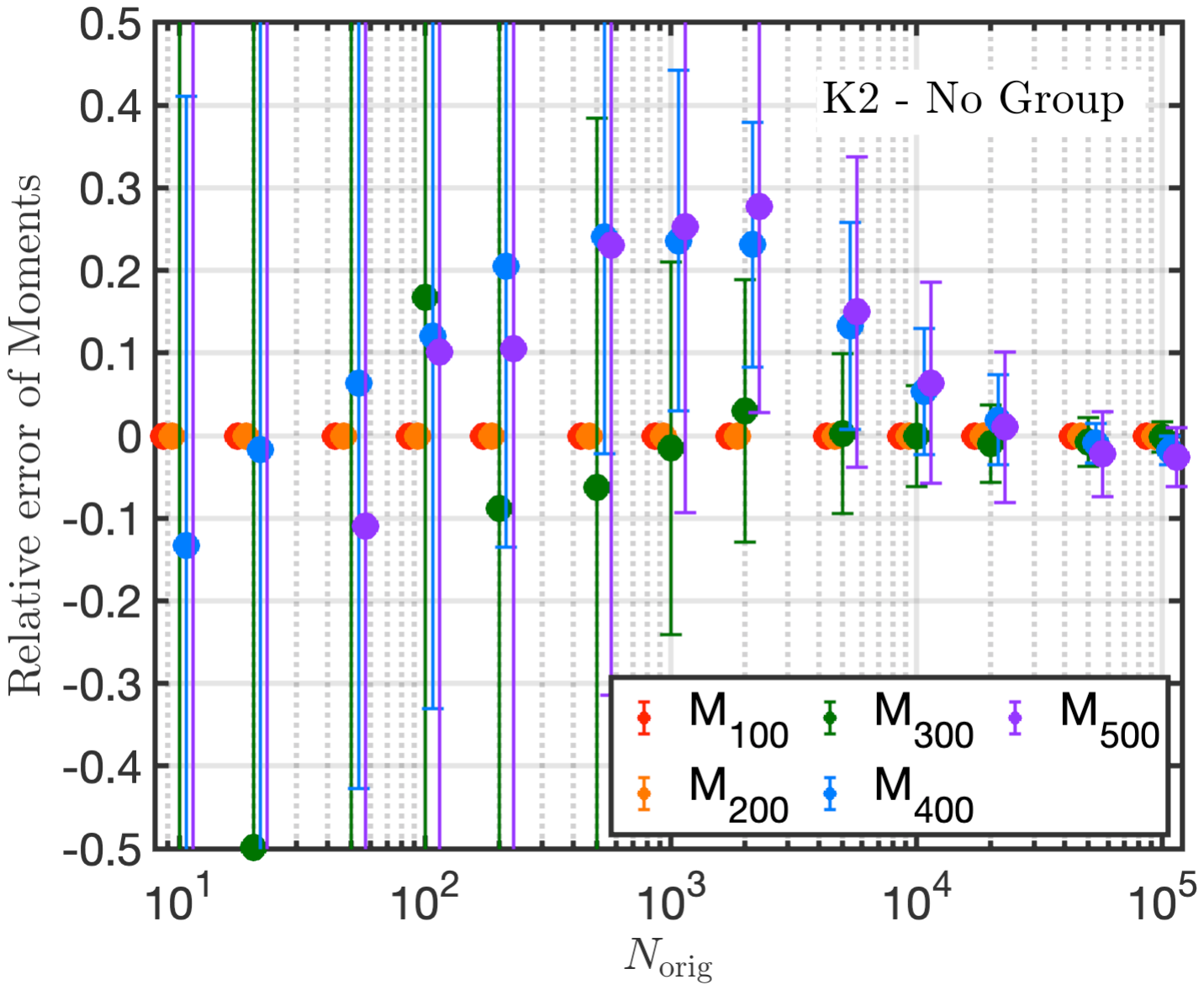}\label{fig:RelMomentErrorNoGroupK2}
        \includegraphics[width=.495\textwidth]{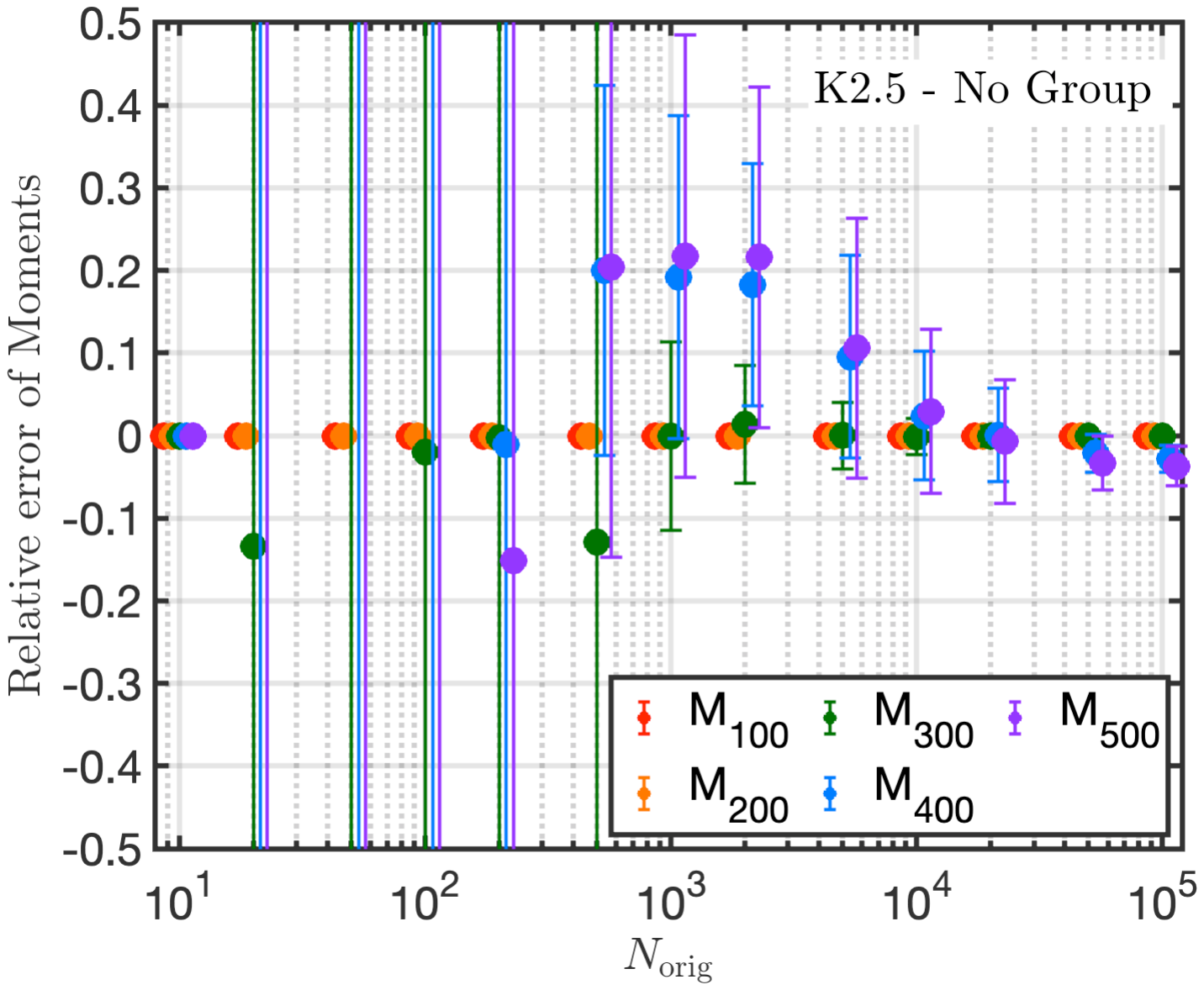}\label{fig:RelMomentErrorNoGroupK2_5}
        \includegraphics[width=.495\textwidth]{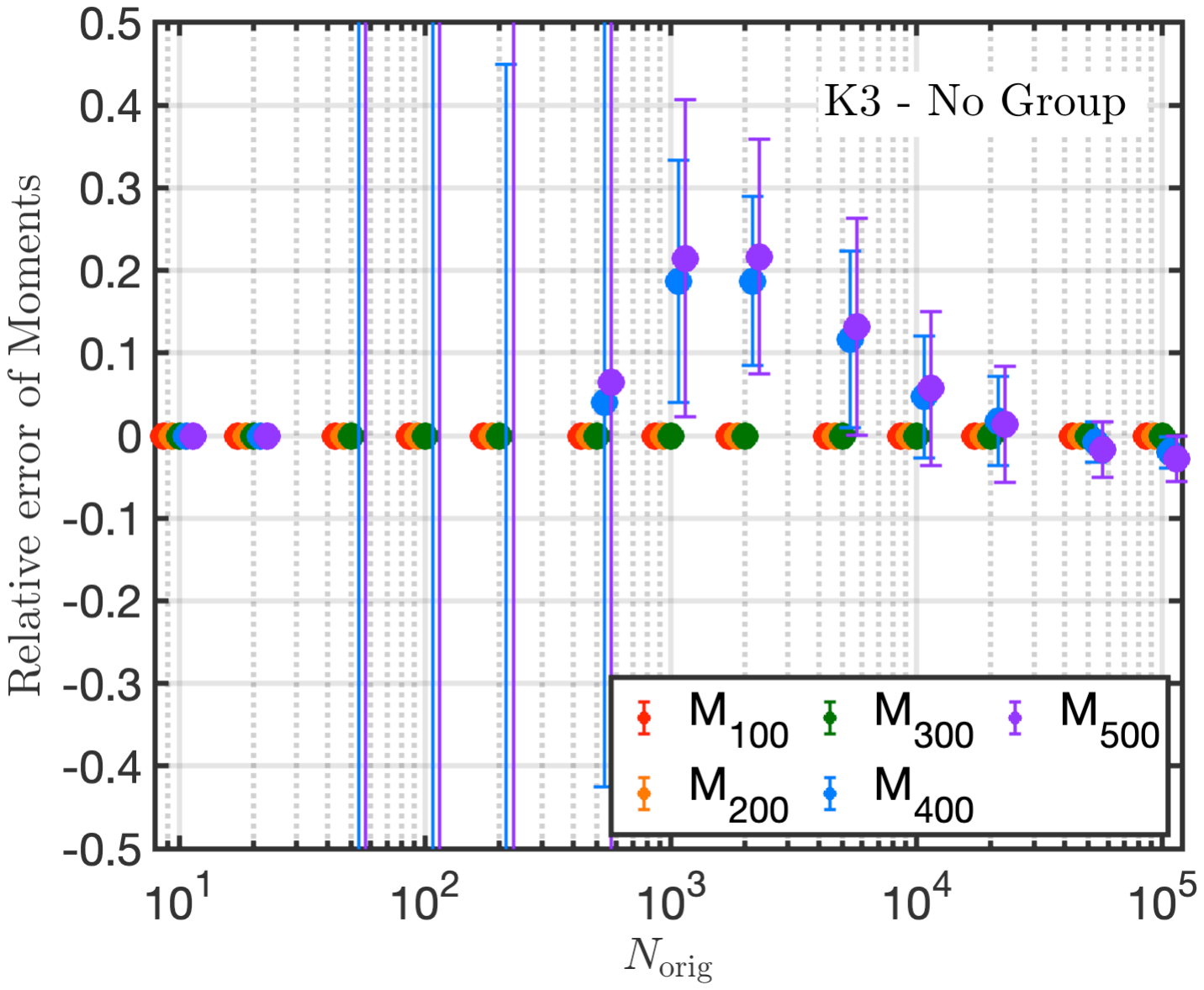}\label{fig:RelMomentErrorNoGroupK3}
\caption{Relative error in the moments for the four reduction schemes without grouping. These results correspond to those for the absolute error shown in Figure~\ref{fig:MomentErrorWithoutGrouping}.}
\label{fig:RelMomentErrorWithoutGrouping}
\end{figure}

\begin{figure}[t]  
        \includegraphics[width=.495\textwidth]{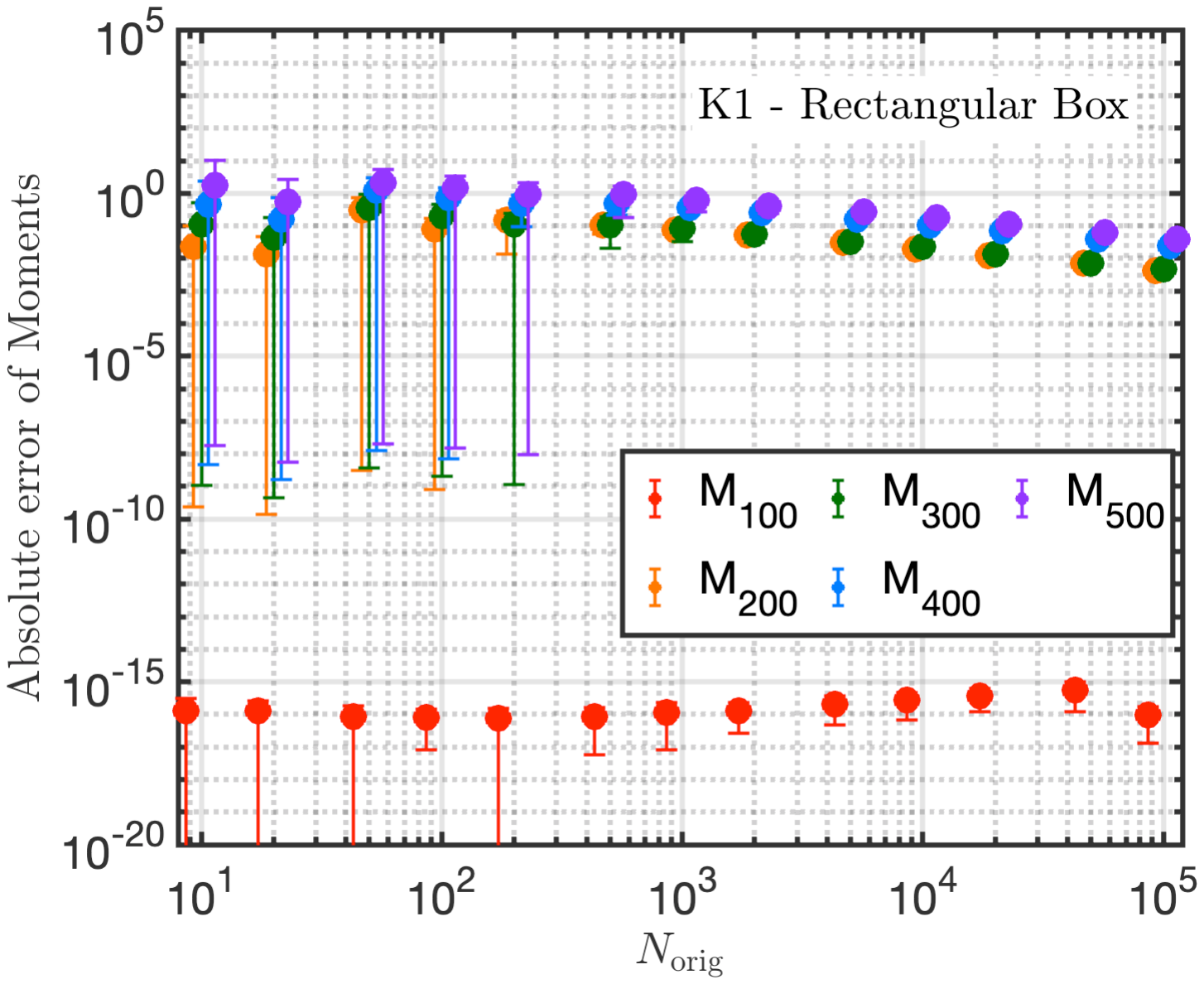}\label{fig:MomentAbsErrorRectK1}
        \includegraphics[width=.495\textwidth]{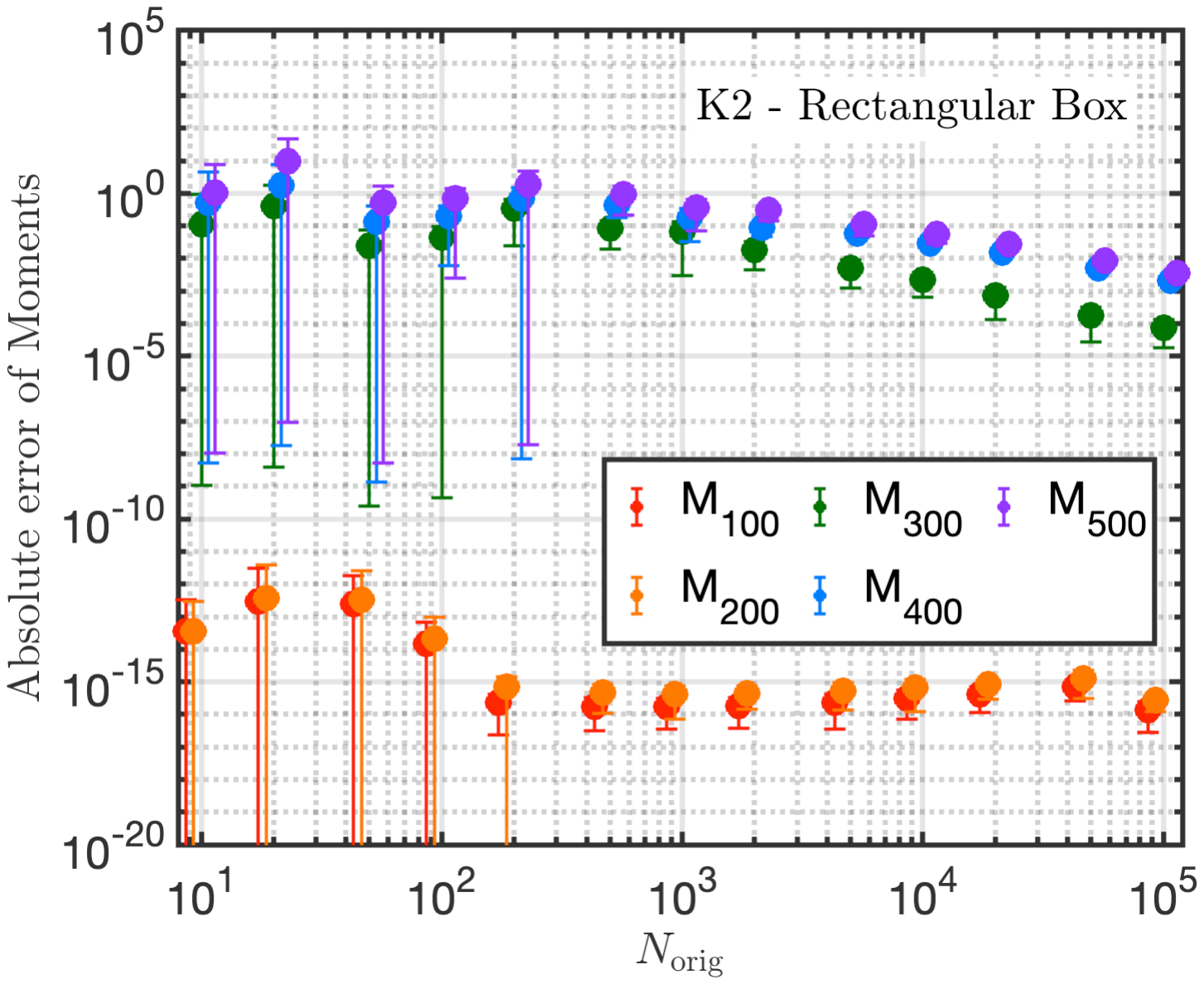}\label{fig:MomentAbsErrorRectK2}
        \includegraphics[width=.495\textwidth]{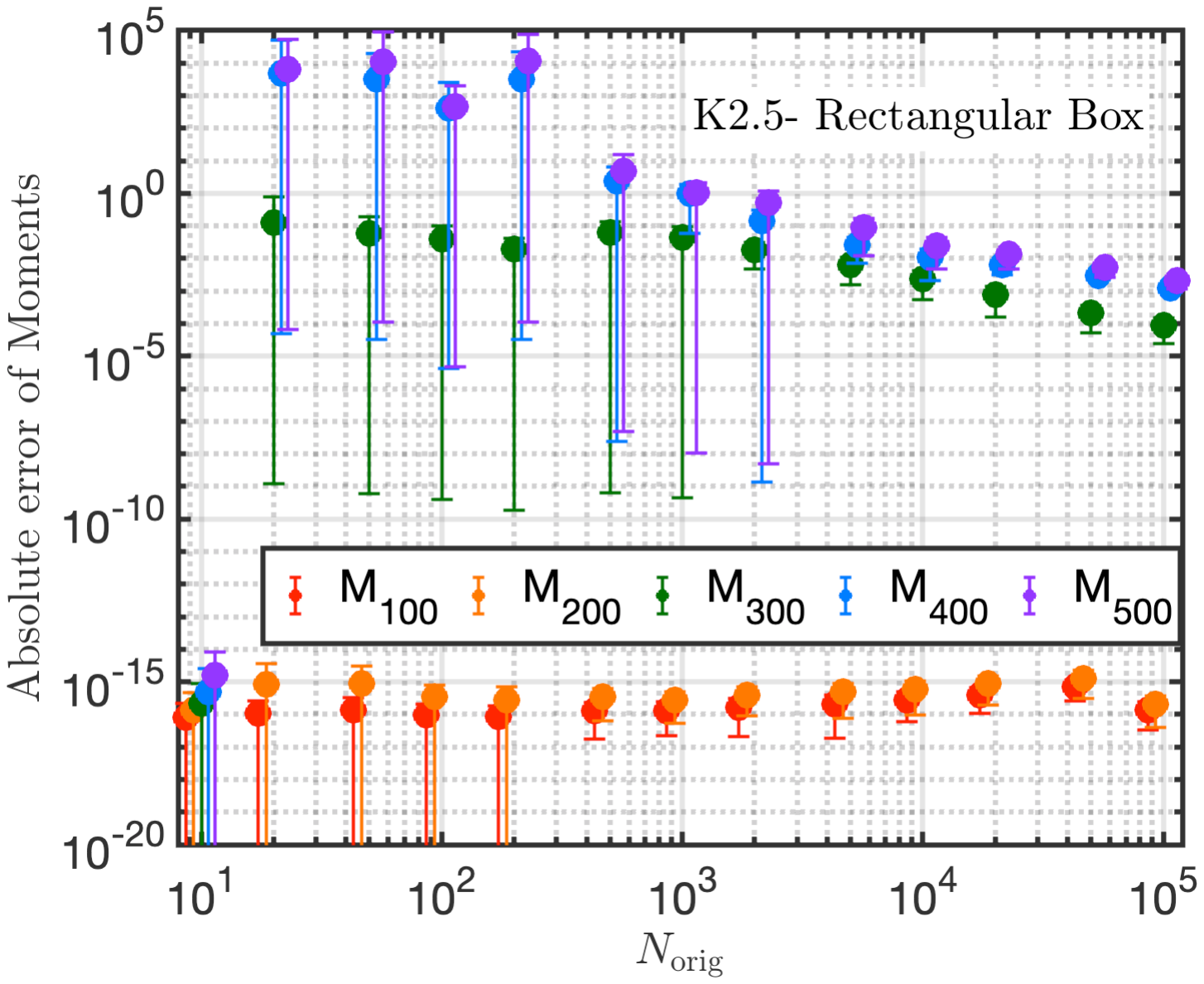}\label{fig:MomentAbsErrorRectK2_5} 
        \includegraphics[width=.495\textwidth]{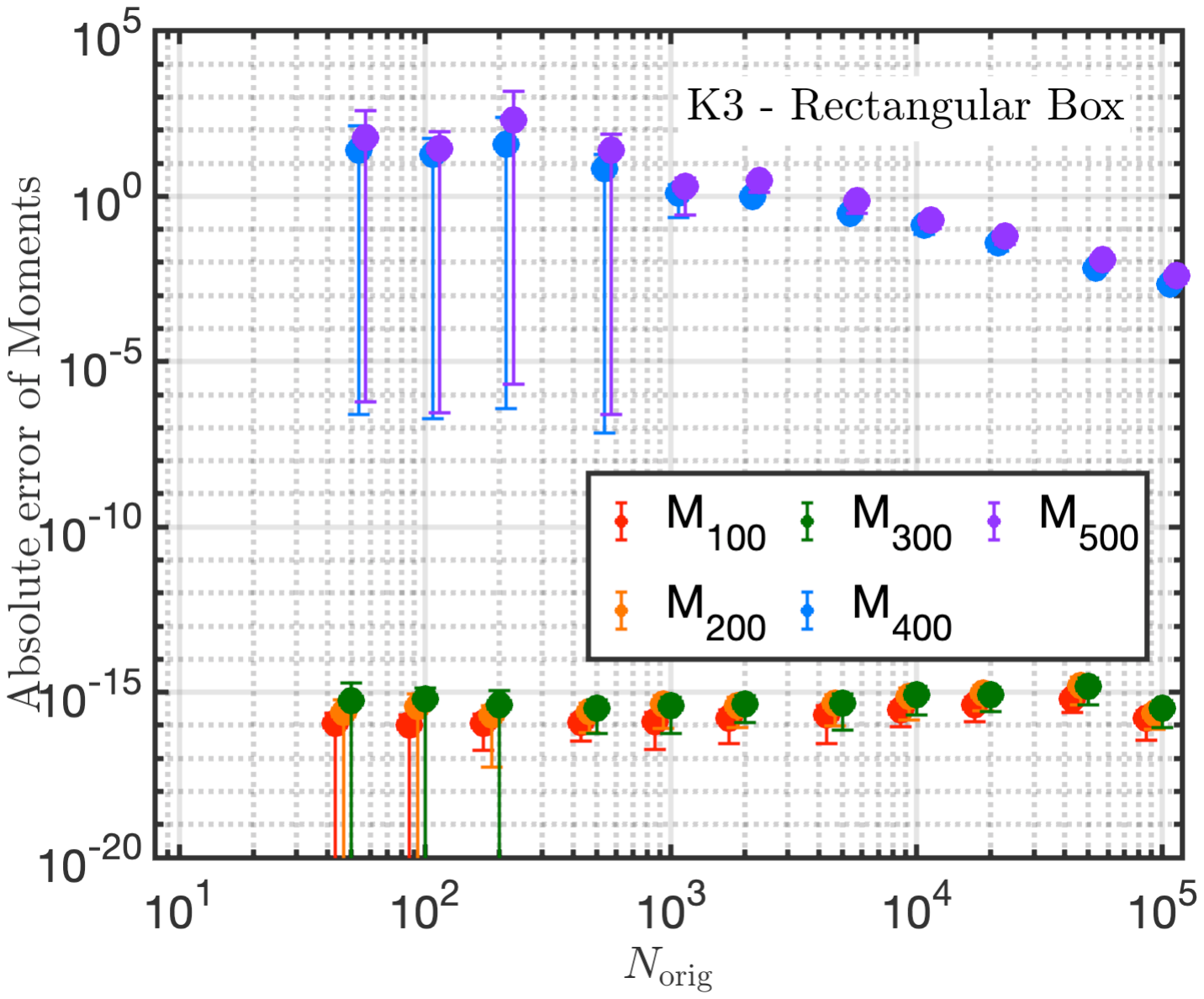}\label{fig:MomentAbsErrorRectK3}
\caption{Absolute error in the moments along the $v_x$ axis due to the four reduction schemes with rectangular box grouping. These results correspond to those in Figure~\ref{fig:MomentErrorWithoutGrouping}.}
\label{fig:MomentErrorWithGrouping}
\end{figure}

\begin{figure}[t]
        \includegraphics[width=.495\textwidth]{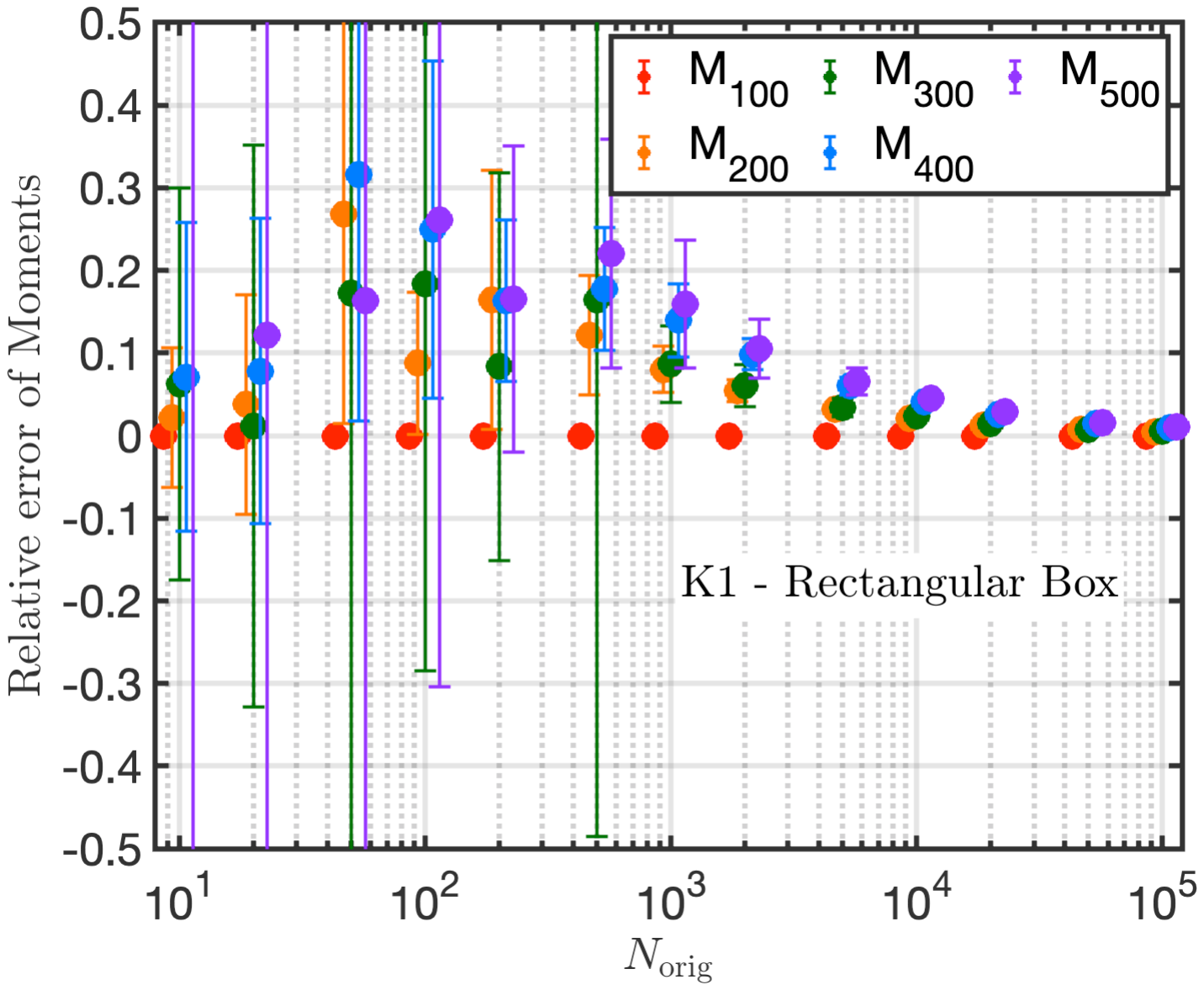}\label{fig:RelMomentErrorK1RectBox}
        \includegraphics[width=.495\textwidth]{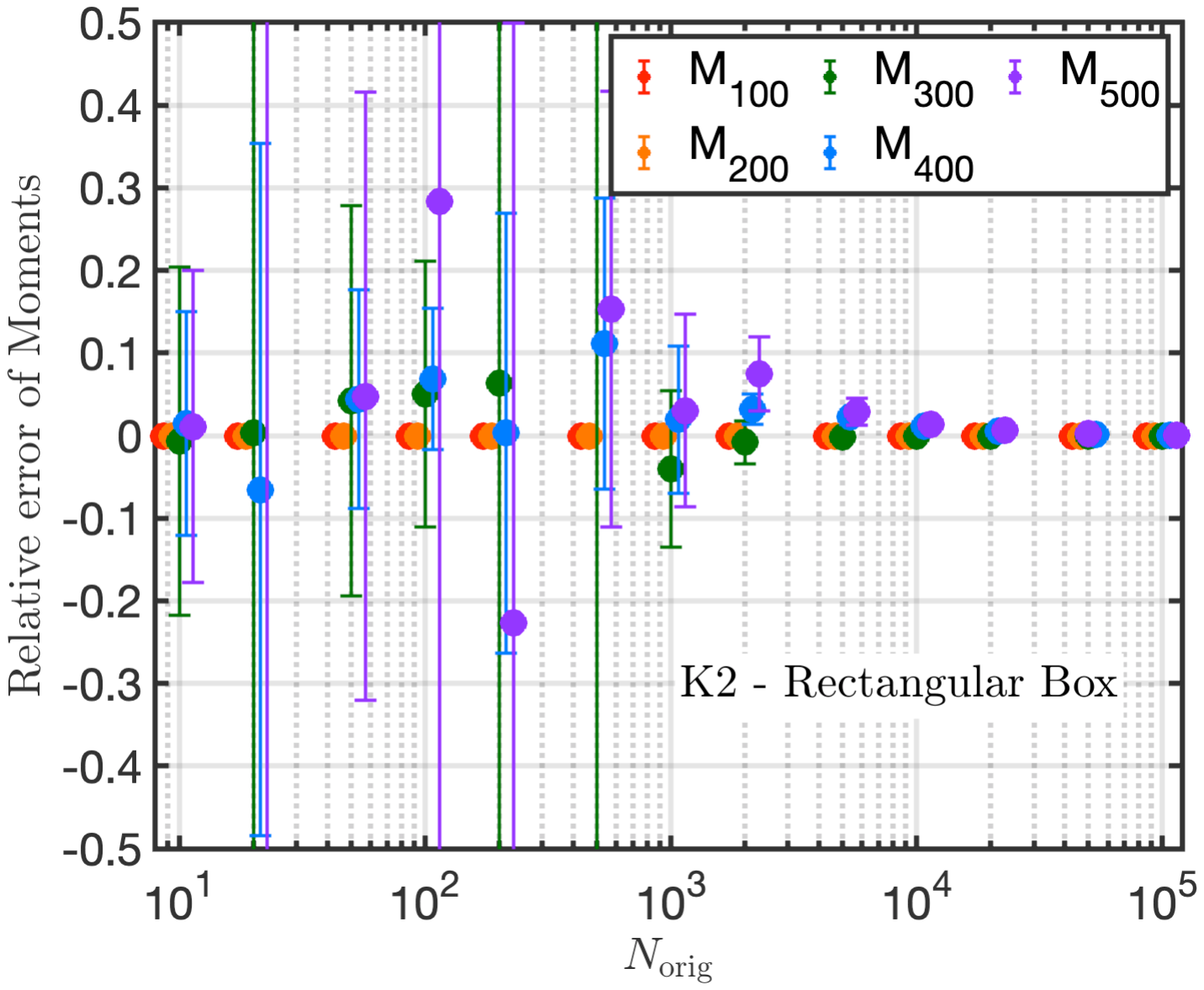}\label{fig:RelMomentErrorK2RectBox}
        \includegraphics[width=.495\textwidth]{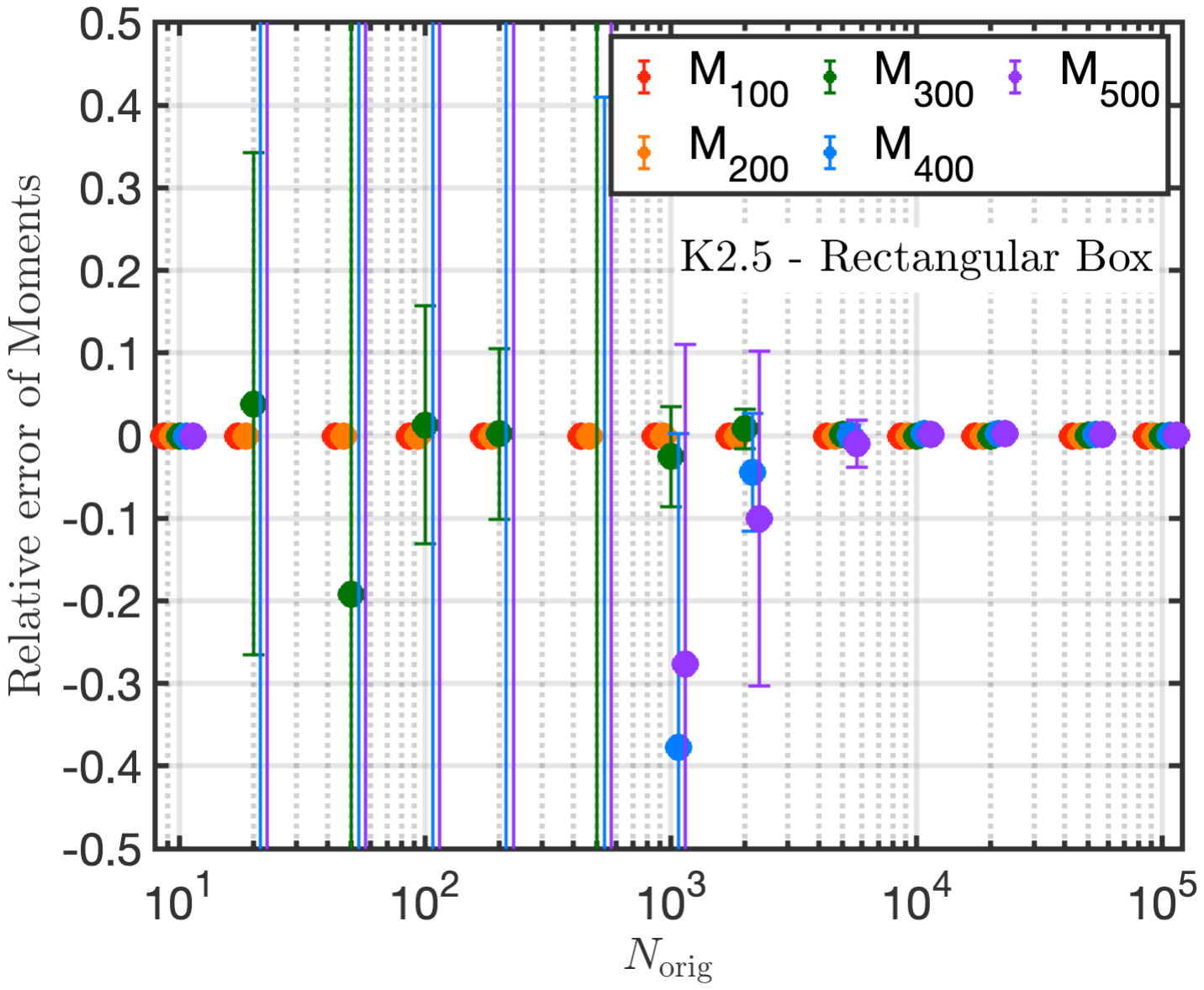}\label{fig:RelMomentErrorK2_5RectBox}
        \includegraphics[width=.495\textwidth]{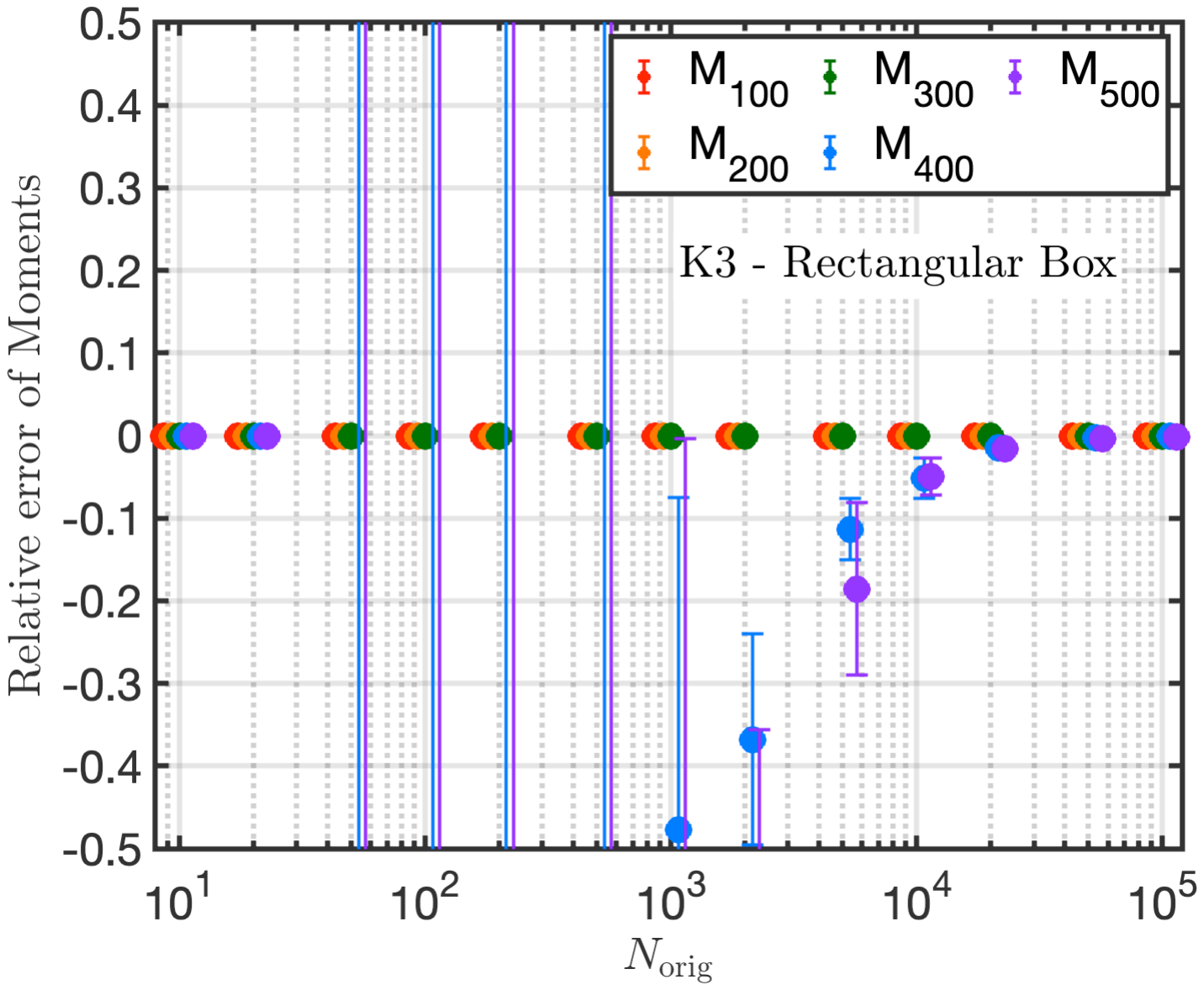}\label{fig:RelMomentErrorK3RectBox}
\caption{Relative error in the moments for the four reduction schemes with rectangular box grouping. These results correspond to those for the absolute error shown in Figure~\ref{fig:MomentErrorWithGrouping}.}
\label{fig:RelMomentErrorWithGrouping}
\end{figure}

\subsection{Uncertainty added by the reduction step}
To examine the uncertainty in the moments due to the reduction step, we first consider several reduction schemes. These reduction schemes can either be applied to the entire set of particles (without grouping), or the particles can first be divided into groups, each of which is then independently reduced (with grouping). We show results for both  reduction without grouping and reduction with grouping. 
For these studies, we set the parameters to be $\boldsymbol{\alpha}=(0.75,0,0)$, $\boldsymbol{\beta}=(0.02,0,0)$, $\delta = \sqrt{2}$, $\gamma = \sqrt{3}$, $l_x = l_y =l_z = 0.5$, and we choose $s$ to be the minimum value necessary for all of the weights to be positive. (Distributions with different levels of skewness, $\boldsymbol{\alpha}$, and kurtosis, $\boldsymbol{\beta}$, show similar results.) 
Specifically, we consider the following cases:
\vskip5pt
\begin{itemize}
    \item[\textbf{K1:}] Here, $K = 1$, which preserves $M_{000}$, $M_{100}$, $M_{010}$ and $M_{001}$. 
    In the standardized coordinate system, this reduction scheme results in a single reduced particle   per group.
    \item[\textbf{K2:}] Here, $K = 2$, which preserves 
     all moments of order less than or equal to two.
 In the standardized coordinate system, this reduction scheme results in $7$ reduced particles per group.
    \item[\textbf{K2.5:}] Here, we preserve all of the moments from case {\bf K2}, together with the third-order moments $M_{300}$, $M_{030}$ and $M_{003}$.  
    In the standardized coordinate system, this reduction scheme results in $10$ reduced particles per group.
    \item[\textbf{K3:}] Here, $K = 3$, which preserves
    all moments of order less than or equal to three. In the standardized coordinate system, this reduction scheme results in $26$ reduced particles per group.
\end{itemize}

\vskip5pt

The first two of these cases have been described elsewhere.  Case {\bf K1} is equivalent to preserving just the drift velocity, as was done by Rjasanow and Wagner~\cite{RWBook2005}.  Case {\bf K2} is equivalent to the work of Lama, Zweck, and Goeckner~\cite{Lama2020}, in which they preserved the full second order tensor.  Cases {\bf K2.5} and {\bf K3} are new.

\begin{figure}[!ht]
        \includegraphics[width=.495\textwidth]{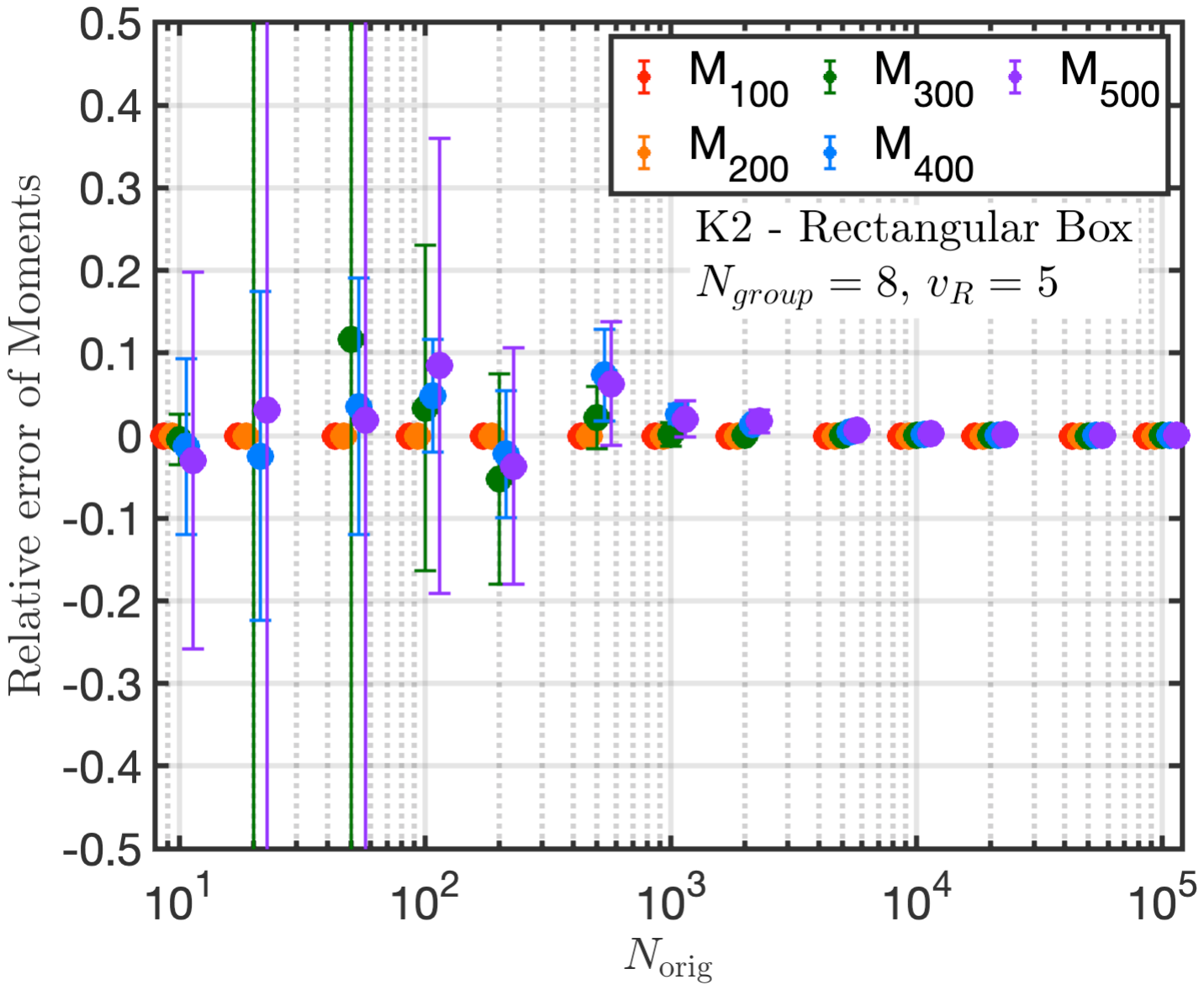}\label{fig:RelMomentErrorK2Mimimized}
        \includegraphics[width=.495\textwidth]{figures/RelErrMK2MinRectBox.pdf}\label{fig:RelMomentErrorK3Mimimized}
\caption{Relative error in moments  for the {\bf K2} (left panel) {\bf K3} (right panel) reduction schemes when $N_{\rm group}$ is set to minimum values of $N_{\rm group} = N_{\rm red} + 1 = 8$ and 27, respectively,  and $v_{R}=5$.}
\label{fig:RelMomentErrorSmallVelVolume}\end{figure}

\subsubsection{Results for reduction without grouping}

Without grouping, we do not consider tail functionals, as the number of stochastic particles after reduction is too small, $N_{\rm red} \le26$, to accurately calculate tail functionals, as is shown in Figure~\ref{fig:InitialTails}.

In Figure~\ref{fig:MomentErrorWithoutGrouping} we show the absolute errors 
\begin{equation*}
    E_{\rm abs} = \left|M_{n00,{\rm red}}-M_{n00,{\rm init}}\right|,
\end{equation*} 
in the moments $M_{100}$, $M_{200}$, $M_{300}$, $M_{400}$ and $M_{500}$ due to the reduction for the four schemes. Results for other moments with the same order as those shown exhibit similar trends. We show the results on a log-log plot to highlight the negligible error in the moments that are explicitly preserved by each of the schemes. In the cases {\bf K1}, {\bf K2} and {\bf K3}, for which the full $K$-th order tensor moments are preserved, the $K$-th order scalar moment is preserved to within numerical error.  Comparing case {\bf K2} to case {\bf K2.5}, for which $M_{300}$, $M_{030}$ and $M_{003}$ are also preserved in the rotated and standardized coordinate system, we find an improvement in $M_{300}$ in the original coordinate system, but not to the same level as when the full moment tensor is preserved in case {\bf K3}. From this we conclude that while it is helpful to preserve individual scalar moments \eqref{eqn:MomentOfTriple}, where feasible, it is best to preserve the full third-order tensor moment. 

In Figure~\ref{fig:RelMomentErrorWithoutGrouping}  we show the corresponding results for the relative errors
\begin{equation*}
    E_{\rm rel} = \frac{M_{n00,{\rm red}}-M_{n00,{\rm init}}}{\left|M_{n00,{\rm init}}\right|}.
\end{equation*}
We show the results on a semilog plot to highlight the error in the non-preserved moments. In case {\bf K1}, the non-preserved moments all show significant error. This is because the higher order moments are under-determined with a single replacement particle.  In cases {\bf K2}-{\bf K3} as $N_{\rm orig}$ is increased, the errors in the non-preserved moments decrease, which is to be expected from the convergence theorem of Rjansanow and Wagner~\cite{RWBook2005}.
We also found that  the accuracy of non-preserved moments can significantly depend on the other simulation parameters.  A  numerical study showed that setting $l$ near 1, for example $l=1.0001$, causes the error of the non-preserved moments to grow by many orders of magnitude. The same study showed that the error of the non-preserved moments was minimized when $l$ was near 2 or 0.5.   Furthermore, if we set $s$ to be much larger than the minimum value needed to ensure weight positivity, then the error in the non-preserved moments also grew by many orders of magnitude.  For these reasons, for the results in the rest of this section we used $l = 0.5$ and values of
$s$ that are very close to the minimum required for weight positivity.  

\subsection{Results for reduction with grouping}\label{sec:grouping_reduction}

In this section, we study the errors for reduction schemes with grouping. Both the grouping technique and the reduction scheme impact the errors in the moments and tail functionals.

To investigate this, here we present preliminary results obtained using a rectangular box grouping technique. (We plan to explore the impact of other grouping techniques in a future article.) With this technique, we create
    \begin{equation}
        n_{\rm groups}=\frac{6}{\pi}\frac{N_{\rm orig}}{N_{\rm group}}
        \label{eq:NumberofGroups}
    \end{equation}
    of rectangular boxes of equal volume, $\Delta V = \Delta v_x \Delta v_y\Delta v_z$,
where $N_{\rm group}$ is the desired number of particles in each group and $\Delta v_x$, $\Delta v_y$ and $\Delta v_z$ are the three side lengths of the boxes. The number of boxes in each direction is chosen so that the union of the boxes covers a cube of volume $V_{\rm cube}=(2 v_{R})^3$, which in turn covers the entire ball of volume $V_{\rm ball}=\frac{4\pi}{3}v_{R}^3$ containing all the stochastic particles. The factor of $\frac{6}{\pi}$ in \eqref{eq:NumberofGroups}  accounts for the volume inside the box that is outside the ball.  Unless otherwise noted, we used $v_{R}=7$. To determine the side lengths of the boxes, we first randomly chose two components of $\bf{v}$ and set  
    \begin{equation*}
        \Delta v_{i,j} = {\rm floor}\left(n_{\rm groups}^{1/3}\right)
        \qquad\text{and}\qquad\Delta v_{k} = {\rm floor}\left(
        \frac{n_{\rm groups}}{\Delta v_{i}\Delta v_{j}}\right),
    \end{equation*}
thereby maximizing the number of groups. With these choices, the boxes that are entirely contained in the ball of  radius $v_{R}$ each  have close to $N_{\rm group}$ particles. 
For {\bf K2}, $N_{\rm group} = 11$ and $N_{\rm red}=7$,
for {\bf K2.5}, $N_{\rm group} = 15$ and $N_{\rm red}=10$, and 
for {\bf K3}, $N_{\rm group} = 39$ and $N_{\rm red}=26$.
In this way, the total number of stochastic particles after the reduction is $N_{\rm red}\approx \frac{2}{3}N_{\rm orig}$.  For {\bf K1}, we chose $N_{\rm group} = 2$ resulting in $N_{\rm red}=1$, as this is the minimum reduction that can be performed with the {\bf K1} scheme. Finally, groups that had $N_{\rm group} \le N_{\rm red}$ were not reduced, and those original particles were added back into the set of all reduced particles.  This is particularly important for boxes that lie partially outside the ball of radius $v_{R}$, and hence statistically may contain few or no particles.


In Figure~\ref{fig:MomentErrorWithGrouping} we show the absolute errors in the moments due to the four reduction schemes  with rectangular box group.  As before, the preserved moments are correct to within numerical error.  We also observe that the absolute errors in the non-preserved moments are smaller for larger $N_{\rm orig}$, and hence larger $n_{\rm groups}$. This result is consistent with the convergence theorem of Rjasanow and Wagner~\cite{RWBook2005}.


In Figure~\ref{fig:RelMomentErrorWithGrouping} we show the corresponding results for the relative error. Just as for the absolute error, there is essentially no error in the preserved moments. For the {\bf K1} reduction scheme with grouping, there is a significant improvement in the non-preserved moments compared to the scheme without grouping. This is simply because the total number of stochastic particles after reduction is much larger with grouping than without. Moreover, the relative error in the higher-order moments due to the grouping and reduction processes is quite large for $N_{\rm orig}\lesssim 1000$. 

As more moments are explicitly preserved, from {\bf K2} up to {\bf K3}, the error in the higher-order moments, $M_{400}$ and $M_{500}$, increases when $N_{\rm orig}$ is between $1000$ and $10,000$. Although this trend may seem counter intuitive,  the convergence theorem of Rjasanow and Wagner~\cite{RWBook2005} provides a clue as to why it occurs. Under mild assumptions, this theorem states that the empirical measure computed using an SWPM simulation converges to the true solution as the velocity radius, $R_{\rm group}$,  of the group decreases. Therefore, the error in the non-preserved moments is expected to increase as $R_{\rm group}$ increases. This is exactly what happens as we go from the {\bf K2} up to the {\bf K3} reduction scheme. The reason is that the value of $R_{\rm group}$ increases as the volume of the simulation ball  increases and as the number of groups, $n_{\rm groups}$, decreases. For the results in Figure~\ref{fig:RelMomentErrorWithGrouping}, the volume of the simulation ball is fixed (as $v_R=7$) and $n_{\rm groups}$ decreases because $N_{\rm group}$ increases (see \eqref{eq:NumberofGroups}).

To improve the accuracy of the non-preserved moments and tail functionals we aim to minimize $R_{\rm group}$. The value of $R_{\rm group}$ can be decreased by choosing the minimum possible value for $N_{\rm group}$, namely $N_{\rm group} = N_{\rm red}+1$, and by decreasing the value of $v_R$. Setting $N_{\rm group}$ at the minimum values of, 8, 11 and 27, decreases $R_{\rm group}$ by roughly 8, 10 and 12\% respectively for {\bf K2}, {\bf K2.5} and {\bf K3}. Moreover, decreasing $v_R$ from $7$ to $5$ reduces $R_{\rm group}$ by an additional 29\%. However, with $v_R=5$ we can no longer calculate $\operatorname{Tail}(5)$ and $\operatorname{Tail}(6)$. That is, there is a trade-off between accurate computation of higher-order moments and the ability to compute the lower-probability tails of the distribution. This trend can be seen in Figure~\ref{fig:RelMomentErrorSmallVelVolume} where we show that the relative error under these conditions is negligible for {\bf K2} once $N_{\rm orig}\ge1000$ and for {\bf K3} once $N_{\rm orig}\ge5000$.

\begin{figure}[!ht]
\begin{center}
    \includegraphics[width=.495\textwidth]{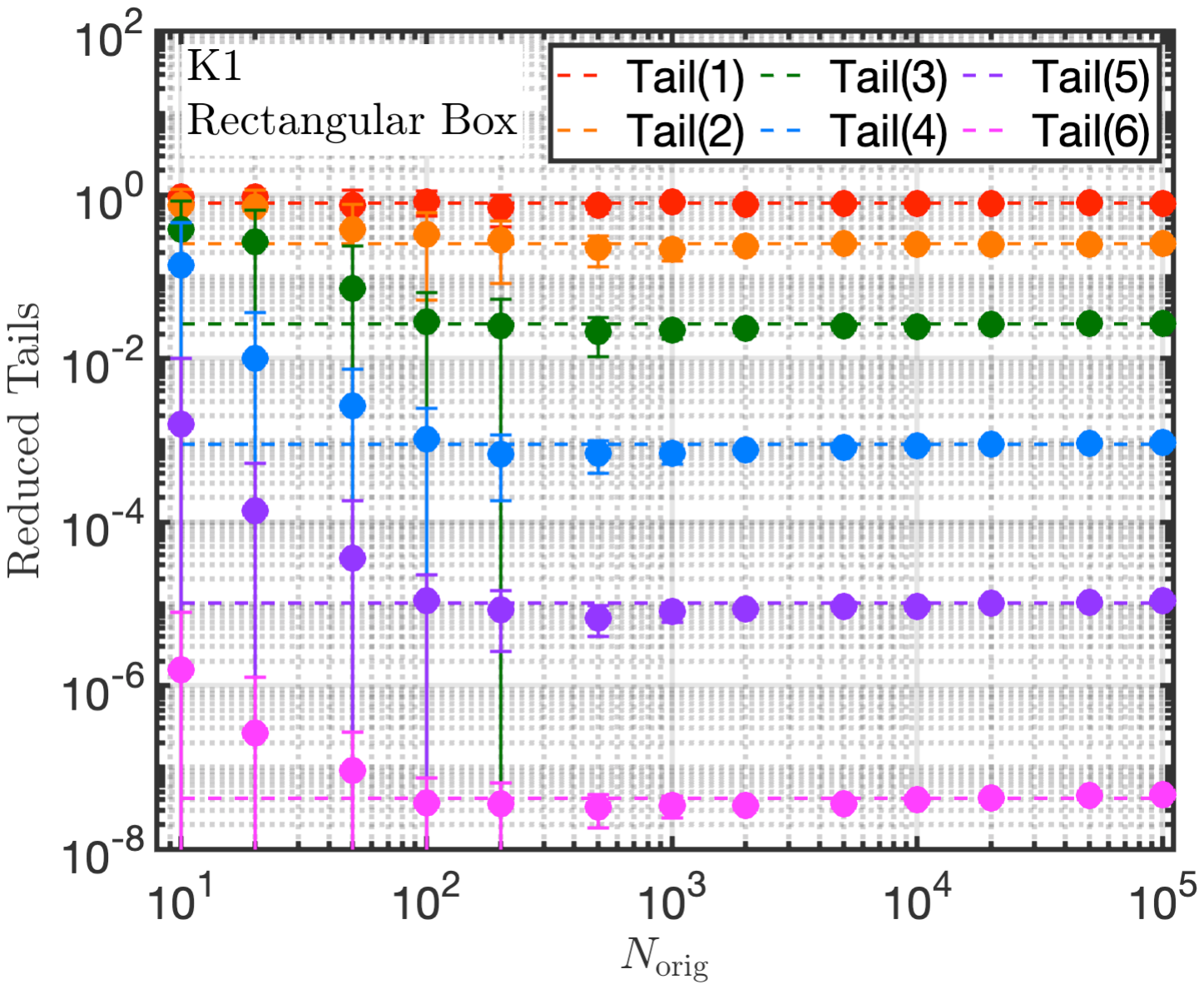}\label{fig:reducedTailK1_RectBox}
    \includegraphics[width=.495\textwidth]{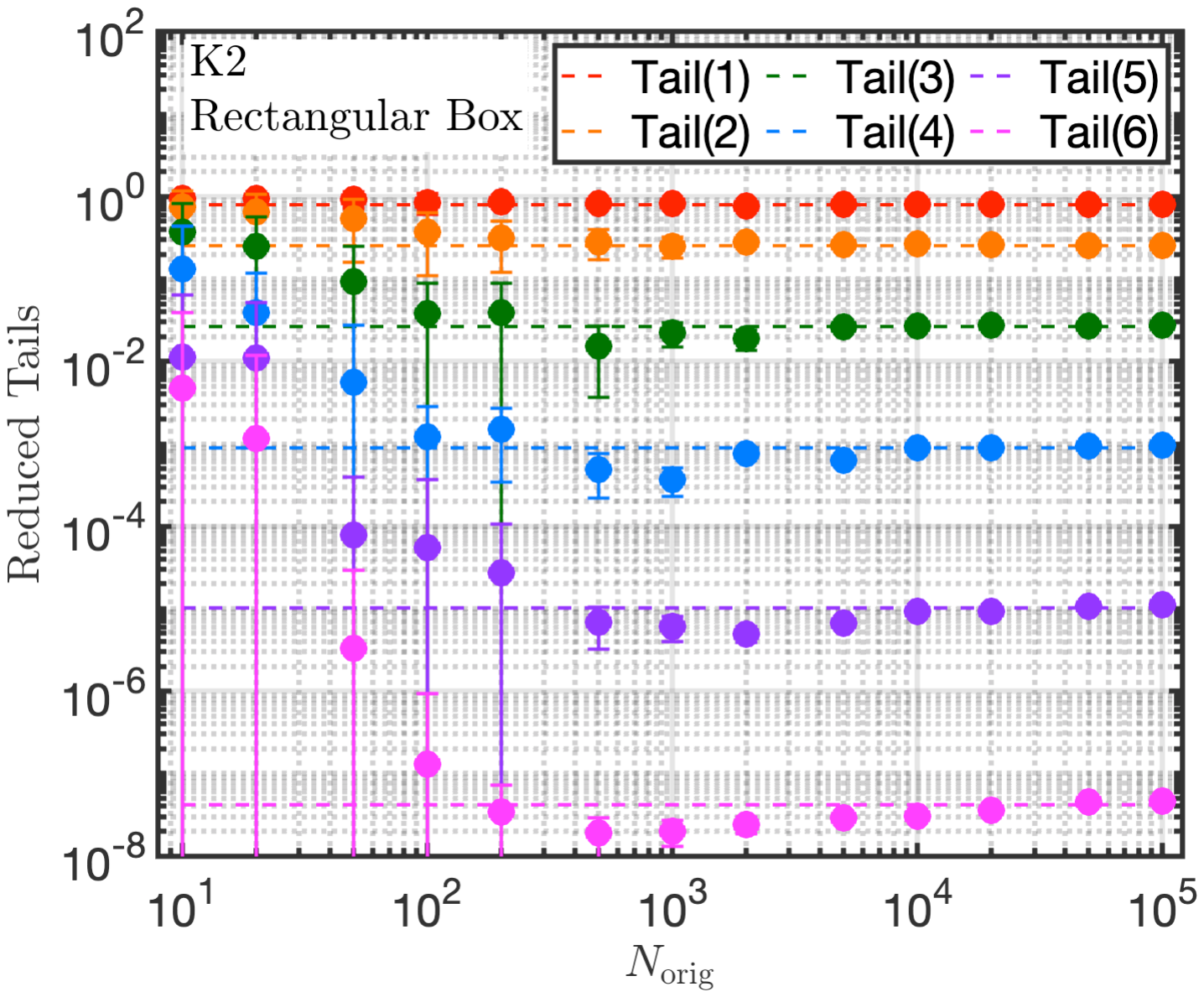}\label{fig:reducedTailK2_RectBox}
    \includegraphics[width=.495\textwidth]{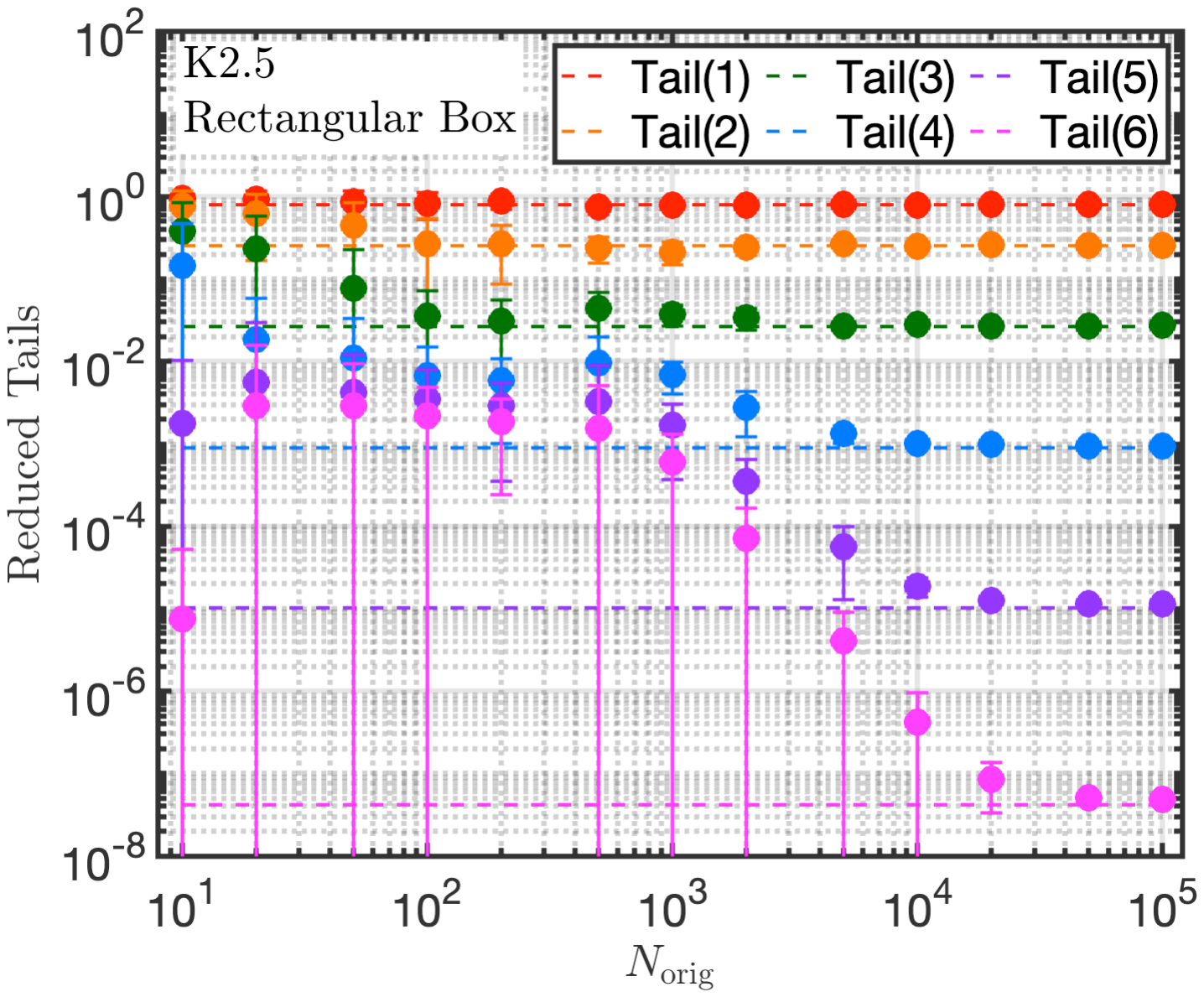}\label{fig:reducedTailK2_5_RectBox}
    \includegraphics[width=.495\textwidth]{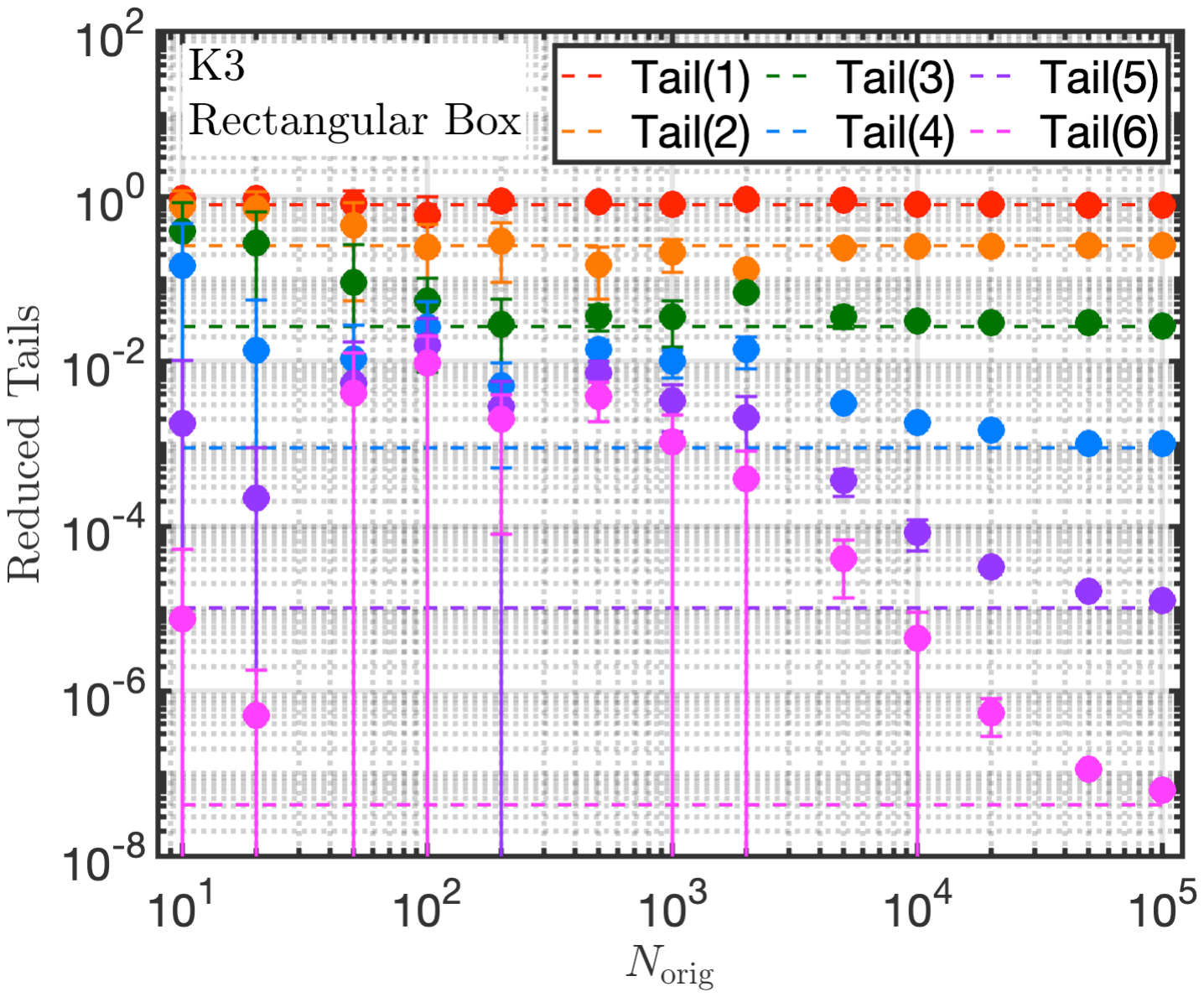}\label{fig:reducedTailK3_RectBox}
\end{center}
\caption{Tail functionals for the four reduction schemes {\bf K1} (top left), {\bf K2} (top right), {\bf K2.5} (bottom left), and {\bf K3} (bottom right) with  rectangular box grouping.}
\label{fig:FinalTailsGrouped}\end{figure}


In Figure~\ref{fig:FinalTailsGrouped}, we show tail functionals after reduction with the rectangular box grouping technique. The results in the top left panel for the {\bf K1} reduction scheme are almost identical to the corresponding results prior to reduction shown in the bottom right panel of Figure~\ref{fig:InitialTails}, with $\operatorname{Tail}(6)$ being preserved with $N_{\rm orig}\ge100$. Therefore if the goal is simply to preserve tail functionals the strategy should be to initialize SWPM using a  uniform distribution of  particles over a  ball whose radius exceeds that of the tail function to be computed, to employ a  grouping technique that (optimally) satisfies the assumptions of the convergence theorem of Rjasanow and Wagner~\cite{Lama2020,RWBook2005}, and, in the standardized coordinate system, to replace each group with a  single particle at the drift velocity whose weight equals that of the group. Additional work is required to identify the optimal grouping technique to best preserve $\operatorname{Tail}(N)$. If instead the goal is to preserve  tail functionals and moments of order three and less, then it is advantageous to use one of the higher-order reduction schemes. For example, with the {\bf K2}, {\bf K2.5}, and {\bf K3}  reduction schemes, $\operatorname{Tail}(6)$ can be preserved with $N_{\rm orig}\ge5,000$, $N_{\rm orig}\ge50,000$, and $N_{\rm orig}\ge100,000$, respectively. Significantly, SWPM with any of these reduction schemes is much more efficient than DSMC, which requires $N_{\rm orig}\gg100,000$  to accurately and precisely compute $\operatorname{Tail}(6)$.

Finally and importantly, because there are fewer groups to populate and reduce, schemes that preserve more higher-order moments generally reduce the computational cost. For example, the cost of running a minimal {\bf K3} reduction is 9 times less than that of a minimal {\bf K2} reduction.

\section{Summary}
\label{sec:Summary}
We developed a general particle reduction scheme for the Stochastic Weighted Particle Method of Rjasanow and Wagner.  The method involves determining the velocities and weights of a reduced system of particles so as to preserve a specified set of moments of an original larger particle system.  The
moment preservation condition results in a linear system 
for the weights in terms of the velocities and moments. 
For reduction schemes that preserve a given collection of moments of order at most three, we show how to choose the velocities using as few particles as possible, so as to ensure that there is an explicit analytical solution with positive weights.

To test the accuracy of these new reduction schemes, we generated
original collections of stochastic particles that would be typically encountered in a SWPM simulation and compared explicitly preserved moments, non-preserved moments, and tail functionals before and after reduction.
As expected, the errors between the preserved moments before
and after reduction are negligible. Moreover, if the original set of particles is divided into groups and each group is reduced
so as to preserve all moments up to order at \emph{least} two, then the 
error in the non-preserved 4-th and 5-th order moments is also
negligible provided that several conditions are met. These conditions are that 
the total number of stochastic particles prior to reduction is large enough, the velocities are all chosen to lie within a small enough ball, and the parameters in the formulae for the velocities of the reduced particles are chosen appropriately.  However, in order to preserve tail functionals out to a given speed, $R$, it is more efficient to preserve all moments up to order at \emph{most} two with a scheme in which all particles lie in a ball of radius at least $v_R=R+1$. Nevertheless, we note the importance of explicitly preserving the third order moment tensor for systems with skewness, such as occur in hypersonic flows and laser ablation into a plasma.

In future publications we will apply these methods to true SWPM simulations that incorporate the physics, more fully investigate the effects that the choice of grouping technique has on the results, and
extend the generalized reduction scheme to spatially dependent
velocity distributions.

\appendix
\section{$K=3$ expanded progenitor matrix solution parameters} 
\label{subsec::parm}
In this appendix we  provide explicit formulae for the parameters in the solutions, $\mathbf{w}$, for the expanded progenitor matrix in subsection \ref{sec:PosWeightsK3}. 

For the block weight vector $\mathbf{w}_x$ found in \eqref{wexdot} the parameters are given by
\begin{equation}
    a_x = 1-|M_{110}| - |M_{101}| ,
\end{equation}
\begin{equation}
    b_x = \delta(|M_{210}|+|M_{120}|+|M_{201}|+|M_{102}|)-\gamma|M_{111}|,
\end{equation}
\begin{equation}
    c_{x\pm} = \beta_x\frac{(M_{120}+M_{102})(\delta^2l_x\mp1)\pm M_{300}}
    {l_x\mp1}.
\end{equation}
Analogous formulae hold for the block weight vectors, $\mathbf w_y$ and
$\mathbf w_z$ in \eqref{wexsoldot++}.  
Finally, for $w_0$ found in \eqref{eqn:w0_formula} the parameters are given by
\begin{equation}
    a_0 =  -3+(2-\delta^2)(|M_{011}|+|M_{101}|+|M_{110}|),
\end{equation}
and
\begin{align}
\begin{split}
    b_0 = &\frac{\beta_x}{l_x}M_{300}+\frac{\beta_y}{l_y}M_{030}+\frac{\beta_z}{l_z}M_{003} +\gamma(3-\gamma^3)|M_{111}|\\
    &+(\delta^2-1)(\frac{\beta_x}{l_x}(M_{102}+M_{120})+\frac{\beta_y}{l_y}(M_{012}+M_{210})+\frac{\beta_z}{l_z}(M_{021}+M_{201}))\\
    &+\delta(2-\delta^2)(|M_{012}|+|M_{021}|+|M_{102}|+|M_{120}|+|M_{201}|+|M_{210}|).
\end{split}
\end{align}


\bibliographystyle{siamplain}
\bibliography{references}

\end{document}